\documentclass[twoside]{amsart}
\usepackage{float}
\usepackage[bookmarksnumbered,plainpages,hypertex]{hyperref}
\usepackage{graphicx}
\usepackage{amssymb}
 \input{xy}
 \xyoption{all}
\newtheorem{theorem}{\sc Theorem}[section]
\newtheorem{proposition}[theorem]{\sc Proposition}
\newtheorem{lemma}[theorem]{\sc Lemma}
\newtheorem{corollary}[theorem]{\sc Corollary}
\theoremstyle{definition}
\newtheorem{definition}[theorem]{\sc Definition}
\newtheorem{definitions}[theorem]{\sc Definitions}
\newtheorem{example}[theorem]{\sc Example}
\newtheorem{examples}[theorem]{\sc Examples}
\theoremstyle{remark}
\newtheorem{remark}[theorem]{\sc Remark}
\newtheorem{remarks}[theorem]{\sc Remarks}
\newtheorem{claim}[theorem]{\sc}

\newtheorem{nonamethm}[theorem]{\sc}
\def\K{\textrm{ker}}
\def\C{\textrm{Coker}}
\def\Id{\textrm{Id}}
\def\M{\mathcal{M}}

\def\H{\mathbb{H}}
\def\T{\mathbb{T}}

\def\ot{\otimes}
\def\cot{\square}
\def\N{\mathbb{N}}
\def\Z{\mathbb{Z}}
\def\Im{\textrm{Im}}
\begin{document}

\pagestyle{headings}
\title{Separable Functors and Formal Smoothness}
\author{Alessandro Ardizzoni}
\address{University of Ferrara, Department of Mathematics, Via Machiavelli 35,
Ferrara, I-44100, Italy} \email{alessandro.ardizzoni@unife.it}
\urladdr{http://www.unife.it/utenti/alessandro.ardizzoni}
\thanks{This paper was written while the author was member of G.N.S.A.G.A. with partial financial support
from M.I.U.R..}

\begin{abstract}
The natural problem we approach in the present paper is to show
how the notion of formally smooth (co)algebra inside monoidal
categories can substitute that of (co)separable (co)algebra in the
study of splitting bialgebra homomorphisms. This is performed
investigating the relation between formal smoothness and
separability of certain functors and led to other results related
to Hopf algebra theory. Between them we prove that the existence
of $ad$-(co)invariant integrals for a Hopf algebra $H$ is
equivalent to the separability of some forgetful functors. In the
finite dimensional case, this is also equivalent to the
separability of the Drinfeld Double $D(H)$ over $H$. Hopf algebras
which are formally smooth as (co)algebras are characterized. We
prove that given a bialgebra surjection $\pi :E\rightarrow H$ with
nilpotent kernel such that $H$ is a Hopf algebra which is formally
smooth as a $K$-algebra, then $\pi $ has a section which is a
right $H$-colinear algebra homomorphism. Moreover, if $H$ is also
endowed with an $ad$-invariant integral, then this section can be
chosen to be $H$-bicolinear. We also deal with the dual case.
\end{abstract}
 \maketitle \tableofcontents
\section*{Introduction}
Separable functors were introduced by C. N\u{a}st\u{a}sescu, M.
Van den Bergh and F. Van Oystaeyen in \cite{NdO}. As highlighted
in \cite{CMZ}, the relevance of these functors lies in a
functorial version of Maschke's theorem they satisfy, namely they
reflect split exact sequences. In \cite[Corollary
2.31]{AMS-Spliting}, this property was applied to the following
situation. Let $H$ be a semisimple and cosemisimple Hopf algebra
over a field $K$ and denote by $^H\mathfrak{M}^H$ the category of
$H$-bicomodules. Then the forgetful functor
$U:{_A^H\mathfrak{M}_A^H}\to {_A\mathfrak{M}_A}$, from the
category of $A$-bimodules in $^H\mathfrak{M}^H$ to the category of
ordinary $A$-bimodules, is a separable functor and hence the
multiplication of $A$ splits as a morphism of $A$-bimodules and
$H$-bicomodules (i.e. $A$ is separable as an algebra in the
monoidal category $(^H\mathfrak{M}^H,\ot_K,K)$) if and only if it
splits as a morphism of $A$-bimodules (i.e. $A$ is separable  as
an ordinary $K$-algebra). The proof of separability of the functor
$U$ relies on the existence of an $ad$-invariant integral
(introduced by D. \c{S}tefan and F. Van Oystaeyen in
\cite[Definition 1.11]{DO}) for any semisimple and cosemisimple
Hopf algebra over a field $K$. The characterization of separable
and formally smooth algebras in the framework of monoidal
categories was developed in \cite{AMS}. The notion of formal
smoothness (or quasi-freeness) for algebras over a field $K$ was
introduced by J. Cuntz and D. Quillen in \cite{CQ} to provide a
natural setting for non-commutative version of certain aspects of
manifolds. A formally smooth algebra $A$ in monoidal categories
behaves like a free algebra with respect to nilpotent extensions
in the sense that, under natural conditions, any algebra
homomorphism $A\to R/I$, where $I$ is a nilpotent ideal of an
algebra $R$, can be lifted to an algebra homomorphism $A\to R$.
This gives a natural way to produce algebra sections in
$\mathcal{M}$ for algebra homomorphisms $E\to A$ which are
epimorphisms with nilpotent kernel in $\mathcal{M}$. Like in the
classical case any separable algebra in a monoidal category is in
particular formally smooth. As a consequence, in
\cite{AMS-Spliting} it was shown that if $E$ is a bialgebra such
that $H=E/J$ is a quotient Hopf algebra of $E$ which is
semisimple, $J$ denoting the Jacobson radical of $E$, then the
canonical Hopf projection $\pi:E\to H$ admits a left $H$-colinear
algebra section $\sigma:H\to E.$ Furthermore this section can be
chosen to be $H$-bicolinear, whenever $H$ is also cosemisimple. In
\cite{AMS-Spliting} also the dual situation of a bialgebra $E$
whose coradical, say $H$, is a Hopf subalgebra is described. In
this case there is a retraction $\pi$ of the canonical injection
$\sigma$ which is a left $H$-linear (bilinear if $H$ is also
semisimple) coalgebra map.

These results fit in the classification of finite dimensional
Hopf algebras problem as follows.

A bialgebra with a projection is a bialgebra $E$ over a field $K$
endowed with a Hopf algebra $H$ and two bialgebra maps
$\sigma:H\to E$ and $\pi:E\to H$ such that $\pi\circ \sigma
=\Id_H$. In \cite{Rad}, M. D. Radford describes the structure of
bialgebras with a projection: $E$ can be decomposed as the smash
product of $H$ with the (right) $H$-coinvariant part of $E$ which
actually comes out to be a braided bialgebra in the monoidal category
${_{H}^{H}\mathcal{YD}}$ of Yetter-Drinfeld modules over $H$. It
is meaningful that, even relaxing some assumption on $\pi$ (as was
done by P. Schauenburg in \cite{Schauenburg1}) or on $\sigma$ (see
\cite{AMS-Spliting}), it is possible to reconstruct $E$ by means
of a suitable bosonization type procedure. An occurrence of this
situation is given by the results in \cite{AMS-Spliting} described
above.\medskip\newline The natural problem we approach in the
present paper is to show how the notion of formally smooth
(co)algebra inside monoidal categories can substitute that of
(co)separable (co)algebra in the study of splitting bialgebra
homomorphisms. This is performed investigating the relation
between formal smoothness and separability of certain functors and
led to other results related to Hopf algebra theory. Between them
we prove that the existence of $ad$-(co)invariant integrals for a
Hopf algebra $H$ is equivalent to the separability of suitable
forgetful functors (Theorem \ref{teo ad-inv}). In the finite
dimensional case, this is also equivalent to the separability of the Drinfeld Double
$D(H)$ as an extension of $H$ (Theorem \ref{teo fd ad-inv}).

Hopf algebras which are formally smooth as (co)algebras are
characterized in Propositions \ref{formal in bicomodules},
\ref{formal in comodules}, \ref{formal in bimodules} and
\ref{formal in modules} (see also \cite[Theorem 1.2]{Mas}). In
particular we obtain that the the underline (co)algebra structures
of a Hopf algebra is formally smooth if and only if it is
hereditary.

As a result we prove that given a bialgebra surjection $\pi
:E\rightarrow H$ with nilpotent kernel such that $H$ is a Hopf
algebra which is formally smooth as a $K$-algebra, then $\pi $ has
a section which is a right $H$-colinear algebra homomorphism
(Theorem \ref{teo: coweak}). Moreover, if $H$ is also endowed with
an $ad$-invariant integral, then this section can be chosen to be
$H$-bicolinear (Theorem \ref{coro 5.23}). Dually, we prove that,
if $H$ is a Hopf subalgebra of a bialgebra $E$ which is formally
smooth as a $K$-coalgebra and such that $Corad(E)\subseteq H$,
then $E$ has a weak projection onto $H$ (Theorem \ref{teo: weak}).
Furthermore, if $H$ is also endowed with an $ad$-coinvariant
integral, then this retraction can be chosen to be $H$-bilinear
(Theorem \ref{coro co5.23}). As an application, in Proposition
\ref{pro: connected} we prove that every connected Hopf algebra
$E$ over a field $K$ with $\mathrm{char} \left( K\right) =0$ has a
weak projection $\pi :E\rightarrow K\left[ x\right] $, for every
$x\in P(E)\backslash\{0\}$.\medskip\newline  The paper is
organized as follows. We begin in Section \ref{sec: Preliminary
results} by recalling the definition of monoidal category and by
listing the most important examples for this paper. We recall the
notion of projectivity (respectively injectivity) of objects in a
category $\mathfrak{C}$ with respect to a class of homomorphisms
in $\mathfrak{C}$ and some general facts about separable functors.
We obtain the main result of this section, Theorem \ref{teo F and
P-project}, providing a diagrammatic method to establish when a
separable functor $F$ preserves or reflects relative projective
(resp. injective) objects.

This technique is applied, in Section \ref{sec: (Co)separable and
formally smooth (co)algebras}, in the case when $F$ is the
forgetful functor ${_{A}\mathcal{M}_{A}}\rightarrow
{_{A}\mathfrak{M}_{A}}$, where $\mathcal{M}$ denotes one of the
categories ${\mathfrak{M}^{H}}$, ${^{H}\mathfrak{M}^{H}}$ of
right, two-sided comodules over a Hopf algebra $H$ respectively,
and $A$ is an algebra in $\mathcal{M}$. In Theorem \ref{teo
separability of F} we prove that, if $F$ is separable, then $A$ is
formally smooth as an algebra in $\mathcal{M}$ if and only if it
is formally smooth as an algebra in ${\mathfrak{M}}_{K}$ (i.e.
regardless the $H$-comodule structure of $A$). A remarkable fact
is that the functor $F$ is separable whenever $H$ has an
$ad$-invariant integral (see Lemma \ref{lem ad-inv}). In
Proposition \ref{pro sep=>left project}, a characterization of
separable algebras in a monoidal category by means of separable
functors is given. We also deal with the dual results.

In Section \ref{sec: Ad-invariant integrals through separable
functors} the existence of $ad$-invariant integrals is related to
separability of suitable functors. In particular $H$ has an
$ad$-invariant integral if and only if the forgetful functor
${_{H}^{H}\mathcal{YD}}\rightarrow {_{H}\mathfrak{M}}$ is
separable (see Theorem \ref{teo ad-inv}). In the finite
dimensional case, this is equivalent to say that the Drinfeld
Double $D(H)$ is a separable extension of $H$ (see Theorem
\ref{teo fd ad-inv}).

Section \ref{sec: Splitting algebra homomorphisms} is devoted to
the study of splitting properties of surjective algebra
homomorphisms by means of the characterization of formally smooth
algebras in monoidal categories given in \cite{AMS}. Using the
results of Section \ref{sec: Ad-invariant integrals through
separable functors}, we prove Theorem \ref{teo ad and fs} that can
be applied to the case $A=H$ where $H$ itself is a formally smooth
algebra in $\mathfrak{M}_K$ which is endowed with an
$ad$-invariant integral (Theorem \ref{coro 5.23}). Theorem
\ref{te:section} deals with the case when $H$ needs not to have an
$ad$-invariant integral but it is formally smooth as an algebra
either in ${\mathfrak{M}^{H}}$ or in ${^{H}\mathfrak{M}^{H}}$.

The main results of Section \ref{sec: Formal Smoothness of a Hopf
algebra as an algebra} are contained in Propositions \ref{formal
in bicomodules} and \ref{formal in comodules}, were we
characterize when $H$ fulfills these properties by means of a
suitable map $\tau :H^+\rightarrow H\otimes H^+$, where $H^+$ is
the augmentation ideal. Moreover a Hopf algebra $H$ comes out to
be formally smooth as a $K$-algebra if and only if it is formally
smooth as an algebra in ${\mathfrak{M}^{H}}$ if and only if it is
a hereditary $K$-algebra (note that a hereditary algebra needs not
to be formally smooth as an algebra in general, while the converse
is always true). In Theorem \ref{teo: KG is fs}, we apply these
facts to the particular case when $H$ is the group algebra $KG$
(compare with \cite[Theorem 2]{LB}).  The main application is
Theorem \ref{teo: coweak} where we prove that given a bialgebra
surjection $\pi :E\rightarrow H$ with nilpotent kernel such that
$H$ is a Hopf algebra which is formally smooth as a $K$-algebra,
then $\pi $ has a section which is a right $H$-colinear algebra
homomorphism. The results of this section are used in Section
\ref{sec: Examples} to handle some particular case related to
group algebras.

Sections \ref{sec: Ad-coinvariant integrals through separable
functors}, \ref{sec: Splitting coalgebra homomorphisms} and
\ref{sec: Formal Smoothness of a Hopf algebra as a coalgebra} are
devoted to the proof of all dual results.\medskip\newline
\textbf{Preliminaries and Notation.} In a category $\mathcal{M}$
the set of morphisms from $X$ to $Y$ will be denoted by
$\mathcal{M}(X,Y).$ If $X$ is an object in $\mathcal{M},$ then the
functor $\mathcal{M}(X,-)$ from $\mathcal{M}$ to $\mathfrak{Sets}$
associates to any morphism $u:U\rightarrow V$ in $\mathcal{M}$ the
map that will be denoted by $\mathcal{M}(X,u).$ We say that a morphism $%
f:X\rightarrow Y$ in $\mathcal{M}$ \emph{splits} (respectively
\emph{cosplits}) or has a section (resp. retraction) in
$\mathcal{M}$ whenever there is a morphism $g:Y\rightarrow X$ such
that $f\circ g=\text{Id}_Y$ (resp. $g\circ f=\text{Id}_X$). In
this case we also say that $f$ is a splitting (resp. cosplitting)
morphism.

Throughout, $K$ is a field and, when working in the category
$\mathfrak{M}=\mathfrak{M}_{K}$ of vector spaces, we write
$\otimes $ for tensor product over $K$. We use Sweedler's notation
for comultiplications $\Delta (c) =c_{(1) }\otimes c_{(2)
}=c_{1}\otimes c_{2},$ and the versions $^{C}\rho (x)
=x_{<-1>}\otimes x_{<0>}=x_{-1}\otimes x_{0}$ and $\rho ^{C}(x)
=x_{<0>}\otimes x_{<1>}=x_{0}\otimes x_{1}$ for left and right
comodules respectively (we omit the summation symbol for the sake
of brevity).

\section{Preliminary results}\label{sec: Preliminary results}

\begin{nonamethm}\textbf{Monoidal Categories.}
  Throughout this paper, the symbol $(\mathcal{M},\otimes,\mathbf{1})$ denotes a strict monoidal category with \emph{unit} $\mathbf{1}\in
\mathcal{M}$ and \emph{tensor
product} $\otimes :\mathcal{M}%
\times \mathcal{M}\rightarrow \mathcal{M}$. See \cite[Chap.
XI]{Kassel}) for a general reference.

The notions of algebra, module over an algebra, coalgebra and
comodule over a coalgebra can be introduced in the general setting
of monoidal categories. Given an algebra $A$ in $\M$ one can
define the categories $_{A}\mathcal{M}$, $\mathcal{M}_{A}$ and
$_{A}\mathcal{M}_{A}$ of left, right and two-sided modules over
$A$ respectively. Similarly, given a coalgebra $C$ in
$\mathcal{M}$, one can define the categories of $C$-comodules $^{C}\mathcal{M},\mathcal{M}%
^{C},{^{C}\mathcal{M}^{C}}$. For more details, the reader is
referred to \cite{AMS}.
\end{nonamethm}

\textbf{The relative tensor and cotensor functors}. Let
$(\mathcal{M},\otimes ,\mathbf{1})$ be a monoidal category. Assume
that $\mathcal{M}$ is abelian and let $A$ be an algebra in
$\mathcal{M}$. It can be proved (see e.g. \cite{Ar1}) that
${_A\mathcal{M}}$ is an abelian category, whenever the functor
$A\otimes(-):\mathcal{M}\to \mathcal{M}$ is additive and right
exact. In the case when both the functors $A\otimes(-):\mathcal{M}
\to\mathcal{M}$ and $(-)\otimes A:\mathcal{M}\to\mathcal{M}$ are
additive and right exact, then the category ${_A\mathcal{M}_A}$ is
abelian too.

Since, sometimes, we have to work with more than one algebra in
$\mathcal{M}$ and its bimodules, it is convenient to assume that
$X\otimes (-):\mathcal{M}\rightarrow \mathcal{M}$ and $(-)\otimes
X:\mathcal{M}\rightarrow \mathcal{M}$ are additive and right
exact, for any $X\in \mathcal{M}.$ Hence we are led to the
following definitions.

\begin{definitions}
\label{abelian assumptions}Let $\M$ be a monoidal category.\newline We say that $\mathcal{%
M}$ is an \textbf{abelian monoidal category }if $\mathcal{M}$ is
abelian and both the functors
$X\otimes (-):\mathcal{M}\rightarrow \mathcal{M}$ and $(-)\otimes X:\mathcal{%
M}\rightarrow \mathcal{M}$ are additive and right exact, for any
$X\in \mathcal{M}.$\\We say that $\mathcal{%
M}$ is a \textbf{coabelian monoidal category }if $\M^o$ is an
abelian monoidal category, where $\mathcal{M}^{o}$ denotes the dual monoidal category of $\M$. Recall that $\mathcal{M}^{o}$ and $%
\mathcal{M}$ have the same objects but $\mathcal{M}^{o}(X,Y)=\mathcal{M}%
(Y,X) $ for any $X,Y$ in $\mathcal{M}$.
\end{definitions}

Given an algebra $A$ in $\M$, there exists a suitable functor $\otimes _{A}:{_{A}\mathcal{M}_{A}}%
\times {_{A}\mathcal{M}_{A}}\rightarrow {_{A}\mathcal{M}}_{A}$ that makes the category $({_{A}\mathcal{M}}%
_{A},\otimes _{A},A)$ monoidal (an algebra in this category will
be called an $A$-\emph{algebra}\textbf{)}: see
\cite[1.11]{AMS}.\newline The tensor product over $A$ in
$\mathcal{M}$ of a right $A$-module $V$ and a left $A$-module $W$
is defined to be the coequalizer:
\begin{equation*}
\xymatrix@C=1cm{
  V\otimes A\otimes W \ar@<.5ex>[rr] \ar@<-.5ex>[rr]&& V\otimes W \ar[rr]^{_{A}\chi _{V,W}} && V\otimes _{A}W \ar[r] & 0 }
\end{equation*}
Note that, since $\otimes $ preserves coequalizers, then $V\otimes
_{A}W$ is also an $A$-bimodule, whenever $V$ and $W$ are
$A$-bimodules.\medskip\\ Dually, let $\M$ be a coabelian monoidal
category.\newline Given a coalgebra $(C,\Delta,\varepsilon)$ in
$\mathcal{M}$, there exists
of a suitable functor $\Box_{C}:{^{C}\mathcal{M}^{C}}\times {^{C}\mathcal{M}%
^{C}}\rightarrow {^{C}\mathcal{M}^{C}}$ that makes the category
$({^{C}\mathcal{M}^{C}},\Box_{C},C)$ monoidal (a coalgebra in this
category will be called a $C$-\emph{coalgebra}).\newline The
cotensor product over $C$ in $\mathcal{M}$ of a right
$C$-bicomodule $V$ and a left $C$-comodule $W$ is defined to be
the equalizer:
\begin{equation*}
\xymatrix@C=1cm{
  0 \ar[r] & V\cot_{C}W \ar[rr]^{_C\varsigma_{V,W}} && V\otimes W \ar@<.5ex>[rr] \ar@<-.5ex>[rr]&&V\ot C\ot W  }
\end{equation*} Note that, since $\otimes $ preserves equalizers, then $V\Box
_{C}W$ is also a $C$-bicomodule, whenever $V$ and $W$ are
$C$-bicomodules.\newline What follows is a list of the most
important monoidal categories meeting our
requirements.\medskip\newline \textbf{Examples of "good" monoidal
categories}.
We provide a list of the monoidal categories we need in this
paper. They are ''good'' in the sense that they are (co)abelian
monoidal categories.\medskip\\ $\bullet $ The category
$(\mathfrak{M}_{K},\otimes _{K},K)$ of all vector spaces over a
field $K$.\medskip\\ Let $(H,m_{H},u_{H},\Delta _{H},\varepsilon
_{H},S)$ be a Hopf algebra over field $K$. Then we have the
following categories (see \cite {Schauenburg2} for more details).
\newline $\bullet $ The category
$_{H}{\mathfrak{M}}=(_{H}{\mathfrak{M}},\otimes _{K},K)$, of all
left modules over $H$: the unit $K$ is a left $H$-module via
$\varepsilon _{H}$ and the tensor $V\otimes W$ of two left
$H$-modules can be regarded as an object in $_{H}{\mathfrak{M}}$
via the diagonal
action. Analogously the category ${\mathfrak{M}_{H}}$ can be introduced.%
\newline
$\bullet $ The category ${_{H}\mathfrak{M}_{H}}=({_{H}\mathfrak{M}_{H}}%
,\otimes _{K},K)$, of all two-sided modules over $H$: the unit $K$ is a $H$%
-bimodule via $\varepsilon _{H}$ and the tensor $V\otimes W$ of two $H$%
-bimodules carries, on both sides, the diagonal action.\medskip\\
We can dualize all the structures given for modules in order to
obtain categories of comodules. \newline $\bullet $ The category
${^{H}\mathfrak{M}}=({^{H}\mathfrak{M}},\otimes
_{K},K)$, of all left comodules over $H$: the unit $K$ is a left $H$%
-comodule via the map $k\mapsto 1_{H}\otimes k$ and the tensor product $%
V\otimes W$ of two left $H$-comodules can be regarded as an object in $^{H}{%
\mathfrak{M}}$ via the codiagonal coaction. Analogously the category ${%
\mathfrak{M}^{H}}$ can be introduced.\newline
$\bullet $ The category ${^{H}\mathfrak{M}^{H}}=({^{H}\mathfrak{M}^{H}}%
,\otimes _{K},K)$ of all two-sided comodules over $H$: the unit $K$ is a $H$%
-bicomodule via the maps $k\mapsto 1_{H}\otimes k$ and $k\mapsto
k\otimes 1_{H}$; the tensor $V\otimes W$ of two $H$-bicomodules
carries, on both sides, the codiagonal coaction.\medskip\\
As observed, given an algebra $A$ in an abelian monoidal category $(\mathcal{M}%
,\otimes,\mathbf{1})$, we can construct the monoidal category of $A$%
-bimodules $({_A\mathcal{M}_A},\otimes_A,A)$. Applying this (in
particular for $A:=H$) to the categories
$(\mathfrak{M}_{K},\otimes _{K},K)$,$({\mathfrak{M}^{H}},\otimes
_{K},K)$,$({^{H}\mathfrak{M}},\otimes _{K},K)$ and
$({^{H}\mathfrak{M}^{H} },\otimes_{K},K)$, we obtain
respectively:\newline $\bullet $ {\small
${_{A}\mathfrak{M}_{A}}=({_{A}\mathfrak{M}_{A}},\otimes _{A},A)$,
${_{A}\mathfrak{M}_{A}^{H}}=({_{A}\mathfrak{M}_{A}^{H}},\otimes _A,A)$, ${%
_{A}^{H}\mathfrak{M}_{A}}=({_{A}^{H}\mathfrak{M}_{A}},\otimes _A,A)$, ${%
_{A}^{H}\mathfrak{M}_{A}^{H}}=({_{A}^{H}\mathfrak{M}_{A}^{H}},\otimes
_A,A) $.}\medskip\\
Given a coalgebra $C$ in a coabelian monoidal category $(\mathcal{M},\otimes,\mathbf{1}%
)$, we can construct the monoidal category of $C$-bicomodules $({^C\mathcal{M%
}^C},\square _C,C)$. Applying this (in particular for $C:=H$) to
the categories $(\mathfrak{M}_{K},\otimes
_{K},K)$,$({\mathfrak{M}_{H}}, \otimes
_{K},K)$,$({_{H}}\mathfrak{M},\otimes _{K},K)$ and
$({_{H}}\mathfrak{M}_{H},\otimes _{K},K)$, we obtain respectively:
\newline {\small $\bullet $
${^{C}\mathfrak{M}^{C}}=({^{C}\mathfrak{M}^{C}},\square _{C},C)$,
${^{C}\mathfrak{M}_{H}^{C}}=({^{C}\mathfrak{M}_{H}^{C}},\square _{C},C)$, ${%
_{H}^{C}\mathfrak{M}^{C}}=({_{H}^{C}\mathfrak{M}^{C}},\square _{C},C)$, ${%
_{H}^{C}\mathfrak{M}_{H}^{C}}=({_{H}^{C}\mathfrak{M}_{H}^{C}},\square
_{C},C) $.}\medskip\\
It is well known that $({_{H}^{H}\mathfrak{M}_{H}^{H},\otimes }_{H},H)$ and $%
({_{H}^{H}\mathfrak{M}_{H}^{H},\square }_{H},H)$ are equivalent
monoidal categories (see \cite[Theorem 5.7]{Schauenburg2}).
\medskip\\ We now consider the categories of Yetter-Drinfeld
modules over $H$. Recall
that a \emph{twisted antipode} for $H$ is an antipode $\overline{S}$ for $%
H^{op}$ (and hence also for $H^{cop}$). One can check that $S^{-1}$ is a
twisted antipode whenever $S$ is bijective. If $H$ is commutative or
cocommutative then $S^{2}=S\circ S=\mathrm{Id}_{H}$ and consequently $%
\overline{S}=S$.\newline
$\bullet $ The category ${_{H}^{H}\mathcal{YD}}=({_{H}^{H}\mathcal{YD}}%
,\otimes _{K},K)$, of all left-left Yetter-Drinfeld modules over
$H$: the unit $K$ is a left $H$-comodule via the map $k\mapsto
1_{H}\otimes k$ and a left $H$-module via $\varepsilon _{H}$; the
tensor product $V\otimes W$ of
two left-left Yetter-Drinfeld modules can be regarded as an object in ${%
_{H}^{H}\mathcal{YD}}$ via the diagonal action and the codiagonal coaction.%
\newline
Recall that an object $V$ in ${_{H}^{H}\mathcal{YD}}$ is a left $H$-module
and a left $H$-comodule satisfying, for any $h\in H,v\in V$, the
compatibility condition:
\begin{equation*}
(h_{1}v)_{<-1>}h_{2}\otimes (h_{1}v)_{<0>}=h_{1}v_{<-1>}\otimes h_{2}v_{<0>}%
\text{ or }(hv)_{<-1>}\otimes (hv)_{<0>}=h_{1}v_{<-1>}S(h_{3})\otimes
h_{2}v_{<0>}.
\end{equation*}
Analogously the categories ${\mathcal{YD}_{H}^{H}}$,
${_{H}\mathcal{YD}}^{H}$ and ${^{H}\mathcal{YD}}_{H}$ can be
defined. The compatibility conditions are respectively:
{\small\begin{gather*}
(vh_{2})_{<0>}\otimes h_{1}(vh_{2})_{<1>}=v_{<0>}h_{1}\otimes v_{<1>}h_{2}%
\text{ or }(vh)_{<0>}\otimes (vh)_{<1>}=v_{<0>}h_{2}\otimes
S(h_{1})v_{<1>}h_{3}, \\
(h_{2}v)_{<0>}\otimes (h_{2}v)_{<1>}h_{1}=h_{1}v_{<0>}\otimes h_{2}v_{<1>}\
\text{or }(hv)_{<0>}\otimes (hv)_{<1>}=h_{2}v_{<0>}\otimes h_{3}v_{<1>}%
\overline{S}(h_{1}), \\
 h_{2}(vh_{1})_{<-1>}\otimes (vh_{1})_{<0>}=v_{<-1>}h_{1}\otimes v_{<0>}h_{2}%
\text{ or }(vh)_{<-1>}\otimes (vh)_{<0>}=\overline{S}(h_{3})v_{<-1>}h_{1}%
\otimes v_{<0>}h_{2},
\end{gather*}}for all $h\in H,v\in V$ and where in the last two cases the right conditions
are available when $H$ has a twisted antipode $\overline{S}$.

\begin{claim}
\label{def 4.2.32} \textbf{Relative Projectivity and Injectivity.}
A main tool for studying (co)separable and formally smooth
(co)algebras is relative projectivity (respectively injectivity).
Most of the material introduced below can be found in \cite[Chap.
IX, page 307-312]{HS} and \cite[Chap. 8, page
279-281]{Weibel}.\medskip\newline Let $\mathfrak{C}$ be an
arbitrary category and let $\mathcal{H}$ be a class of
homomorphisms in $\mathfrak{C}$. An object $P\in \mathfrak{C}$ is called $f$-\emph{projective} where $%
f:C_{1}\rightarrow C_{2}$ is a morphism, if $\mathfrak{C}(P,f):\mathfrak{C}%
(P,C_{1})\rightarrow \mathfrak{C}(P,C_{2})$ is surjective. $P$ is $\mathcal{H%
}$\textbf{-}\emph{projective} if it is $f$-projective for every
$f\in \mathcal{H}$. Dually, an object $I\in \mathfrak{C}$ is called $f$-\emph{injective}, where $%
f:C_{1}\rightarrow C_{2}$ is a morphism, if and only if,
considered as an object in the
opposite category $\mathfrak{C}^{op}$, it is $f^{op}$-projective, where $%
f^{op}:C_{2}\rightarrow C_{1}$ is in $\mathfrak{C}^{op}$. $I$ is called $%
\mathcal{H}$\textbf{-}\emph{injective }if it is $f$-injective for every $%
f\in \mathcal{H}$.\newline All the results we will obtain for
projectivity, can be dualized to get their analogues for
injectivity.
\end{claim}

\begin{theorem}
\label{teo 4.2.33}Let $\mathbb{H}:\mathfrak{B}\rightarrow
\mathfrak{A}$ be a covariant functor and consider:
\begin{equation*}
\mathcal{E}_{\mathbb{H}}:=\{f\in \mathfrak{B}\mid
\mathbb{H}(f)\text{ splits in }\mathfrak{A}\}.
\end{equation*}
Let $\mathbb{T}:\mathfrak{A}\rightarrow \mathfrak{B}$ be a left adjoint of $%
\mathbb{H}$ and let $\varepsilon :\mathbb{TH}\rightarrow \mathrm{Id}_{%
\mathfrak{B}}$ be the counit of the adjunction.\newline Then, for
any object $P\in \mathfrak{B}$, the following assertions are
equivalent:

$(a)$ $P$ is $\mathcal{E}_{\mathbb{H}}$-projective.

$(b)$ Every morphism $f:B\rightarrow P$ in
$\mathcal{E}_{\mathbb{H}}$ has a section.

$(c)$ $\varepsilon _{P}:\mathbb{TH}P\rightarrow P$ has a section.

$(d)$ There is a splitting morphism $\pi :\mathbb{T}X\rightarrow
P$ for a suitable object $X\in \mathfrak{A}$.\newline In
particular all objects of the form $\mathbb{T}X$, $X\in
\mathfrak{A}$, are $\mathcal{E}_{\mathbb{H}}$-projective.
\end{theorem}

\begin{proof}
Let $\eta :\mathrm{Id}_{\mathfrak{A}}\rightarrow \mathbb{HT}$ be
the unit of the adjunction. \newline
$(a)\Rightarrow (b).$ Assume that $P\in \mathfrak{B}$ is $\mathcal{E}_{%
\mathbb{H}}$-projective i.e. that for every $f:B\rightarrow B_{2}$ in $%
\mathcal{E}_{\mathbb{H}}$ and for every morphism $\gamma
:P\rightarrow B_{2}, $ there exists a morphism $\beta
:P\rightarrow B$ such that $\gamma =f\circ \beta $. In particular,
for $B_{2}:=P$ and $\gamma :=\mathrm{Id}_{P},$ there exists a
morphism $\beta :P\rightarrow B$ such that $\mathrm{Id}_{P}=f\circ
\beta $.\newline
$(b)\Rightarrow (c).$ Since $\mathbb{H}(\varepsilon _{B})\circ \eta _{%
\mathbb{H}B}=\mathrm{Id}_{\mathbb{H}B}$, we infer that $\mathbb{H}%
(\varepsilon _{B})$ splits and hence the counit $\varepsilon _{B}:\mathbb{TH}%
B\rightarrow B$ belongs to $\mathcal{E}_{\mathbb{H}}$ for any
$B\in \mathfrak{B}$.$\newline (c)\Rightarrow (d).$
Obvious.\newline
$(d)\Rightarrow (a)$. Let $f:B_{1}\rightarrow B_{2}$ be in $\mathcal{E}_{%
\mathbb{H}}$ and denote by $g:\mathbb{H}B_{2}\rightarrow
\mathbb{H}B_{1}$ the section of $\mathbb{H}(f)$. Let $\gamma
:P\rightarrow B_{2}$. Assume that $\pi :\mathbb{T}X\rightarrow P$
is a split morphism for a suitable object $X\in \mathfrak{A.}$ Let
$\sigma :P\rightarrow \mathbb{T}X$ be a
section of $\pi $ and $\tau :P\rightarrow B_{1}$ be defined by%
\begin{equation*}
P\overset{\sigma }{\rightarrow
}\mathbb{T}X\overset{\mathbb{T}\left( \eta _{X}\right)
}{\longrightarrow }\mathbb{THT}X\overset{\mathbb{TH}\left( \pi
\right) }{\longrightarrow }\mathbb{TH}P\overset{\mathbb{TH}\left(
\gamma \right) }{\longrightarrow
}\mathbb{TH}B_{2}\overset{\mathbb{T}\left(
g\right) }{\rightarrow }\mathbb{TH}B_{1}\overset{\varepsilon _{B_{1}}}{%
\rightarrow }B_{1}.
\end{equation*}%
We have%
\begin{eqnarray*}
f\circ \tau  &=&f\circ \varepsilon _{B_{1}}\circ \mathbb{T}\left(
g\right) \circ \mathbb{TH}\left( \gamma \right) \circ
\mathbb{TH}\left( \pi \right)
\circ \mathbb{T}\left( \eta _{X}\right) \circ \sigma  \\
&=&\varepsilon _{B_{2}}\circ \mathbb{T}\left[ \mathbb{H}\left(
f\right) \circ \left( g\right) \right] \circ \mathbb{TH}\left(
\gamma \circ \pi
\right) \circ \mathbb{T}\left( \eta _{X}\right) \circ \sigma  \\
&=&\varepsilon _{B_{2}}\circ \mathbb{TH}\left( \gamma \circ \pi
\right) \circ \mathbb{T}\left( \eta _{X}\right) \circ \sigma
=\gamma \circ \pi \circ \varepsilon _{\mathbb{T}X}\circ
\mathbb{T}\left( \eta _{X}\right) \circ \sigma =\gamma \circ \pi
\circ \sigma =\gamma
\end{eqnarray*}%
and hence $P$ is $\mathcal{E}_{\mathbb{H}}$-projective.

Since $\Id _{\mathbb{T}X}:\mathbb{T}X\to \mathbb{T}X$ is an
isomorphism, by $(d)\Rightarrow (a),$ we have that $\mathbb{T}X$ is $%
\mathcal{E}_{\mathbb{H}}$-projective.
\end{proof}

For completeness we include the dual statement of Theorem \ref{teo 4.2.33}.

\begin{theorem}
\label{teo dual inj}Let $\ \mathbb{T}:\mathfrak{A}\rightarrow
\mathfrak{B}$ be a covariant functor and consider:
\begin{equation*}
\mathcal{I}_{\mathbb{T}}:=\{g\in \mathfrak{A}\mid
\mathbb{T}(g)\text{ cosplits in }\mathfrak{B}\}.
\end{equation*}
Let $\mathbb{H}:\mathfrak{B}\rightarrow \mathfrak{A}$ be a right
adjoint of $\mathbb{T}$ and let $\eta
:\mathrm{Id}_{\mathfrak{A}}\rightarrow \mathbb{HT} $ be the unit
of the adjunction.\newline Then, for any object $I\in
\mathfrak{A}$, the following assertions are equivalent:

$(a)$ $I$ is $\mathcal{I}_{\mathbb{T}}$-injective.

$(b)$ Every morphism $f:I\to A$ in $\mathcal{I}_{\mathbb{T}}$ has
a retraction.

$(c)$ $\eta _{I}:I\rightarrow \mathbb{HT}I$ has a retraction.

$(d)$ There is a cosplitting morphism $i:I\rightarrow \mathbb{H}Y$
for a suitable object $Y\in \mathfrak{B}$.\newline In particular
all objects of the form $\mathbb{H}Y$, $Y\in \mathfrak{B}$, are
$\mathcal{I}_{\mathbb{T}}$-injective.
\end{theorem}

\begin{claim}\label{faithful}\textbf{Separable Functors}. Let $\mathbb{U}:\mathfrak{B}\rightarrow
\mathfrak{A}$ be a covariant functor. We have functors
\begin{equation*}
Hom_{\mathfrak{B}}(\bullet ,\bullet ),Hom_{\mathfrak{A}}({\mathbb{U}}%
(\bullet ),{\mathbb{U}}(\bullet )):\mathfrak{B}^{op}\times \mathfrak{B}%
\rightarrow \mathfrak{Sets}
\end{equation*}
and a natural transformation
\begin{equation*}
\mathcal{U}:Hom_{\mathfrak{B}}(\bullet ,\bullet )\rightarrow Hom_{%
\mathfrak{A}}({\mathbb{U}}(\bullet ),{\mathbb{U}}(\bullet )),\quad\mathcal{U}_{{B}%
_{1},{B}_{2}}({f}):={\mathbb{U}}({f})\text{ for all objects }{B}_{1},{B}%
_{2}\in \mathfrak{B}.
\end{equation*}
The functor ${\mathbb{U}}$ is called \emph{separable }if
$\mathcal{U}$ cosplits, that is there is a natural
transformation $$\mathcal{P}:Hom_{\mathfrak{A}}({\mathbb{U}}(\bullet ),{%
\mathbb{U}}(\bullet ))\rightarrow Hom_{\mathfrak{B}}(\bullet
,\bullet )$$
such that $\mathcal{P\circ {\mathcal{U}}}=\mathbf{1}_{Hom_{\mathfrak{B}%
}(\bullet ,\bullet )},$ the identity natural transformation on $Hom_{%
\mathfrak{B}}(\bullet ,\bullet ).$ \newline It is proved in
\cite[page 1446]{Rafael} that this definition is consistent with
the one given in \cite{NdO}.
%
%
\end{claim}

\begin{remark}
\label{rem faithful} Let $\alpha :X\rightarrow Y$ be a morphism in $%
\mathfrak{B}$. If $\mathbb{U}$ is a faithful functor, then,
$\alpha $ is an epimorphism (resp. monomorphism) whenever
$\mathbb{U}(\alpha )$ is.
\end{remark}

Let us recall some well known property on separable functors.

\begin{lemma}\cite[Proposition 1.2]{NdO}
\label{lem separable}Let $\mathbb{U}:\mathfrak{B}\rightarrow
\mathfrak{A}$ be a covariant separable functor and let
$\alpha:X\rightarrow Y$ be a morphism in $\mathfrak{B}$. If
${\mathbb{U}}(\alpha )$ has a section $h$ (resp. a retraction $l$)
in $\mathfrak{A}$, then $\alpha $ has a section (retraction) in
$\mathfrak{B}.$
\end{lemma}

\begin{lemma}
\label{lemma non so}Let $F:\mathfrak{A}\rightarrow \mathfrak{B}$ and $G:%
\mathfrak{B}\rightarrow \mathfrak{C}$ be covariant functors. Then $\mathcal{E%
}_{F}\subseteq \mathcal{E}_{GF}$ and $\mathcal{I}_{F}\subseteq \mathcal{I}%
_{GF}.$ Moreover the equalities hold whenever $G$ is separable.
\end{lemma}

\begin{theorem}
\label{teo compos of separable}Consider functors $\mathbb{T}:\mathfrak{A}%
\rightarrow \mathfrak{B}$ and $\mathbb{H}:\mathfrak{B}\rightarrow %
\mathfrak{C}.$ Then, we have that:

$1)$ If $\mathbb{T}$ and $\mathbb{H}$ are separable, then $\mathbb{H}\circ
\mathbb{T}$ is also separable.

$2)$ If $\mathbb{H}\circ \mathbb{T}$ is separable, then $\mathbb{T}$ is
separable.

$3)$ If $\mathfrak{C}=\mathfrak{A}$ and that $(\mathbb{T},\mathbb{H}) $ is a
category equivalence, then $\mathbb{T}$ and $\mathbb{H}$ are both separable.
\end{theorem}

\begin{proof}
See \cite[Proposition 46 and Corollary 9]{CMZ}.
\end{proof}

We quote from \cite{Rafael} the so-called Rafael Theorem:

\begin{theorem}
\label{teo Rafael}\cite[Theorem 1.2]{Rafael} Let
$(\mathbb{T},\mathbb{H})$ be an
adjunction, where $\mathbb{T}:\mathfrak{A}\rightarrow \mathfrak{B}$ and $%
\mathbb{H}:\mathfrak{B}\rightarrow \mathfrak{A}$. Then we have:

1) $\mathbb{T}$ is separable if and only if the unit $\eta :\mathrm{Id}_{\mathfrak{A}%
}\rightarrow \mathbb{H}\mathbb{T}$ of the adjunction cosplits, i.e. there
exists a natural transformation $\mu :\mathbb{H}\mathbb{T}\rightarrow
\mathrm{Id}_{\mathfrak{A}}$ such that $\mu \circ \eta =\mathrm{Id}_{\mathrm{%
Id}_{\mathfrak{A}}},$ the identity natural transformation on $\mathrm{Id}_{%
\mathfrak{A}}$.

2) $\mathbb{H}$ is separable if and only if the counit $\varepsilon :\mathbb{T}\mathbb{H%
}\rightarrow \mathrm{Id}_{\mathfrak{B}}$ of the adjunction splits, i.e.
there exists a natural transformation $\sigma :\mathrm{Id}_{\mathfrak{B}%
}\rightarrow \mathbb{T}\mathbb{H}$ such that $\varepsilon \circ \sigma =%
\mathrm{Id}_{\mathrm{Id}_{\mathfrak{B}}},$ the identity natural
transformation on $\mathrm{Id}_{\mathfrak{B}}$.
\end{theorem}

\begin{corollary}
\label{coro Rafael}Let $(\mathbb{T},\mathbb{H})$ be an adjunction, where $%
\mathbb{T}:\mathfrak{A}\rightarrow \mathfrak{B}$ and $\mathbb{H}:\mathfrak{B}%
\rightarrow \mathfrak{A}$. Then we have:

1) $\mathbb{H}$ separable $\Rightarrow $ any object in $\mathfrak{B}$ is $%
\mathcal{E}_{\mathbb{H}}$-projective.

2) $\mathbb{T}$ separable $\Rightarrow $ any object in $\mathfrak{A}$ is $%
\mathcal{I}_{\mathbb{T}}$-injective.
\end{corollary}

\begin{proof}1) Let $B$ be an object in $\mathfrak{B}$. Since $\mathbb{H}(\varepsilon _{B})\circ \eta _{\mathbb{H}B}=%
\mathrm{Id}_{\mathbb{H}B}$ and $\mathbb{H}$ is separable, by Lemma
\ref{lem separable}, $\varepsilon _{B}$ has a section in
$\mathfrak{B}$. By Theorem \ref{teo 4.2.33}, $B$ is
$\mathcal{E}_{\mathbb{H}}$-projective.\newline 2) It follows
analogously by Lemma \ref{lem separable} and Theorem \ref{teo dual
inj} once observed that $\varepsilon_{\mathbb{T}A}\circ
\mathbb{T}(\eta_A)=\textrm{Id}_{\mathbb{T}A}$ for any $A\in
\mathfrak{A}$.
\end{proof}

We are now ready to prove the main theorem of this section, that
investigates whether a functor $F$ (resp. $F^{\prime }$) preserves and
reflects relative projective (resp. injective) objects.

\begin{theorem}
\label{teo F and P-project} Let $(\mathbb{T},\mathbb{H})$ and $(\mathbb{%
T^{\prime}},\mathbb{H^{\prime}})$ be adjunctions and assume that,
in the
following diagrams, ${\mathbb{T}^{\prime}}\circ F^{\prime}$ and $F\circ {%
\mathbb{T}}$ (and also $F^{\prime}\circ {\mathbb{H}}$ and ${\mathbb{H}%
^{\prime}}\circ F$) are naturally equivalent:
\begin{equation*}
\xymatrix@R=15pt@C=50pt{
  \mathfrak{A} \ar[d]_{\mathbb{T}} \ar[r]^{F'} & \mathfrak{A'} \ar[d]^{\mathbb{T}'} \\
  \mathfrak{B} \ar[r]_{F} & \mathfrak{B'}   }
\text{ \qquad } \xymatrix@R=15pt@C=50pt{
  \mathfrak{A} \ar[r]^{F'} & \mathfrak{A'}  \\
  \mathfrak{B} \ar[u]^{\mathbb{H}} \ar[r]_{F} & \mathfrak{B'} \ar[u]_{\mathbb{H}'}  }
\end{equation*}
Let $P$ be an object in $\mathfrak{B}$ and let $I$ be an object in $%
\mathfrak{A}$. We have:

a) $P$ is $\mathcal{E}_{\mathbb{H}}$-projective $\Longrightarrow $ $F(P)$ is
$\mathcal{E}_{\mathbb{H}^{\prime }}$-projective; the converse is true
whenever $F$ is separable.

a$^{op}$) $I$ is $\mathcal{I}_{\mathbb{T}}$-injective $\Longrightarrow $ $%
F^{\prime }(I)$ is $\mathcal{I}_{\mathbb{T}^{\prime }}$-injective; the
converse is true whenever $F^{\prime }$ is separable.
\end{theorem}

\begin{proof}
a) Let $\varepsilon :\mathbb{T}\mathbb{H}\rightarrow
\mathrm{Id}_{\mathfrak{B}}$ be the counit of the adjunction
$(\mathbb{T},\mathbb{H})$.\newline Assume that $P$ is
$\mathcal{E}_{\mathbb{H}}$-projective. Then, by Theorem \ref {teo
4.2.33}, $\varepsilon _{P}:\mathbb{T}\mathbb{H}P\rightarrow P$ has
a section $\beta :P\rightarrow \mathbb{T}\mathbb{H}P,$ i.e.
$\varepsilon _{P}\circ \beta =\mathrm{Id}_{P}.$ Since $F(\beta )$
is a section of $F(\varepsilon _{P}):\mathbb{T}^{\prime
}\mathbb{H}^{\prime }FP\sim F\mathbb{T}\mathbb{H}P\rightarrow FP$,
by applying Theorem \ref{teo 4.2.33} to the adjunction
$(\mathbb{T}^{\prime },\mathbb{H}^{\prime }) $ in the case when
$X=\mathbb{H}^{\prime }FP$ and to the split morphism
$F(\varepsilon _{P})$, we conclude that $FP$ is
$\mathcal{E}_{\mathbb{H}^{\prime }}$-projective.
\newline Conversely, assume that $FP$ is $\mathcal{E}_{\mathbb{H}^{\prime
}}$-projective and that $F$ is separable. Let $\eta
:\mathrm{Id}_{\mathfrak{B}}\rightarrow \mathbb{H}\mathbb{T}$ be
the unit of the adjunction $(\mathbb{T},\mathbb{H}).$ Thus
$\mathbb{H}(\varepsilon _{P})\circ \eta
_{\mathbb{H}P}=\mathrm{Id}_{\mathbb{H}P}$ and hence $F^{\prime
}(\eta _{\mathbb{H}P})$ is a section of $F'\mathbb{H}(\varepsilon
_{P})$. Then also $\mathbb{H}^{\prime }F(\varepsilon _{P})$ has a
section, so that $F(\varepsilon
_{P}):F\mathbb{T}\mathbb{H}P\rightarrow FP$ belongs to
$\mathcal{E}_{\mathbb{H}^{\prime }}$. As $FP$ is
$\mathcal{E}_{\mathbb{H}^{\prime }}$-projective, by Theorem \ref
{teo 4.2.33}, we get a section in $\mathfrak{B}^{\prime }$ of
$F(\varepsilon _{P}).$ Since $F$ is separable, by Lemma \ref{lem
separable}, we conclude that $\varepsilon _{P}$ splits in\
$\mathfrak{B}$: hence $P$ is
$\mathcal{E}_{\mathbb{H}}$-projective.\newline a$^{op}$) It
follows by duality.
\end{proof}

\section{(Co)separable and formally smooth (co)algebras}\label{sec: (Co)separable and formally smooth
(co)algebras}


\begin{claim}
\label{5} \label{cinque}Let $(A,m,u)$ be an algebra in a monoidal
category $(\mathcal{M},\otimes ,\mathbf{1}).$ We have the functors
\begin{gather*}
{_{A}\mathbb{T}}:\mathcal{M}\rightarrow {_{A}\mathcal{M}}\text{ where }{_{A}%
\mathbb{T}}(X):=A\otimes X\text{ and }{_{A}\mathbb{T}}(f):=A\otimes f, \\
{\mathbb{T}_{A}}:\mathcal{M}\rightarrow \mathcal{M}_{A}\text{ where }{%
\mathbb{T}_{A}}(X):=X\otimes A\text{ and }{\mathbb{T}_{A}}(f):=f\otimes A, \\
{_{A}\mathbb{T}_{A}}:\mathcal{M}\rightarrow {_{A}\mathcal{M}}_{A}\text{
where }{_{A}\mathbb{T}_{A}}(X):=A\otimes (X\otimes A)\text{ and }{_{A}%
\mathbb{T}_{A}}(f):=A\otimes (f\otimes A),
\end{gather*}
with their right adjoint (see \cite[Proposition 1.6]{AMS}) ${_{A}}\mathbb{%
H},{\mathbb{H}_{A}},{_{A}\mathbb{H}_{A}}$, respectively, that
forget the
module structures. \label{4.2.39}Then the adjunctions $({\mathbb{T}_{A}},%
\mathbb{H}_{A})$, $(_{A}{\mathbb{T}},{_{A}\mathbb{H}})$ and $(_{A}{\mathbb{T}%
}_{A},_{A}\mathbb{H}_{A})$, give rise to the following classes:
\begin{gather*}
\mathcal{E}_{A} :=\mathcal{E}_{\mathbb{H}_{A}}=\{g\in \mathcal{M}_{A}\mid g%
\text{ splits in }\mathcal{M}\}, \\
_{A}\mathcal{E} :=\mathcal{E}_{_{A}\mathbb{H}}=\{g\in {_{A}\mathcal{M}}%
\mid g\text{ splits in }\mathcal{M}\}, \\
{_{A}\mathcal{E}_{A}} :=\mathcal{E}_{_{A}\mathbb{H}_{A}}=\{g\in {_{A}%
\mathcal{M}_{A}}\mid g\text{ splits in }\mathcal{M}\}.
\end{gather*}
\label{def separable algebra}Recall
that an algebra $(A,m,u)$ is called \emph{separable\ }in {$\mathcal{M}$}%
\emph{\ }whenever the multiplication $m$ admits a section
$A\rightarrow A\otimes A$ in ${_{A}}${$\mathcal{M}$}${_{A}}$.

Assume that $\mathcal{M}$ is an abelian monoidal category. Then
$
(\Omega ^{1}A,j):=\ker m
$
carries a natural $A$-bimodule structure that makes it the kernel
of $m$ in the category ${_{A}}${$\mathcal{M}$}${_{A}}$. We say
that $A$ is \emph{formally smooth} in {$\mathcal{M}$ }(see
\cite[Corollary 3.12]{AMS}) if and only if $\Omega ^{1}A$ is an
${_{A}\mathcal{E}_{A}}$-projective $A$-bimodule.
\end{claim}

Let us recall the following result that holds true for unitary
rings.

\begin{proposition}\cite[Proposition 1.3]{NdO}
\label{lem sep=>left project for rings}For any ring homomorphism $%
i:S\rightarrow R,$ the following are equivalent:

$(1) $ $R$ is separable in $({_{S}}\mathfrak{M}_{S},\otimes
_{S},S)$, i.e. $R/S$ is separable.

$(2)$ The restriction of scalars functor ${_{R}}\mathfrak{M}\rightarrow {_{S}%
\mathfrak{M}}$ is separable.

$(3)$ The restriction of scalars functor $\mathfrak{M}_{R}\rightarrow %
\mathfrak{M}_{S}$ is separable.
\end{proposition}

As we will explain in Remark \ref{rem: one side funtor}, the
previous result, in general, can not be extended to algebras in a
monoidal category.

\begin{lemma}
\label{lem one side}Let $A$ be a separable algebra in a monoidal category {$%
\mathcal{M}$}. The following assertions hold true:

1) The forgetful functor $_{A}\mathbb{H}:{_{A}\mathcal{M}}\rightarrow
\mathcal{M}$ is separable. In particular, any left $A$-module $(M{,{^{A}}}{%
\mu }_{M})$ is $_{A}\mathcal{E}$-projective. Moreover if $M$ is an $A$%
-bimodule, the multiplication ${{^{A}}}{\mu }_{M}:A\otimes M\rightarrow M$
has a section $^{A}\sigma_M $ which is $A$-bilinear and natural in $M$.

2) The forgetful functor $\mathbb{H}_{A}:{\mathcal{M}_{A}}\rightarrow
\mathcal{M}$ is separable. In particular, any right $A$-module $(M{,\mu }%
_{M}^{A})$ is $\mathcal{E}_{A}$-projective. Moreover if $M$ is an $A$%
-bimodule, the multiplication ${\mu }_{M}^{A}:M\otimes A\rightarrow M$ has a
section $\sigma ^{A}_M$ which is $A$-bilinear and natural in $M$.
\end{lemma}

\begin{proof}
1) By assumption, the multiplication $m$ of $A$ admits a section
$\nu
:A\rightarrow A\otimes A$ in ${_{A}}${$\mathcal{M}$}${_{A}.}$ Let $(M,{%
^{A}\mu _{M}})$ be a left $A$-module and consider the morphism
$^{A}\sigma
_{M}:M\rightarrow A\otimes M$ defined by $^{A}\sigma _{M}:=(A\otimes {{^{A}}}%
\mu _{M})\circ (\nu u\otimes M)\circ l_{M}^{-1}$, where $u:\mathbf{1}%
\rightarrow A$ is the unit of $A$. It is straightforward to check that $%
^{A}\sigma _{M}$ is a left $A$-linear section of ${^{A}\mu _{M}}$
which is $A $-bilinear whenever $M\in {_{A}\mathcal{M}_{A}}$ (see
the left handed version of \cite[Lemma 1.29]{AMS}). Since
${{^{A}}}{\mu }$ is the counit of the adjunction
$(_{A}{\mathbb{T}},_{A}\mathbb{H}),$ and ${^{A}}\sigma _{M}$
defines a natural transformation ${^{A}}\sigma :\mathrm{Id}_{{_{A}\mathcal{M}%
}}\rightarrow \left( {_{A}{\mathbb{T}}}\right) \left( {_{A}\mathbb{H}}%
\right) $, we get, by Theorem \ref{teo Rafael}, that
$_{A}\mathbb{H}$ is separable. Note that, by Corollary \ref{coro
Rafael}, if the forgetful functor
$_{A}\mathbb{H}:{_{A}\mathcal{M}}\rightarrow \mathcal{M}$ is
separable, then any left $A$-module is
$_{A}\mathcal{E}$-projective.\newline 2) It follows analogously.
\end{proof}


\begin{proposition}
\label{2 sep functors}Let $H$ be a Hopf algebra over a
field $K$. The forgetful functors $\mathfrak{M}^H_H\rightarrow \mathfrak{M}%
^H $ and $^H\mathfrak{M}^H_H\rightarrow {^H\mathfrak{M}^H}$ are
separable.
\end{proposition}

\begin{proof}
Composing the functor $\left( -\right)
^{coH}:\mathfrak{M}^{H}\rightarrow
\mathfrak{M}_{K}$ with the forgetful functor $\mathfrak{M}%
_{H}^{H}\rightarrow \mathfrak{M}^{H},$ one gets the Sweedler's
equivalence of categories $\left( -\right)
^{coH}:\mathfrak{M}_{H}^{H}\rightarrow \mathfrak{M}_{K}.$ Since,
by Theorem \ref{teo Rafael}, this functor is
separable, by Theorem \ref{teo compos of separable}, the forgetful functor $%
\mathfrak{M}_{H}^{H}\rightarrow \mathfrak{M}^{H}$ is separable
too.\newline
Composing the functor $\left( -\right) ^{coH}:{^{H}\mathfrak{M}^{H}}%
\rightarrow {^{H}\mathfrak{M}}$ with the forgetful functor $^{H}\mathfrak{M}%
_{H}^{H}\rightarrow {^{H}\mathfrak{M}^{H}},$ one gets the
Sweedler's
equivalence of categories $\left( -\right) ^{coH}:{^{H}\mathfrak{M}_{H}^{H}}%
\rightarrow {^{H}\mathfrak{M}}.$ As in the first part, we conclude
that the forgetful functor $^{H}\mathfrak{M}%
_{H}^{H}\rightarrow {^{H}\mathfrak{M}^{H}}$ is separable.
\end{proof}

\begin{remark}\label{rem: one side funtor}
By Lemma \ref{lem one side}, the forgetful functor $\mathbb{H}_{A}:{\mathcal{%
M}_{A}}\rightarrow \mathcal{M}$ is separable for any separable
algebra $A$ in a monoidal category $\mathcal{M}$. The converse
does not hold true. In
fact, when $\mathcal{M}=\mathfrak{M}^H$ and $A=H$, the functor $%
\mathbb{H}_{A}$ is always separable (Proposition \ref{2 sep
functors}), but $A$ is separable in $\mathcal{M}$ if and only if
$H$ is a semisimple algebra (\cite[Proposition
2.11]{AMS-Spliting}).
\end{remark}

\begin{proposition}
\label{pro sep=>left project}Let $A$ be an algebra in a monoidal category {$%
\mathcal{M.}$ }The following assertions are equivalent:

$(a)$ $A$ is separable in $\mathcal{M}$.

$(b)$ The forgetful functor $_{A}\mathbb{H}_{A}:{_{A}\mathcal{M}_{A}}%
\rightarrow \mathcal{M}$ is separable.

$(c)$ Any $A$-bimodule is $_{A}\mathcal{E}_{A}$-projective.

$(d)$ The $A$-bimodule $A$ is $_{A}\mathcal{E}_{A}$-projective.
\end{proposition}

\begin{proof}
$(a)\Rightarrow (b)$ If $(M,{{^{A}}}{\mu }_{M},\mu _{M}^{A})$ is an $A$%
-bimodule, by Lemma \ref{lem one side}, there are $A$-bilinear natural
sections ${^{A}}\sigma _{M}$ and $\sigma _{M}^{A},$ respectively of ${{^{A}}}%
{\mu }_{M}$ and $\mu _{M}^{A}$. The morphism $\sigma _{M}:=({^{A}}\sigma
_{M}\otimes A)\circ \sigma _{M}^{A}:M\rightarrow A\otimes M\otimes A$ is a
section in ${_{A}}${$\mathcal{M}$}${_{A}}$ of the counit $%
\varepsilon _{M}:=\mu _{M}^{A}\circ ({{^{A}}\mu }_{M}\otimes A):A\otimes
M\otimes A\rightarrow M$ of the adjunction $(_{A}{\mathbb{T}}_{A},_{A}%
\mathbb{H}_{A}).$ Since $\sigma _{M}$ is natural in $M,$ we get a
natural transformation $\sigma
:\mathrm{Id}_{\mathcal{M}}\rightarrow
{_{A}\mathbb{T}_{A}}{_{A}\mathbb{H}_{A}}$ such that $\varepsilon
\circ \sigma =\mathrm{Id}_{\mathrm{Id}_{\mathcal{M}}}.$ We
conclude by Theorem \ref{teo Rafael}.\newline $(b)\Rightarrow (c)$
It follows by Corollary \ref{coro Rafael}.\newline $(c)\Rightarrow
(d)$ Obvious.\newline $(d)\Rightarrow (a)$ Since $A$ is
$_{A}\mathcal{E}_{A}$-projective, the multiplication $m:A\otimes
A\rightarrow A$, that is a morphism in $_{A}\mathcal{E}_{A},$
admits a section $\sigma :A\rightarrow A\otimes A$ in ${_{A}}${$\mathcal{M}$}%
${_{A}.}$
\end{proof}

\begin{corollary}\label{coro: sep=>fs}
Any separable algebra in an abelian monoidal category
$\mathcal{M}$ is formally smooth.
\end{corollary}

\begin{corollary}
\label{coro proj}Let $A$ be a separable algebra in $\mathfrak{M}_K$. Then
any left $A$-module is projective in $_A\mathfrak{M}$. Hence any left $A$%
-module is also injective in $_A\mathfrak{M}$ and $A$ is
semisimple.
Moreover any $A$-bimodule is projective in $_A\mathfrak{M}_A$ and hence any $%
A$-bimodule is injective in $_A\mathfrak{M}_A$.
\end{corollary}

\begin{proof}
Since $\mathcal{M}=\mathfrak{M}_K$, any epimorphism in
$\mathcal{M}$ splits. So a left $A$-module is
$_A\mathcal{E}$-projective if and only if it is projective in
$_A\mathfrak{M}$ in the usual sense. The right and two-sided cases
follow analogously.
\end{proof}

\begin{claim}
\label{ex: monoidal functor} Let $(F^{\prime },\phi _{0},\phi _{2}):(%
\mathcal{M},\otimes ,\mathbf{1)\rightarrow (}\mathcal{M}^{\prime }\mathfrak{,%
}\otimes ^{\prime },\mathbf{1}^{\prime }\mathbf{)}$ be a monoidal
functor between two monoidal categories, where $\phi
_{2}(U,V):F^{\prime }(U\otimes
V)\rightarrow F^{\prime }(U)\otimes ^{\prime }F^{\prime }(V),$ for any $%
U,V\in \mathcal{M}$, and $\phi _{0}:\mathbf{1}^{\prime
}\rightarrow F^{\prime }(\mathbf{1}).$ Let $(A,m,u)$ be an algebra
in $\mathcal{M}$. It is well known that $(A^{\prime },m_{A^{\prime
}},u_{A^{\prime }}):=(F^{\prime
}(A),m_{F^{\prime }(A)},u_{F^{\prime }(A)})$ is an algebra in $\mathcal{M}%
^{\prime }$, where
\begin{equation*}
\begin{tabular}{l}
$m_{F^{\prime }(A)}:=F^{\prime }(A)\otimes ^{\prime }F^{\prime }(A)\overset{%
\phi _{2}(A,A)}{\longrightarrow }F^{\prime }(A\otimes
A)\overset{F^{\prime
}(m)}{\longrightarrow }F^{\prime }(A)$ \\
$u_{F^{\prime }(A)}:=\mathbf{1}^{\prime }\overset{\phi
_{0}}{\longrightarrow }F^{\prime }(\mathbf{1})\overset{F^{\prime
}(u)}{\longrightarrow }F^{\prime
}(A)$.%
\end{tabular}%
\end{equation*}%
Consider the functor $F:{_{A}\mathcal{M}_{A}}\rightarrow {_{A^{\prime }}%
\mathcal{M}_{A^{\prime }}^{\prime }}$ defined by {\small $F((M,{^{A}\mu _{M}}%
,{\mu _{M}^{A}}))=(F^{\prime }(M),{^{A^{\prime }}\mu _{F^{\prime
}(M)}},{\mu _{F^{\prime }(M)}^{A^{\prime }}}),$} where
\begin{equation*}
\begin{tabular}{l}
${^{A^{\prime }}\mu _{F^{\prime }(M)}}:=F^{\prime }(A)\otimes
^{\prime }F^{\prime }(M)\overset{\phi _{2}(A,M)}{\longrightarrow
}F^{\prime
}(A\otimes M)\overset{F^{\prime }({^{A}\mu _{M}})}{\longrightarrow }%
F^{\prime }(M)$ \\
${\mu _{F^{\prime }(M)}^{A^{\prime }}}:=F^{\prime }(M)\otimes
^{\prime }F^{\prime }(A)\overset{\phi _{2}(M,A)}{\longrightarrow
}F^{\prime
}(M\otimes A)\overset{F^{\prime }({\mu _{M}^{A}})}{\longrightarrow }%
F^{\prime }(M)$.%
\end{tabular}%
\end{equation*}
\end{claim}

Let us study a particular case of Theorem \ref{teo F and
P-project}.

\begin{proposition}
\label{pro separability of F gen}Let ${\mathcal{M}}$ and ${\mathcal{M}}%
^{\prime }$ be abelian monoidal categories. Let $A$, $A^{\prime}$, $%
F^{\prime}$ and $F$ as in \ref{ex: monoidal functor}. Then, in the
following diagrams, ${\mathbb{T}^{\prime}}\circ F^{\prime}$ and
$F\circ {\mathbb{T}}$
are naturally equivalent and $F^{\prime}\circ {\mathbb{H}}={\mathbb{H}%
^{\prime}}\circ F$:
\begin{equation*}
\xymatrix@R=15pt@C=50pt{
  {\mathcal{M}} \ar[d]_{\mathbb{T}} \ar[r]^{F'} & {\mathcal{M}'} \ar[d]^{\mathbb{T}'} \\
  {_{A}\mathcal{M}_{A}} \ar[r]_{F} & {_{A'}\mathcal{M}'_{A'}}   }
\text{ \qquad } \xymatrix@R=15pt@C=50pt{
  {\mathcal{M}} \ar[r]^{F'} & {\mathcal{M}'}  \\
  {_{A}\mathcal{M}_{A}} \ar[u]^{\mathbb{H}} \ar[r]_{F} & {_{A'}\mathcal{M}'_{A'}} \ar[u]_{\mathbb{H}'}  }
\end{equation*}
where $(\mathbb{T},\mathbb{H})$ is the adjunction $({_{A}\mathbb{T}_{A}},{%
_{A}\mathbb{H}_{A}})$ defined in \ref{5}, and
$(\T^{\prime},\H^{\prime})$ is analogously defined. \newline We
have that:

$P\in {_{A}}${$\mathcal{M}$}${_{A}}$ is $\mathcal{E}_{\mathbb{H}}$%
-projective $\Longrightarrow F(P)\in {_{A^{\prime }}\mathcal{M^{\prime }}%
_{A^{\prime }}}$ is $\mathcal{E}_{\mathbb{H^{\prime
}}}$-projective; the converse is true whenever $F$ is separable.
\newline In particular we obtain that:

i) $A$ is separable in ${\mathcal{M}}\Longrightarrow A^{\prime }$
is separable in $\mathcal{M^{\prime }}$ (i.e. $\mathbb{H}$ is
separable $\Longrightarrow \mathbb{H^{\prime }}$ is separable);
the converse is true whenever $F$ is separable.

ii) If $F^{\prime }$ preserves kernels, then: $A$ is formally smooth in ${%
\mathcal{M}}\Longrightarrow A^{\prime }$ is formally smooth in ${\mathcal{%
M^{\prime }}}$; the converse is true whenever $F$ is separable.
\end{proposition}

\begin{proof}
Define $\alpha _{M}:F^{\prime }(A)\otimes' F^{\prime }(M)\otimes'
F^{\prime }(A)\rightarrow F^{\prime }(A\otimes M\otimes A)$ by
$\alpha _{M}=\phi _{2}(A\otimes M,A)[\phi _{2}(A,M)\otimes'
F^{\prime }(A)]$, for any $M\in \mathcal{M}$. Then $(\alpha
_{M})_{M\in \mathcal{M}}$ defines a natural equivalence $\alpha
:\mathbb{T}^{\prime }\mathbb{F^{\prime }}\rightarrow
\mathbb{F}\mathbb{T}$.\newline The first assertion holds by
Theorem \ref{teo F and P-project}. \newline
i) By Proposition \ref{pro sep=>left project}, $A$ is separable in ${%
\mathcal{M}}$ if and only if $A\in {_{A}\mathcal{M}_{A}}$ is $\mathcal{E}_{\mathbb{H}}$%
-projective if and only if the functor $\mathbb{H}$ is separable. Analogously $%
A^{\prime }$ is separable in $\mathcal{M^{\prime }}$ if and only if $A^{\prime }\in {%
_{A^{\prime }}\mathcal{M^{\prime }}_{A^{\prime }}}$ is $\mathcal{E}_{\mathbb{%
H^{\prime }}}$-projective if and only if the functor
$\mathbb{H}^{\prime }$ is separable. Since $A^{\prime }=F(A),$ we
conclude by the first part.\newline
ii) Let $(\Omega ^{1}(A),j)=\text{ker}(m_{A})$ in $\mathcal{M}$. Since $%
F^{\prime }$ preserves kernels, we get that
\begin{equation*}
(\Omega ^{1}(A^{\prime }),j^{\prime }):=\ker (m_{A^{\prime
}})=(F^{\prime }(\Omega ^{1}(A),\phi _{2}(A,A)F^{\prime }(j))
\end{equation*}%
in $\mathcal{M}^{\prime }$. Observe that, $\Omega ^{1}(A^{\prime })=\text{ker%
}(m_{A^{\prime }})=\text{ker}[{F^{\prime }(m)}{\phi _{2}(A,A)}]$.
Now, if we regard regard $\Omega ^{1}(A)$ as an $A$-bimodule via
the structures induced by $m_{A}$ and $\Omega ^{1}(A^{\prime })$
as an $A^{\prime }$-bimodule via the structures induced by
$m_{A^{\prime }}$, we obtain that $\Omega ^{1}(A^{\prime
})=F(\Omega ^{1}(A))$.\newline By definition, $A$ is formally
smooth in ${\mathcal{M}}$ if and only if $\Omega
^{1}A\in {_{A}}${$\mathcal{M}$}${_{A}}$ is $\mathcal{E}_{\mathbb{H}}$%
-projective. Analogously $A^{\prime }$ is formally smooth in $\mathcal{%
M^{\prime }}$ if and only if $\Omega ^{1}(A^{\prime })\in {_{A^{\prime }}\mathcal{%
M^{\prime }}_{A^{\prime }}}$ is $\mathcal{E}_{\mathbb{H^{\prime }}}$%
-projective. Since $\Omega ^{1}(A^{\prime })=F(\Omega ^{1}(A)),$
we conclude by the first part.
\end{proof}

\begin{examples}
\label{ex of F}Let $H$ be a Hopf algebra over a field $K.$ With hypotheses
and notations of Proposition \ref{pro separability of F gen}, let ${\mathcal{%
M}}^{\prime }:=\mathfrak{M}_{K}.$ We want to apply the previous
result to the particular case when ${\mathcal{M}}$ is either
$({^{H}\mathfrak{M}^{H}},\otimes ,K)$ or
$(\mathfrak{M}^{H},\otimes ,K)$. Let $A$ be an algebra in
$\mathcal{M}$.\newline
1) ${\mathcal{M}}:={^{H}\mathfrak{M}^{H}}$. The forgetful functor $F_{1}:{%
_{A}^{H}\mathfrak{M}_{A}^{H}}\rightarrow {_{A}\mathfrak{M}_{A}}$ has a right
adjoint $G_{1}:{_{A}\mathfrak{M}_{A}}\rightarrow {_{A}^{H}\mathfrak{M}%
_{A}^{H}}$, $G_{1}(M)=H\otimes M\otimes H$, where $G_{1}(M)$ is a bicomodule
via $\Delta _{H}\otimes M\otimes H$ and $H\otimes M\otimes \Delta _{H}$, and
it is a bimodule with diagonal actions. For any $M\in $ ${{_{A}^{H}%
\mathfrak{M}_{A}^{H}}}$ the unit of the adjunction is the map $\eta
_{M}:M\rightarrow H\otimes M\otimes H,\eta _{M}=(^{H}\rho _{M}\otimes
H)\circ \rho _{M}^{H}.$\newline
2) ${\mathcal{M}}:={\mathfrak{M}^{H}}$. The forgetful functor $F_r:{_{A}%
\mathfrak{M}_{A}^{H}}\rightarrow {_{A}\mathfrak{M}_{A}}$ has a right adjoint
$G_r:{_{A}\mathfrak{M}_{A}}\rightarrow {_{A}\mathfrak{M}_{A}^{H}}$, $%
G_r(M)=M\otimes H$, where $G_{r}(M)$ is a comodule via $M\otimes \Delta _{H}$%
, and it is a bimodule with diagonal actions. For any $M\in $ ${{_{A}%
\mathfrak{M}_{A}^{H}}}$ the unit of the adjunction is the map
$\eta _{M}:M\rightarrow M\otimes H,\eta _{M}=\rho
_{M}^{H}.$\newline In the case $A=H$ we set
$(F_{2},G_{2}):=(F_{1},G_{1})$.\newline
The forgetful functor $F_{b}:{_{H}^{H}\mathfrak{M}_{H}^{H}}\rightarrow {_{H}%
\mathfrak{M}_{H}^{H}}$ has a right adjoint $G_{b}:{_{H}\mathfrak{M}_{H}^{H}}%
\rightarrow {_{H}^{H}\mathfrak{M}_{H}^{H}}$, $G_{b}(M)=H\otimes M$, where $%
G_{b}(M)$ is a bicomodule via $\Delta _{H}\otimes M$ and $M\otimes \rho
_{M}^{H}$, and it is a bimodule with diagonal action.\newline
The forgetful functor ${F}_{a}:{{_{H}\mathfrak{M}_{H}^{H}}\rightarrow {_{H}%
\mathfrak{M}_{H}},}$ is nothing but $F_r$ in the case $A=H$. Then it has a
right adjoint $G_{a}:{_{H}\mathfrak{M}_{H}}\rightarrow {_{H}\mathfrak{M}%
_{H}^{H}}$, which is $G_r$ for $A=H$.\newline
Note that the forgetful functor $F_{2}:{_{H}^{H}\mathfrak{M}_{H}^{H}}%
\rightarrow {_{H}\mathfrak{M}_{H}}$ can be decomposed as $F_{2}=F_{a}\circ
F_b$.
\end{examples}

In view of Examples \ref{ex of F}, we obtain the following important result:

\begin{theorem}
\label{teo separability of F}Let $H$ be a Hopf algebra over a
field $K$ and
let $\mathcal{M}$ denote either ${^{H}\mathfrak{M}^{H}}$ or ${%
\mathfrak{M}^{H}}$. Let $A$ be an algebra in $\mathcal{M%
}$ and consider the forgetful functors $\mathbb{H}:{_{A}\mathcal{M}_{A}}%
\rightarrow \mathcal{M}$, $\mathbb{H^{\prime}}:{_{A}\mathfrak{M}_{A}}%
\rightarrow \mathfrak{M}_K$ and $F:{_{A}\mathcal{M}_{A}}\rightarrow {_{A}%
\mathfrak{M}_{A}}$. \newline We have that:

$P\in {_{A}\mathcal{M}_{A}}$ is $\mathcal{E}_{\mathbb{H}}$-projective $%
\Longrightarrow $ $P$ is $\mathcal{E}_{\mathbb{H^{\prime
}}}$-projective as an object in ${_{A}\mathfrak{M}_{A}}$; the
converse is true whenever $F$ is separable.\newline In particular
we obtain that:

i) $A$ is separable as an algebra in $\mathcal{M}$
$\Longrightarrow $ $A$ is separable as an algebra in
$\mathfrak{M}_{K}$; the converse is true whenever $F$ is
separable.

ii) $A$ is formally smooth as an algebra in $\mathcal{M}$
$\Longrightarrow $ $A$ is formally smooth as an algebra in
${\mathfrak{M}}_{K}$; the converse is true whenever $F$ is
separable.
\end{theorem}

\begin{proof}Apply Proposition \ref{pro separability of F gen} in
the case when $\mathcal{M}'=\mathfrak{M}_K$, and
$F':\mathcal{M}\rightarrow \mathfrak{M}_K$ is the forgetful
functor.
\end{proof}

Dually we have.

\begin{claim}
Let $(C,\Delta ,\varepsilon )$ be a coalgebra in a monoidal category $(\mathcal{M},\otimes ,%
\mathbf{1}).$ Like in the dual case, we have the functors
\begin{gather*}
{^{C}}\mathbb{H}:\mathcal{M}\rightarrow {^{C}\mathcal{M}}\text{ where }{^{C}}%
\mathbb{H}(X):=C\otimes X\text{ and }{^{C}}\mathbb{H}(f):=C\otimes f, \\
\mathbb{H}^{C}:\mathcal{M}\rightarrow \mathcal{M}^{C}\text{ where }\mathbb{H}%
^{C}(X):=X\otimes C\text{ and }\mathbb{H}^{C}(f):=f\otimes C, \\
{^{C}\mathbb{H}^{C}}:\mathcal{M}\rightarrow {^{C}\mathcal{M}^{C}}\text{
where }{^{C}\mathbb{H}^{C}}(X):=C\otimes (X\otimes C)\text{ and }{^{C}%
\mathbb{H}^{C}}(f):=C\otimes (f\otimes C),
\end{gather*}
with their left adjoint ${^{C}}\mathbb{T},\mathbb{T}^{C},{^{C}\mathbb{T}^{C}}
$, respectively, that forget the comodule structures. \label{Co-adjunctions}%
Then the adjunctions $({^{C}}\mathbb{T},{^{C}}\mathbb{H})$, $(\mathbb{T}^{C},%
\mathbb{H}^{C})$ and $({^{C}\mathbb{T}^{C}},{^{C}\mathbb{H}^{C}})$ gives
rise to the following classes:
\begin{gather*}
{^{C}}\mathcal{I} :=\mathcal{I}_{{^{C}}\mathbb{T}}=\{g\in {^{C}\mathcal{M}%
}\mid g\text{ is a cosplits in }\mathcal{M}\}, \\
\mathcal{I}^{C} :=\mathcal{I}_{\mathbb{T}^{C}}=\{g\in
\mathcal{M}^{C}\mid
g\text{ is a cosplits in }\mathcal{M}\}, \\
{^{C}\mathcal{I}^{C}} :=\mathcal{I}_{{^{C}\mathbb{T}^{C}}}=\{g\in {^{C}%
\mathcal{M}^{C}}\mid g\text{ is a cosplits in }\mathcal{M}\}.
\end{gather*}
By duality we can obtain the definition of coseparability and
formal smoothness for a coalgebra $(C,\Delta ,\varepsilon )$ in a
monoidal category $\mathcal{M}$. We say that $C$ is
\emph{coseparable} whenever the comultiplication $\Delta $
cosplits in ${^{C}\mathcal{M}^{C}.}$

Assume that $\mathcal{M}$ is a coabelian monoidal category. Then
$
\mho _{1}C:=\C\Delta _{C}
$
carries a natural $C$-bicomodule structure that makes it the
cokernel of $\Delta$ in the category
${^{C}}${$\mathcal{M}$}${^{C}}$. We say that $C$ is \emph{formally
smooth }in ${\mathcal{M}}$ if $\mho _{1}C$ is
${^{C}\mathcal{I}^{C}}$-injective. By duality, from Lemma \ref{lem
one side} and Proposition \ref{2 sep functors}, we obtain the
following two results.
\end{claim}

\begin{lemma}
\label{lem dual one side}Let $C$ be a coseparable coalgebra in a
monoidal category {$\mathcal{M}$}. The following assertions hold
true:

1) The forgetful functor ${^{C}}\mathbb{T}:{^{C}\mathcal{M}}\rightarrow {%
\mathcal{M}}$ is separable. In particular, any left $C$-comodule $(M{,{^{C}}%
\rho }_{M})$ is $^{C}\mathcal{I}$-injective and if $M$ is a $C$-bicomodule,
the comultiplication ${{^{C}}\rho }_{M}:M\rightarrow C\otimes M$ has a
retraction $^{C}\mu $ which is $C$-bicolinear and natural.

2) The forgetful functor $\mathbb{T}^{C}:\mathcal{M}^{C}\rightarrow \mathcal{%
M}$ is separable. In particular, any right $C$-comodule $(M{,\rho }_{M}^{R})$
is ${\mathcal{I}^{C}}$-injective and if $M$ is a $C$-bicomodule, the
comultiplication ${\rho }_{M}^{C}:M\rightarrow M\otimes C$ has a retraction $%
\mu ^{R}$ which is $C$-bicolinear and natural.
\end{lemma}

\begin{proposition}
\label{2 sep functors dual}Let $H$ be a Hopf algebra with antipode $S$ over
a field $K$. The forgetful functors $\mathfrak{M}^H_H\rightarrow \mathfrak{M}%
_H$ and $_H\mathfrak{M}^H_H\rightarrow {_H\mathfrak{M}_H}$ are separable.
\end{proposition}

\begin{proof} It is dual to Proposition \ref{2 sep functors}.
\end{proof}

\begin{remark}
By Lemma \ref{lem dual one side}, the forgetful functor $\mathbb{T}^{C}:{%
\mathcal{M}^{C}}\rightarrow \mathcal{M}$ is separable for any
coseparable coalgebra $C$ in a monoidal category $\mathcal{M}$.
The converse does not hold true. In fact, in the case when
$\mathcal{M}=\mathfrak{M}_H$ and $C=H$, the functor
$\mathbb{T}^{C}$ is always separable (Lemma \ref{2 sep functors
dual}), but $C$ is coseparable in $\mathcal{M}$ if and only if $H$
is a cosemisimple coalgebra (\cite[Proposition
2.11]{AMS-Spliting}).
\end{remark}

\begin{proposition}
\label{pro cosep=>left inject}Let $C$ be a coalgebra in a monoidal category {%
$\mathcal{M}$. }The following assertions are equivalent:

$(a)$ $C$ is coseparable in $\mathcal{M}$.

$(b)$ The forgetful functor ${^{C}\mathbb{T}^{C}}:{^{C}\mathcal{M}^{C}}%
\rightarrow \mathcal{M}$ is separable.

$(c)$ Any $C$-bicomodule is ${^{C}\mathcal{I}^{C}}$-injective.

$(d)$ The $C$-bicomodule $C$ is ${^{C}\mathcal{I}^{C}}$-injective.
\end{proposition}

\begin{corollary}
Any coseparable coalgebra in a coabelian monoidal category
$\mathcal{M}$ is formally smooth.
\end{corollary}

\begin{corollary}
\label{coro inj}Let $C$ be a coseparable coalgebra in
$\mathfrak{M}_{K}$. Then any left $C$-comodule is injective in
$^{C}\mathfrak{M}$. Hence any left $C$-comodule is also projective
in $^{C}\mathfrak{M}$ and $C$ is
cosemisimple. Moreover any $C$-bicomodule is injective in $^{C}\mathfrak{M}%
^{C}$ and hence any $C$-bicomodule is projective in $^{C}\mathfrak{M}^{C}$.
\end{corollary}

\begin{claim}
\label{ex: co monoidal functor} Let $(F^{\prime },\phi _{0},\phi _{2}):(%
\mathcal{M},\otimes ,\mathbf{1)\rightarrow (}\mathcal{M}^{\prime }\mathfrak{,%
}\otimes ^{\prime },\mathbf{1}^{\prime }\mathbf{)}$ be a monoidal
functor between two monoidal categories, where $\phi
_{2}(U,V):F^{\prime }(U\otimes
V)\rightarrow F^{\prime }(U)\otimes ^{\prime }F^{\prime }(V),$ for any $%
U,V\in \mathcal{M}$, and $\phi _{0}:\mathbf{1}^{\prime
}\mathbf{\rightarrow
F^{\prime }(1)}.$ Let $(C,\Delta ,\varepsilon )$ is a coalgebra in $\mathcal{%
M}$. It is well known that $(F^{\prime }(C),\Delta _{F^{\prime
}(C)},\varepsilon _{F^{\prime }(C)})$ is a coalgebra in
$\mathcal{M}^{\prime }$, where
\begin{equation*}
\begin{tabular}{l}
$\Delta _{F^{\prime }(C)}:=F^{\prime }(C)\overset{F^{\prime }(\Delta )}{%
\longrightarrow }F^{\prime }(C\otimes C)\overset{\phi _{2}^{-1}(C,C)}{%
\longrightarrow }F^{\prime }(C)\otimes ^{\prime }F^{\prime }(C)$ \\
$\varepsilon _{F^{\prime }(C)}:=F^{\prime }(C)\overset{F^{\prime
}(\varepsilon )}{\longrightarrow }F^{\prime
}(\mathbf{1})\overset{\phi
_{0}^{-1}}{\longrightarrow }\mathbf{1}^{\prime }.$%
\end{tabular}%
\end{equation*}%
Consider the functor $F:{^{C}\mathcal{M}^{C}}\rightarrow {^{C^{\prime }}%
\mathcal{M}^{\prime }{}^{C^{\prime }}}$ defined by {\small
$F((M,{^{C}\rho
_{M}},{\rho _{M}^{C}}))=(F^{\prime }(M),{^{C^{\prime }}\rho _{F^{\prime }(M)}%
},{\rho _{F^{\prime }(M)}^{C^{\prime }}}),$} where
\begin{equation*}
\begin{tabular}{l}
${^{C^{\prime }}\rho _{F^{\prime }(M)}}:=F^{\prime }(C)\overset{F^{\prime }({%
^{C}\rho _{M}})}{\longrightarrow }F^{\prime }(C\otimes
M)\overset{\phi _{2}^{-1}(C,M)}{\longrightarrow }F^{\prime
}(C)\otimes ^{\prime }F^{\prime
}(M)$ \\
${\rho _{F^{\prime }(M)}^{C^{\prime }}}:=F^{\prime }(C)\overset{F^{\prime }({%
\rho _{M}^{C}})}{\longrightarrow }F^{\prime }(M\otimes
C)\overset{\phi _{2}^{-1}(M,C)}{\longrightarrow }F^{\prime
}(M)\otimes ^{\prime }F^{\prime
}(C)$.%
\end{tabular}%
\end{equation*}
\end{claim}

\begin{proposition}
\label{pro coseparability of F gen}Let ${\mathcal{M}}$ and $\mathcal{%
M^{\prime }}$ be coabelian monoidal categories. Let $C$, $C^{\prime }$, $%
F^{\prime }$ and $F$ as in Example \ref{ex: co monoidal functor}.
Then, in
the following diagrams, ${\mathbb{H}^{\prime }}\circ G^{\prime }$ and $%
G\circ {\mathbb{H}}$ are naturally equivalent and $G^{\prime }\circ {\mathbb{%
T}}={\mathbb{T}^{\prime }}\circ G$:
\begin{equation*}
\xymatrix@R=15pt@C=50pt{
  {^C\mathcal{M}^C} \ar[d]_{\mathbb{T}} \ar[r]^{G} & {^{C'}\mathcal{M'}^{C'}} \ar[d]^{\mathbb{T}'} \\
  \mathcal{M}  \ar[r]_{G'} & \mathcal{M'}   }
\text{ \qquad } \xymatrix@R=15pt@C=50pt{
  {^C\mathcal{M}^C} \ar[r]^{G} & {^{C'}\mathcal{M'}^{C'}}  \\
  \mathcal{M}  \ar[u]^{\mathbb{H}} \ar[r]_{G'} & \mathcal{M'} \ar[u]_{\mathbb{H}'}  }
\end{equation*}
where $(\mathbb{T},\mathbb{H})$ is the adjunction $({^{C}\mathbb{\mathbb{T}}%
^{C}},{^{C}\mathbb{\mathbb{H}}^{C}})$ defined in \ref{Co-adjunctions}, and $(%
\mathbb{T^{\prime }},\mathbb{H^{\prime }})$ is analogously
defined.  \newline We have that:

$I\in {^{C}\mathcal{M}^{C}}$ is $\mathcal{I}_{\mathbb{T}}$-injective $%
\Longrightarrow G(I)\in {^{C^{\prime }}\mathcal{M^{\prime
}}^{C^{\prime }}}$ is $\mathcal{I}_{\mathbb{T^{\prime
}}}$-injective; the converse is true whenever $G$ is separable.
\newline In particular we obtain that:

i) $C$ is coseparable in ${\mathcal{M}}$ $\Longrightarrow $
$C^{\prime }$ is
coseparable in $\mathcal{M^{\prime }}$ (i.e. $\mathbb{T}$ is separable $%
\Longrightarrow $ $\mathbb{T^{\prime }}$ is separable); the
converse is true whenever $G$ is separable.

ii) If $G^{\prime }$ preserves cokernels, then: $C$ is formally smooth in ${%
\mathcal{M}}$ $\Longrightarrow $ $C^{\prime }$ is formally smooth in ${%
\mathcal{M^{\prime }}}$; the converse is true whenever $G$ is
separable.
\end{proposition}

\begin{proof}
It is dual to Proposition \ref{pro separability of F gen}.
\end{proof}

\begin{examples}
\label{ex of co F}Let $H$ be a Hopf algebra over a field $K.$ With
hypotheses and notations of Proposition \ref{pro coseparability of F gen},
let ${\mathcal{M}}^{\prime }:=\mathfrak{M}_{K}.$ We want to apply the
previous result to the particular case when ${\mathcal{M}}$ is ether $({_{H}%
\mathfrak{M}_{H}},\otimes ,K)$ or $(\mathfrak{M}_{H},\otimes ,K)$.
Let $C$ be a coalgebra in $\mathcal{M}$.\newline
1) ${\mathcal{M}}:={_{H}\mathfrak{M}_{H}}$. The forgetful functor $G^{1}:{%
_{H}^{C}\mathfrak{M}_{H}^{C}}\rightarrow {^{C}\mathfrak{M}^{C}}$ has a left
adjoint $F^{1}:{^{C}\mathfrak{M}^{C}}\rightarrow {_{H}^{C}\mathfrak{M}%
_{H}^{C}}$, $F^{1}(M)=H\otimes M\otimes H$, where $F^{1}(M)$ is a bimodule
via $m_{H}\otimes M\otimes H$ and $H\otimes M\otimes m_{H}$, and it is a
bicomodule with codiagonal coactions. For any $M\in $ ${_{H}^{C}\mathfrak{M}%
_{H}^{C}}$ the counit of the adjunction is the map $\varepsilon
_{M}:H\otimes M\otimes H\rightarrow M,\varepsilon _{M}=\mu _{M}^{H}\circ
(^{H}\mu _{M}\otimes H).$\newline
2) ${\mathcal{M}}:=\mathfrak{M}_{H}$. The forgetful functor $G^r:{^{C}%
\mathfrak{M}_{H}^{C}}\rightarrow {^{C}\mathfrak{M}^{C}}$ has a left adjoint $%
F^r:{^{C}\mathfrak{M}^{C}}\rightarrow {^{C}\mathfrak{M}_{H}^{C}}$, $%
F^r(M)=M\otimes H$, where $F^r(M)$ is a module via $M\otimes m_{H}$, and it
is a bicomodule with codiagonal coactions. For any $M\in $ ${{^{C}%
\mathfrak{M}_{H}^{C}}}$ the counit of the adjunction is the map
$\varepsilon _{M}:M\otimes H\rightarrow M,\varepsilon _{M}=\mu
_{M}^{H}$. \newline In the case $C=H$ we set
$(F^{2},G^{2}):=(F^{1},G^{1}).$\newline
The forgetful functor $G^{a}:{^{H}\mathfrak{M}_{H}^{H}}\rightarrow {^{H}%
\mathfrak{M}^{H}}$ is nothing but $G^r$ in the case $C=H$. Then it has a
left adjoint $F^{a}:{^{H}\mathfrak{M}^{H}}\rightarrow {^{H}\mathfrak{M}%
_{H}^{H}}$, which is $F^r$ for $C=H$.\newline
The forgetful functor $G^{b}:{_{H}^{H}\mathfrak{M}_{H}^{H}}\rightarrow {^{H}%
\mathfrak{M}_{H}^{H}}$ has a left adjoint $F^{b}:{^{H}\mathfrak{M}_{H}^{H}}%
\rightarrow {_{H}^{H}\mathfrak{M}_{H}^{H}}$, $F^{b}(M)=H\otimes M$, where $%
F^{b}(M)$ is a bimodule via $m_{H}\otimes M$ and $H\otimes \mu _{M}^{H}$,
and it is a bicomodule with codiagonal coactions.\newline
Note that the forgetful functor $G^{2}:{_{H}^{H}\mathfrak{M}_{H}^{H}}%
\rightarrow {^{H}\mathfrak{M}^{H}}$ can be decomposed as $G^{2}=G^{a}\circ
G^{b}.$
\end{examples}

In view of Examples \ref{ex of co F}, we obtain the following important
result:

\begin{theorem}
\label{teo coseparability of F}Let $H$ be a Hopf algebra over a
field $K$
and let $\mathcal{M}$ denote either $_{H}\mathfrak{M}_{H}$ or $%
\mathfrak{M}_{H}$. Let $C$ be a coalgebra in $\mathcal{M}$
and consider the forgetful functors $\mathbb{T}:{^{C}\mathcal{M}^{C}}%
\rightarrow \mathcal{M}$, $\mathbb{T^{\prime }}:{^{C}\mathfrak{M}^{C}}%
\rightarrow \mathfrak{M}_{K}$ and $G:{^{C}\mathcal{M}^{C}}\rightarrow {^{C}%
\mathfrak{M}^{C}}$. \newline We have that:

$I\in {^{C}\mathcal{M}^{C}}$ is $\mathcal{I}_{\mathbb{T}}$-injective $%
\Longrightarrow $ $I$ is $\mathcal{I}_{\mathbb{T^{\prime
}}}$-injective as an object in ${^{C}\mathfrak{M}^{C}}$; the
converse is true whenever $G$ is separable.\newline In particular
we obtain that:

i) $C$ is coseparable as a coalgebra in $\mathcal{M}$
$\Longrightarrow $ $C$ is coseparable as a coalgebra in
${\mathfrak{M}}_{K}$; the converse is true whenever $G$ is
separable.

ii) $C$ is formally smooth as a coalgebra in $\mathcal{M}$
$\Longrightarrow $ $C$ is formally smooth as a coalgebra in
${\mathfrak{M}}_{K}$; the converse is true whenever $G$ is
separable.
\end{theorem}

\section{Ad-invariant integrals through separable functors}\label{sec: Ad-invariant integrals through separable
functors}

\begin{remark}
Let $H$ be a Hopf algebra over a field $K$. For sake of brevity
many results will be stated only for the category
${\mathfrak{M}}_{H}$. Clearly all the results still hold true for
$_{H}{\mathfrak{M}}$ (as $_{H}{\mathfrak{M}}\simeq
{\mathfrak{M}}_{H^{op}}$). Similar arguments apply to the
categories ${\mathfrak{M}}^{H}$ and $^{H}{\mathfrak{M}}$.
\end{remark}

\begin{claim}
\label{adj actions}Let $H$ be a Hopf algebra with antipode $S$ over a field $%
K$ and set:
\begin{gather*}
h\vartriangleright x:=h_{1}xS(h_{2})\qquad \text{and\qquad }%
x\vartriangleleft h:=S(h_{1})xh_{2} \\
{^{H}\varrho }(h):=h_{1}S(h_{3})\otimes h_{2}\qquad \text{and\qquad }{%
\varrho }^{H}(h):=h_{2}\otimes S(h_{1})h_{3}
\end{gather*}
for all $h,x\in H.$ It is easy to check that $\vartriangleright $ defines a
left module action of $H$ on itself called \emph{left adjoint action} and
that ${^{H}\varrho }$ defines a left comodule coaction of $H$ on itself
called \emph{left adjoint coaction}.\emph{\ }Analogously $\vartriangleleft $
gives rise to the right adjoint action and ${\varrho }^{H}$ to the right
adjoint coaction. \newline
If $S$ is bijective, we can consider the following actions and coactions of $%
H$ on itself:
\begin{gather*}
h\blacktriangleright x:=h_{2}xS^{-1}(h_{1})\qquad \text{and } \qquad%
x\blacktriangleleft h:=S^{-1}(h_{2})xh_{1} \\
\overline{\varrho }^{H}(h)=h_{2}\otimes h_{3}S^{-1}(h_{1})\qquad \text{%
and } \qquad{^{H}}\overline{{\varrho
}}(h):=S^{-1}(h_{3})h_{1}\otimes h_{2}.
\end{gather*}
The structures above provide two different ways of looking at $H$ as an
object in the categories of Yetter-Drinfeld modules. In fact, if ${%
\Delta }_{H}$ is the comultiplication and $m_{H}$ is the
multiplication
of $H$, then $H$ can be regarded as an object in ${_{H}^{H}\mathcal{YD}}$, ${%
\mathcal{YD}_{H}^{H}}$, ${_{H}\mathcal{YD}^{H}}$, ${^{H}\mathcal{YD}_{H}}$
respectively via:
\begin{equation*}
(\vartriangleright ,\Delta _{H}),(\vartriangleleft ,{\Delta _{H}}%
),(\blacktriangleright ,{\Delta _{H}}),(\blacktriangleleft
,{\Delta
_{H}})\text{\qquad or}\qquad ({m}_{H},{^{H}\varrho }),({m}_{H},{\varrho }%
^{H}),({m}_{H},{\overline{\varrho }^{H}}),({m}_{H},{^{H}}\overline{{\varrho }%
}).
\end{equation*}
\end{claim}

\begin{claim}
\label{ex. of adjoint}\textbf{The adjunctions.}\newline The
actions recalled in \ref{adj actions} are closely linked to the
categories of Yetter-Drinfeld modules. We now consider some
adjunctions involving these modules that will be very useful in
finding equivalent conditions to the existence of an
$ad$-invariant integral.\newline $1)$ The forgetful functor
$F_{3}:{_{H}^{H}\mathcal{YD}}\rightarrow {_{H} \mathfrak{M}}$ has
a right adjoint $G_{3}:{_{H}\mathfrak{M}}\rightarrow
{_{H}^{H}\mathcal{YD}},{G}(M)=H\otimes M$, where $G(M)$ is a
comodule via $\Delta _{H}\otimes M$ and a module via the action:
$h\cdot (l\otimes m)=h_{1}lS(h_{3})\otimes h_{2}m.$ For any $M\in
$ ${_{H}^{H}\mathcal{YD}}$ the unit of the adjunction is the map
$\eta _{M}:M\rightarrow H\otimes M,\eta _{M}={^{H}\rho }_{M}.$
\newline
$2)$ The forgetful functor
$F_{4}:{\mathcal{YD}_{H}^{H}}\rightarrow \mathfrak{M}{_{H}}$ has a
right adjoint $G_{4}:\mathfrak{M}{_{H}}\rightarrow
{\mathcal{YD}_{H}^{H}},{G}_{4}(M)=M\otimes H$, where $G_{4}(M)$ is
a comodule via $M\otimes \Delta _{H}$ and a module via the action:
$(m\otimes l)\cdot h=mh_{2}\otimes S(h_{1})lh_{3}.$ For any $M\in
$ ${\mathcal{YD}_{H}^{H}}$ the unit of the adjunction is the map
$\eta _{M}:M\rightarrow M\otimes H,\eta _{M}={\rho
_{M}^{H}}.$\newline $3)$ Assume $H$ has bijective antipode. The
forgetful functor $F_{5}:{_{H}\mathcal{YD}^{H}}\rightarrow
{_{H}\mathfrak{M}}$ has a right adjoint
$G_{5}:{_{H}\mathfrak{M}}\rightarrow
{_{H}\mathcal{YD}^{H}},{G}_{5}(M)=M\otimes H$, where $G_{5}(M)$ is
a comodule via $M\otimes \Delta _{H}$ and a module via the action:
$h\cdot (l\otimes m)=h_{2}l\otimes h_{3}mS^{-1}(h_{1}).$ For any
$M\in $ ${_{H}\mathcal{YD}^{H}}$ the unit of the adjunction is the
map $\eta _{M}:M\rightarrow M\otimes H,\eta _{M}={\rho
_{M}^{H}}.$\newline $4)$ Assume $H$ has bijective antipode. The
forgetful functor $F_{6}:{^{H}\mathcal{YD}_{H}}\rightarrow
{\mathfrak{M}_{H}}$ has a right adjoint
$G_{6}:{\mathfrak{M}_{H}}\rightarrow
{^{H}\mathcal{YD}_{H}},{G}_{6}(M)=H\otimes M$, where $G_{6}(M)$ is
a comodule via $\Delta _{H}\otimes M$ and a module via the action:
$(l\otimes m)\cdot h=S^{-1}(h_{3})lh_{1}\otimes mh_{2}.$ For any
$M\in $ ${^{H}\mathcal{YD}_{H}}$ the unit of the adjunction is the
map $\eta _{M}:M\rightarrow H\otimes M,\eta _{M}={^{H}\rho }_{M}.$
\medskip\\ Consider now the dual version of this
functors.\medskip\\ $1^{op})$ The forgetful functor
$G^{3}:{_{H}^{H}\mathcal{YD}}\rightarrow {^{H}\mathfrak{M}}$ has a
left adjoint $F^{3}:{^{H}\mathfrak{M}}\rightarrow
{_{H}^{H}\mathcal{YD}},{F}^{3}(M)=H\otimes M$, where $F^3(M)$ is a
module via $m_{H}\otimes M$ and a comodule via the coaction:
${^{H}\rho }(h\otimes m)=h_{1}m_{-1}S(h_{3})\otimes h_{2}\otimes
m_{0}$. For any $M\in $ ${_{H}^{H}\mathcal{YD}}$ the counit of the
adjunction is the map $\varepsilon _{M}:H\otimes M\rightarrow
M,\varepsilon _{M}={^{H}\mu _{M}}.$
\newline
$2^{op})$ The forgetful functor
$G^{4}:{\mathcal{YD}_{H}^{H}}\rightarrow {\mathfrak{M}^{H}} $ has
a left adjoint $F^{4}:{\mathfrak{M}^{H}}\rightarrow
{\mathcal{YD}_{H}^{H}},{F}^{4}(M)=M\otimes H$, where $F^{4}(M)$ is
a module via $M\otimes m_{H}$ and a comodule via the coaction:
${\rho }^{H}(m\otimes h)=m_{0}\otimes h_{2}\otimes
S(h_{1})m_{1}h_{3}.$ For any $M\in $ ${\mathcal{YD}_{H}^{H}}$ the
counit of the adjunction is the map $\varepsilon _{M}^{4}:M\otimes
H\rightarrow M,\varepsilon _{M}^{4}={\mu _{M}^{H}}.$
\newline
$3^{op})$ Assume $H$ has bijective antipode. The forgetful functor
$G^{5}:{_{H}\mathcal{YD}^{H}}\rightarrow {\mathfrak{M}^{H}}$ has a
left adjoint $F^{5}:{\mathfrak{M}^{H}}\rightarrow
{_{H}\mathcal{YD}^{H}},{F}^{5}(M)=H\otimes M$, where $F^{5}(M)$ is
a module via $m_{H}\otimes M$ and a comodule via the coaction:
${\rho }^{H}(h\otimes m)=h_{2}\otimes m_{0}\otimes
h_{3}m_{1}S^{-1}(h_{1})$. For any $M\in $ ${_{H}\mathcal{YD}^{H}}$
the counit of the adjunction is the map $\varepsilon _{M}:H\otimes
M\rightarrow M,\varepsilon _{M}={^{H}\mu _{M}}.$\newline $4^{op})$
Assume $H$ has bijective antipode. The forgetful functor
$G^{6}:{^{H}\mathcal{YD}_{H}}\rightarrow {^{H}\mathfrak{M}}$ has a
left adjoint $F^{6}:{^{H}\mathfrak{M}}\rightarrow
{^{H}\mathcal{YD}_{H}},F^{6}(M)=M\otimes H$, where $F^{6}(M)$ is a
module via $M\otimes m_{H}$ and a comodule via the coaction:
${^{H}\rho }(m\otimes h)=S^{-1}(h_{3})m_{-1}h_{1}\otimes
m_{0}\otimes h_{2}$. For any $M\in $ ${^{H}\mathcal{YD}_{H}}$ the
counit of the adjunction is the map $\varepsilon _{M}:M\otimes
H\rightarrow M,\varepsilon _{M}={\mu _{M}^{H}}.$
\end{claim}

\begin{claim}
\label{claim integrals}Let $K$ be any field. An augmented
$K$-algebra $(A,m,u,p)$ is a $K$-algebra $(A,m,u)$ endowed with an
algebra homomorphism $p:A\to K$ called augmentation of $A$. An
element $x\in A$ is a \emph{left integral in }$A$, whenever
$a\cdot _{A}x=p \left( a\right) x$, for every $a\in A$. The
definition of a right integral in $A$ is analogous. $A$ is called
\emph{unimodular}, whenever the space of left and right integrals
in $A$ coincide. A (left or right) integral $x$ in $A$ is called a
\emph{total integral in} $ A $, whenever $p \left( x\right)
=1_{K}$.\medskip\newline Let $\left( H,m_{H},u_{H},\Delta
_{H},\varepsilon _{H}\right) $ be a bialgebra.\newline 1) $\left(
H,m_{H},u_{H},\varepsilon _{H}\right) $ is an augmented algebra.
Then a left integral in $H$ is an element $t\in H$ such that
$h\cdot _{H}t=\varepsilon _{H}\left( h\right) t$, for every $h\in
H$. Moreover $t$ is total whenever $\varepsilon _{H}\left(
t\right) =1_{K}. $\newline 2) $\left( H^{\ast },m_{H^{\ast
}},u_{H^{\ast }},\varepsilon _{H^{\ast }}\right) $ is an augmented
algebra. Then a left integral in $H^{\ast }$ is an element
$\lambda \in H^{\ast },$ that is a $K$-linear map $ f \lambda
=f\left( 1_{H}\right) \lambda$, for every $f\in H^{\ast }$.
Moreover $\lambda $ is total, whenever $\lambda \left(
1_{H}\right) =1_{K}.$ It is clear that $\lambda \in H^{\ast }$ is
a left (resp. right) integral in $H^{\ast }$ if and only if
$h_{1}\lambda (h_{2}) =1_{H}\lambda (h) $ (resp. $\lambda (h_{1})
h_{2}=\lambda (h) 1_{H}$) for every $h\in H$.\newline If $H$ is
finite dimensional, $H^{\ast }$ becomes a Hopf algebra: in
particular one can consider the notion of left integral in
$(H^{\ast  })^{\ast}$ in the sense of 2). By means of the
isomorphism \begin{equation*} H\rightarrow H^{\ast \ast
}:h\longmapsto \left(
\begin{tabular}{l}
$H^{\ast }\rightarrow K$ \\
$f\longmapsto f\left( h\right) $%
\end{tabular}%
\right),
\end{equation*}%
one can check that a left integral in $H^{\ast \ast }$ is nothing
but a left integral in $H$ in the sense of 1): thus there is no
danger of confusion.
\end{claim}

For the reader's sake, we outline the following facts.

\begin{theorem}
\label{teo integrals}Let $H$ be a Hopf algebra with antipode $S$ over any
field $K.$ Then we have:

1) There exists a total integral $t\in H$ (i.e. $H$ is semisimple)
if and only if $H$ is separable.

2) There exists a total integral $\lambda\in H^{\ast}$ (i.e. $H$
is cosemisimple) if and only if $H$ is coseparable.

\end{theorem}

\begin{proof}
1) $"\Leftarrow "$ Let $\sigma :H\rightarrow H\otimes H\ $an
$H$-bilinear section of the multiplication $m$ and set $t_{\sigma
}:=(H\otimes \varepsilon _{H}) \sigma (1_{H}) \in H.$ Then
$t_{\sigma }$ is a total integral.\newline $"\Rightarrow "$ Let
$t\in H$ be a total integral. Since $t$ is a left integral and
$\Delta _{H}$ is an homomorphism of algebras, we have:
\begin{equation}
ht_{1}\otimes S(t_{2}) =h_{1}t_{1}\otimes S(h_{2}t_{2}) h_{3}=\varepsilon
_{H}(h_{1}) t_{1}\otimes S(t_{2}) h_{2}=t_{1}\otimes S(t_{2}) h,\forall h\in
H,  \label{left integ
1}
\end{equation}
so that the map $\sigma _{t}:H\rightarrow H\otimes H:h\mapsto ht_{1}\otimes
S(t_{2}) $ is $H$-bilinear. Moreover $m_{H}\sigma _{t}(h) =ht_{1}S(t_{2})
=h\varepsilon _{H}(t) =h,$ so that $\sigma _t$ is an $H$-bilinear section of $%
m_{H}$ and $H$ is separable by definition.\newline
2) $"\Leftarrow "$ Let $\theta :H\otimes H\rightarrow H\ $an $H$%
-bicolinear retraction of the comultiplication $\Delta$ and set $%
\lambda _{\theta }:=\varepsilon _{H}\theta (-\otimes 1_{H}) \in H^{\ast }.$
Then $\lambda _{\theta }$ is a total integral.\newline
$"\Rightarrow "$ Let integral $\lambda \in H^{\ast }$ be a left integral
such that $\lambda (1_{H}) =1$. Since $\lambda $ is a left integral and $m$
is an homomorphism of coalgebras, we have:
\begin{equation}
x_{1}\lambda (x_{2}S(y))=x_{1}S(y_{2}) \lambda
(x_{2}S(y_{1}))y_{3}=(xS(y_{1}) ) _{1}\lambda ((xS(y_{1}))
_{2})y_{2}=\lambda (xS(y_{1}))y_{2},\forall x,y\in H,  \label{left
integ}
\end{equation}
so that the map $\theta _{\lambda }:H\otimes H\rightarrow
H:x\otimes y\mapsto x_{1}\lambda (x_{2}S(y))$ is $H$-bicolinear.
Moreover $\theta _{\lambda }\Delta (h) =h_{1}\lambda
(h_{2}S(h_{3}))=h\lambda (1_{H}) =h,$ so that $\theta _{\lambda }$
is an $H$-bicolinear retraction of the comultiplication $\Delta$
and $H$ is coseparable by definition.
\end{proof}

Our next aim is to characterize the existence of a so-called
$ad$-invariant integral.\newline A remarkable fact is that any
semisimple and cosemisimple Hopf algebra $H$ over a field $K$
admits such an integral (see \cite[Theorem 2.27]{AMS-Spliting}).

\begin{definition}\cite[Definition 1.11]{DO}
\label{ad invariant integ}Let $H$ be a Hopf algebra with antipode
$S$ over any field $K$ and let $\lambda \in H^{\ast }$. $\lambda $
will be called an $\emph{ad}$\emph{-invariant} \emph{integral}
whenever:

$a$) $h_{1}\lambda (h_{2})=1_{H}\lambda (h)$ for all $h\in H$ (i.e. $\lambda
$ is a left integral in $H^{\ast }$);

$b$) $\lambda (h_{1}xS(h_{2}))=\varepsilon (h)\lambda (x),$ for all $h,x\in
H $ (i.e. $\lambda $ is left linear with respect to $\vartriangleright $)$;$

$c$) $\lambda (1_{H})=1_{K}$.
\end{definition}

\begin{lemma}
\label{ad-inv = retraction} An element $\lambda \in H^{\ast }$ is an $ad$%
-invariant integral if and only if it is a retraction of the unit $%
u_{H}:K\to H$ of $H$ in ${_{H}^{H}\mathcal{YD,}}$ where $H$ is regarded as
an object in the category via the left adjoint action $\vartriangleright $
and the comultiplication $\Delta _{H}$.
\end{lemma}

\begin{example}
\label{ex ad inv of KG} 1) Let $G$ be an arbitrary group an let
${K}G$ {be the
group algebra associated. Let }$\lambda :KG\rightarrow K$ be defined by $%
\lambda (g)=\delta _{e,g}$ (the Kronecker symbol), where $e$
denotes the
neutral element of $G.$ Then $\lambda $ is an $ad$-invariant integral for $%
KG $ (see \cite[Corollary 2.8]{DO}).\\2) Every commutative
cosemisimple Hopf algebra has an $ad$-invariant integral.
\end{example}

\begin{remark}
If $H$ is a Hopf algebra with a nonzero integral then the left and
right integral spaces are both one-dimensional \cite[Theorem
5.4.2]{DNR}. If $H$ has a total integral $\lambda \in H^*$ (i.e.
$H$ is cosemisimple), then the left and right integral space
coincide \cite[Exercise 5.5.10]{DNR}, and are generated by
$\lambda$. Hence there can be at most one ad-invariant integral,
namely the unique total integral.
\end{remark}

The following lemma shows that in the definition of $ad$-invariant integral
we can choose $\vartriangleleft ,\blacktriangleright$ or $\blacktriangleleft
$ instead of $\vartriangleright .$ Since $\lambda $ is in particular a total
integral, it is both a left and a right integral. Thus it is the same to have a retraction of $u_{H}$ in ${%
_{H}^{H}\mathcal{YD}},{\mathcal{YD}_{H}^{H}},{_{H}\mathcal{YD}^{H}}$ or ${%
^{H}\mathcal{YD}_{H}}$.

\begin{lemma}
\label{lem ad-invariant}Let $H$ be a Hopf algebra with antipode $S$ over any
field $K$ and let $\lambda \in H^{\ast }$ be a total integral. Then the
following are equivalent:

$(1) $ $\lambda $ is left linear with respect to $\vartriangleright$.

$(2) $ $\lambda $ is right linear with respect to $\vartriangleleft$.

$(3) $ $\lambda $ is left linear with respect to $\blacktriangleright$.

$(4) $ $\lambda $ is right linear with respect to $\blacktriangleleft$.
\end{lemma}

\begin{proof}
We have that $\lambda $ is both a left integral and a right
integral for $H^{\ast }.$\newline Since $\lambda $ is a total
integral $S$ is bijective (see \cite[Corollary 5.4.6]{DNR}) and
hence it makes sense to consider $S^{-1}.$ \newline
$(1)\Rightarrow (2)$ Observe that: $ S(x\vartriangleleft
h)=S(S(h_{1})xh_{2})=S(h)_{1}S(x)S[ S(h)_{2}]
=S(h)\vartriangleright S(x)$.\newline Thus, since $\lambda
=\lambda S$ and $\lambda $ is left linear with respect to
$\vartriangleright ,$ we get $\lambda (x\vartriangleleft
h)=\lambda S(x\vartriangleleft h)=\lambda (S(h)\vartriangleright
S(x))=$ $\varepsilon S(h)\lambda (S(x))=\varepsilon (h)\lambda
S(x)=\varepsilon (h)\lambda (x)$ that is $\lambda $ is right
linear with respect to $\vartriangleleft .$
\newline
$(2)\Rightarrow (1)$ follows analogously once proved the relation $%
S(h\vartriangleright x)=S(x)\vartriangleleft S(h).$ \newline
$(1)\Rightarrow (3)$ We have: $S[ h\blacktriangleright S^{-1}(x)]
=S[ h_{2}S^{-1}(x)S^{-1}(h_{1})]
=h_{1}xS(h_{2})=h\vartriangleright x$.\newline Then, since
$\lambda =\lambda S\ $and $\lambda $ is left linear with respect
to $\vartriangleright $, we have $\lambda (h\blacktriangleright
x)=\lambda S(h\blacktriangleright S^{-1}S(x))=\lambda
(h\vartriangleright S(x))=\varepsilon (h)\lambda S(x)=\allowbreak
\varepsilon (h)\lambda (x)$ i.e. $\lambda $ is left linear with
respect to $\blacktriangleright $.
\newline
$(3)\Rightarrow (1)$ Since $\lambda $ is left linear with respect to $%
\blacktriangleright $ one has $\lambda (h\vartriangleright x)=\lambda S[
h\blacktriangleright S^{-1}(x)] =\lambda [ h\blacktriangleright S^{-1}(x)]
=\varepsilon (h)\lambda SS^{-1}(x)=\varepsilon (h)\lambda (x)$ i.e. $\lambda
$ is left linear with respect to $\vartriangleright .$\newline
$(1)\Leftrightarrow (4)$ Analogous to $(1)\Leftrightarrow (3)$ by means of $%
S^{-1}[ S(x)\vartriangleleft h] =x\blacktriangleleft h.$
\end{proof}

The following result improves \cite[Theorem 2.29]{AMS-Spliting}.

\begin{lemma}
\label{lem ad-inv}Let $H$ be a Hopf algebra with antipode $S$ over a field $%
K.$ Assume there exists an $ad$-invariant integral $\lambda \in H^{\ast }.$
Then we have that:

$i)$ The forgetful functor ${_{A}\mathfrak{M}_{A}^{H}}\rightarrow {{_{A}%
\mathfrak{M}_{A}}}$ is separable for any algebra $A$ in ${\mathfrak{M}^{H}}$.

$ii)$ The forgetful functor ${_{A}^{H}\mathfrak{M}_{A}}\rightarrow {{_{A}%
\mathfrak{M}_{A}}}$ is separable for any algebra $A$ in ${^{H}\mathfrak{M}}$.

$iii)$ The forgetful functor ${_{A}^{H}\mathfrak{M}_{A}^{H}}\rightarrow {{%
_{A}\mathfrak{M}_{A}}}$ is separable for any algebra $A$ in ${^{H}%
\mathfrak{M}^{H}}$.
\end{lemma}

\begin{proof}
$i)$ By Examples \ref{ex of F}, the forgetful functor $F_r:{_{A}%
\mathfrak{M}_{A}^{H}\rightarrow }${{$_{A}\mathfrak{M}_{A}$ has a right
adjoint }}$G_r:{{_{A}\mathfrak{M}_{A}}\rightarrow {_{A}\mathfrak{M}%
_{A}^{H},G}}_{8}(M)=M\otimes H${.} Thus by Theorem \ref{teo
Rafael}, $F_r$
is separable if and only if the unit $\eta ^{H}:\mathrm{Id}_{{_{A}%
\mathfrak{M}_{A}^{H}}}\rightarrow G_rF_r$ of the adjunction
cosplits, i.e. there exists a natural transformation $\mu
^{H}:G_rF_r\rightarrow \mathrm{Id}_{{_{A}\mathfrak{M}_{A}^{H}}}$
such that $\mu _{M}^{H}\circ \eta _{M}^{H}=\mathrm{Id}_{M}$ for
any $M$ in ${_{A}\mathfrak{M}_{A}^{H}.}$ Let us define:
\begin{equation*}
\mu _{M}^{H}:M\otimes H\rightarrow M,\quad \mu _{M}^{H}(m\otimes
h)=m_{0}\lambda (m_{1}S(h)).
\end{equation*}
Obviously $(\mu _{M}^{H})_{M\in {{_{A}\mathfrak{M}_{A}^{H}}}}$ is a
functorial morphism. \newline
Let us check that $\mu _{M}^{H}$ is a morphism in ${{_{A}\mathfrak{M}_{A}^{H}%
}}$, i.e. a morphism of $A$-bimodules and of $H$-bicomodules.
Since $\mu _{M}^{A}\in {{\mathfrak{M}^{H},}}$ we have: $\mu
_{M}^{H}((m\otimes h)a)=m_{0}a_{0}\lambda
(m_{1}a_{1}S(a_{2})S(h))=\mu _{M}^{H}(m\otimes h)a$.\newline Since
$^{A}\mu _{M}\in {{\mathfrak{M}^{H}}}$ and as $\lambda $ satisfies
relation b) of Definition \ref {ad invariant integ}, we get that
$\mu _{M}^{H}$ is also left $A$-linear: $\mu _{M}^{H}(a(m\otimes
h))=a_{0}m_{0}\lambda (a_1\vartriangleright m_{1}S(h))=a\mu
_{M}^{H}(m\otimes h).$\newline By (\ref{left integ}), we have:
$\lambda (xS(y_{1}))y_{2}=x_{1}\lambda (x_{2}S(y)),\forall x,y\in
H.$ Thus we get also the right $H$-collinearity of $\mu _{M}^H$:
$(\mu _{M}^{H}\otimes H)\rho ^{H}(m\otimes h)=m_{0}\otimes \lambda
(m_{1}S(h_{1}))h_{2}=m_{0}\otimes m_{1}\lambda (m_{2}S(h))=\rho
^{H}\mu _{M}^{H}(m\otimes h).$\newline It remains to prove that
$\mu _{M}^{H}$ is a retraction of $\eta _{M}^{H}$: $\mu
_{M}^{H}\eta _{M}^{H}(m)= m_{0}\lambda (m_{1}S(m_{2}))=m\lambda
(1_{H})=m$.\newline $ii)$ It is analogous to $i)$ by setting
${^{H}\mu _{M}}(h\otimes m)=\lambda (hS(m_{-1}))m_{0}.$\newline
$iii)$ We have to construct a functorial retract of $(\eta _{M})_{M\in {{%
_{A}^{H}\mathfrak{M}_{A}^{H}}}},$where $\eta _{M}=(_{M}^{H}\eta \otimes
H)\circ \eta _{M}^{H}.$ By the previous part, there are a functorial
retraction $(\mu _{M}^{H})_{M\in {{_{A}\mathfrak{M}_{A}^{H}}}}$ of $(\sigma
_{M}^{H})_{M\in {_{A}\mathfrak{M}_{A}^{H}}}$ and a functorial retract $({%
^{H}\mu _{M}})_{M\in {_{A}^{{H}}{{\mathfrak{M}_{A}}}}}$ of $({^{H}\sigma _{M}%
})_{M\in {_{A}^{H}\mathfrak{M}_A}}.$ Let us define the morphism $\mu
_{M}:H\otimes M\otimes H\rightarrow M$ by $\mu _{M}=\mu _{M}^{H}\circ ({%
^{H}\mu _{M}\otimes H}).$ Obviously it is a retraction of $\sigma
_{M}$ in ${_{A}\mathfrak{M}_{A}^{H}}$. It is easy to prove that
$\mu _{M}={^{H}\mu _{M}}\circ (H\otimes {\mu_{M}^{H}})$: hence one
gets that $\mu _{M}$ is a morphism in
${{_{A}^{H}\mathfrak{M}_{A}^{H}.}}$
\end{proof}

We can now consider the main result concerning $ad$-invariant
integrals. The equivalence $(1)\Leftrightarrow (3b)$ was proved in
a different way in \cite[Proposition 2.11]{AMS-Spliting}.

\begin{theorem}
\label{teo ad-inv}Let $H$ be a Hopf algebra over a field $%
K.$ The following assertions are equivalent:

$(1) $ There is an $ad$-invariant integral $\lambda \in H^{\ast }.$

$(2) $ The forgetful functor ${_{A}^{H}\mathfrak{M}_{A}^{H}}\rightarrow {{%
_{A}\mathfrak{M}_{A}}}$ is separable for any algebra $A$ in ${^{H}%
\mathfrak{M}^{H}}${.}

$(3) $ The forgetful functor ${_{H}^{H}\mathfrak{M}_{H}^{H}}\rightarrow {_{H}%
\mathfrak{M}_{H}}$ is separable.

$(3b) $ $H$ is coseparable in $({_{H}\mathfrak{M}_{H},\otimes ,K})$.

$(4) $ The forgetful functor ${_{H}^{H}\mathcal{YD}}\rightarrow {_{H}%
\mathfrak{M}}$ is separable.

$(4b)$ $K$ is $\mathcal{I}_{F}$-injective where $F$ is the forgetful functor
of $(4)$.
\end{theorem}

\begin{proof}
$(1)\Rightarrow (2)$ It follows by Lemma \ref{lem ad-inv}.\newline
$(2)\Rightarrow (3).$ Obvious.\newline $(3)\Leftrightarrow (3b).$
It is just Proposition \ref{pro cosep=>left inject} applied to
$\mathcal{M}=({_{H}\mathfrak{M}_{H},\otimes ,K})$.\newline
$(3)\Rightarrow (4).$ Take the notations of Examples \ref{ex of F}
and \ref{ex. of adjoint}. Since
$F_{2}:{_{H}^{H}\mathfrak{M}_{H}^{H}}\rightarrow
{_{H}\mathfrak{M}_{H}}$ is separable and $F_{2}=F_a\circ F_b$,
where ${F_b}:{_{H}^{H}\mathfrak{M}_{H}^{H}}\rightarrow {_{H}%
\mathfrak{M}_{H}^{H}}$ and {$F_a$}${:{_{H}\mathfrak{M}_{H}^{H}}\rightarrow {%
_{H}\mathfrak{M}_{H}},}$ then, by Theorem \ref{teo compos of
separable}, $F_b$ is separable. Consider the inverses $(F')^{-1}$
and $F^{-1}$ respectively of the functors $F'=(-)^{coH}
:{_{H}^{H}\mathfrak{M}_{H}^{H}}\rightarrow {_{H}^{H}\mathcal{YD}}$
and $F=(-)^{coH}:{_{H}\mathfrak{M}_{H}^{H}}\rightarrow
{_{H}\mathfrak{M}}$ (these are category equivalences; see
\cite[Theorem 5.7] {Schauenburg2}). One can easily check that
$F^{-1}\circ F_3=F_b\circ (F')^{-1}$. By Theorem \ref{teo compos
of separable}, $(F')^{-1}$ is separable so that $F^{-1}\circ F_3$,
and hence $F_3$, is a separable functor.\newline
$(4)\Rightarrow (4b).$ By Corollary \ref{coro Rafael} the separability of $%
F_{3}:{_{H}^{H}\mathcal{YD}}\rightarrow {_{H}\mathfrak{M}}$ (that
has $G_{3}$ as a right adjoint) implies that any object in
${_{H}^{H}\mathcal{YD}},$ in particular $K,$ is
$\mathcal{I}_{F_{3}}$-injective. \newline $(4b)\Rightarrow (1).$
Observe that $u_H$ can be regarded as a morphism in
${_{H}^{H}\mathcal{YD}}$, once $H$ is regarded as an object in
${_{H}^{H}\mathcal{YD}}$ via the action $\vartriangleright $
(defined in \ref {adj actions}) and the coaction given by the
comultiplication ${\Delta }$. In particular, $u_H$ belongs to $\mathcal{I}_{F_{3}}$%
: in fact the counit $\varepsilon _{H}$ of $H$ is a left linear
retraction
of $F_{3}(u_{H})$. Hence, since $K$ is $\mathcal{I}_{F_{3}}$%
-injective, there is $\lambda :H\rightarrow K$ in
${_{H}^{H}\mathcal{YD}}$ such that $\lambda\circ u_H=\Id_K,$ i.e.,
by Lemma \ref{ad-inv = retraction}, an $ad$-invariant integral.
\end{proof}

\begin{remark}
The following assertions are all equivalent to the existence of an $ad$%
-invariant integral $\lambda \in H^{\ast }.$

$(5)$ The forgetful functor ${\mathcal{YD}_{H}^{H}}\rightarrow
\mathfrak{M}%
{_{H}}$ is separable.

$(6) $ The forgetful functor ${_{H}\mathcal{YD}^{H}}\rightarrow {_{H}%
\mathfrak{M}}$ is separable and $S$ is bijective.

$(7) $ The forgetful functor ${^{H}\mathcal{YD}_{H}}\rightarrow {\mathfrak{M}%
_{H}}$ is separable and $S$ is bijective.

$(8)$ $K$ is $\mathcal{I}_{F}$-injective where $F$ is the forgetful functor
of $(5)$,$(6)$ or $(7)$.\newline
In fact, note that ${_{H}^{H}\mathfrak{M}_{H}^{H}\simeq \mathcal{YD}_{H}^{H}.%
}$ Since $\lambda $ is in particular a total integral, the antipode $S$ is
bijective and hence, by \cite[Corollary 6.4]{Schauenburg2}, we can also
assume ${{^{H}\mathcal{YD}_{H}}}\simeq {_{H}^{H}\mathfrak{M}_{H}^{H}\simeq {%
_{H}\mathcal{YD}}}^{H}.$ Now, by means of Lemma \ref{lem ad-invariant}, one
can proceed like in the proof of Theorem \ref{teo ad-inv}.
\end{remark}

\begin{theorem}
\label{teo fd ad-inv}Let $H$ be a finite dimensional Hopf algebra
over a field $K$ and let $D(H)$ be the Drinfeld Double. The
following assertions are equivalent:

$(i)$ There is an $ad$-invariant integral $\lambda \in H^{\ast }.$

$(ii)$ The forgetful functor ${_{D(H)}\mathfrak{M}}\rightarrow {_{H}%
\mathfrak{M}}$ is separable.

$(iii) $ $D(H)$ is separable in $({_{H}\mathfrak{M}_{H},\otimes
}_{H},H)$, i.e. $D(H)/H$ is separable.
\end{theorem}

\begin{proof}
$(i) \Leftrightarrow (ii)$. Since $H$ is finite dimensional, it
has bijective antipode. Hence we have
${_{H}^{H}\mathcal{YD}}\simeq$ ${_{H}\mathcal{YD}^{H}}\simeq
{_{D(H)}\mathfrak{M}}$. By Theorem \ref{teo ad-inv}, $(i) $ holds
if and only if the forgetful functor
${_{H}^{H}\mathcal{YD}}\rightarrow {_{H}\mathfrak{M}}$ is
separable if and only if ${_{D(H)}\mathfrak{M}}\rightarrow
{_{H}\mathfrak{M}}$ is separable.\newline $(ii) \Leftrightarrow
(iii) .$ It follows by Proposition \ref{lem sep=>left project for
rings} applied to the ring homomorphism $H\rightarrow D(H)
=H^{\ast cop}\bowtie H:h\mapsto \varepsilon _{H}\bowtie h.$
\end{proof}

\begin{proposition}
\label{coro Fsep}Let $H$ be a Hopf algebra with an $ad$-invariant integral $%
\lambda \in H^{\ast }$ and let $\mathcal{M}$ be either
${\mathfrak{M}^{H}}$ or ${^{H}\mathfrak{M}^{H}} $. For any algebra
$A$ in $\mathcal{M}$, we have:

i) $A$ is separable as an algebra in $\mathcal{M}$ if and only if
it is separable as an algebra in $\mathfrak{M}_{K}$.

ii) $A$ is formally smooth as an algebra in $\mathcal{M}$ if and
only if it is formally smooth as an algebra in
${\mathfrak{M}}_{K}$.
\end{proposition}

\begin{proof}
Since $H$ has an $ad$-invariant integral $\lambda$, by Lemma \ref{lem ad-inv}%
, the forgetful functor $F:{_{A}\mathcal{M}_{A}}\rightarrow {_{A}\mathfrak{M}%
_{A}}$ is separable. By Theorem \ref{teo separability of F} we
conclude.
\end{proof}

\section{Splitting algebra homomorphisms}\label{sec: Splitting algebra
homomorphisms}

We recall the following important result.

\begin{theorem}
(see \cite[Theorem 3.13]{AMS} ) \label{X formally smooth teo}Let
$(A,m,u)$ be an algebra in an abelian monoidal category
$(\mathcal{M},\otimes , \mathbf{1})$. Then the following
assertions are equivalent:

(a) $A$ is formally smooth as an algebra in ${\mathcal{M}}$.

(b) Let $\pi :E\rightarrow A$ be an algebra homomorphism in
$\mathcal{M}$ which is an epimorphism in $\mathcal{M}$ and let $I$
denote the kernel of $\pi $. Assume that there is $n\in
\mathbb{N}$ so that $I^{n}=0$ ($I$ is nilpotent). If for any
$r=1,\cdots ,n-1$
the canonical projection $p_{r}:E/I^{r+1}\rightarrow E/I^{r}$ splits in $%
\mathcal{M}$, then $\pi $ has a section which is an algebra
homomorphism in $\mathcal{M}$.
\end{theorem}

When ${\mathcal{M}}$ is ${\mathfrak{M}^{H}}$ the previous theorem
has the following application.

\begin{proposition}
\label{lem canonical map}Let $H$ be a Hopf algebra and let $A$ and
$E$ be algebras in ${\mathfrak{M}^{H}}$. Let $\pi :E\rightarrow A$
be an algebra homomorphism in ${\mathfrak{M}^{H}}$ which is
surjective. Assume that $A$ is formally smooth as an algebra in
${\mathfrak{M}^{H}}$ and that the kernel of $\pi $ is a nilpotent
ideal. Given an algebra homomorphism $f:H\rightarrow A$ in
${\mathfrak{M}^{H}}$, then $\pi $ has a section which is an
algebra homomorphism in ${\mathfrak{M}^{H}}$.
\end{proposition}

\begin{proof} It is similar to \cite[Theorem
2.13]{AMS-Spliting}, where $A=H=E/Rad(E)$ is semisimple and
$f=\Id_H$.\newline Let $I$ denote the kernel of $\pi $ and assume
there is an $n\in \mathbb{N}$ such that $I^{n}=0.$ First of all
let us observe that, since $\pi $ is a morphism in
${\mathfrak{M}^{H}}$, $I$ is
a subobject of $E$ in ${\mathfrak{M}^{H}}$ . Hence, for every $r>0$, ${I}%
^{r} $ is a subobject of $E$ and the canonical maps $E/{I}^{r+1}\rightarrow
E/{I}^{r}$ are morphisms in ${\mathfrak{M}^{H}}$. \newline
Now, the object ${I}^{r}/{I}^{r+1}$ has a natural module structure over $E/{I%
}\simeq A$, and hence, via $f$, a module structure over $H$. With respect to
this structure ${I}^{r}/{I}^{r+1}$ is an object in ${\mathfrak{M}}_{H}^{H}$.
Via the category equivalences ${\mathfrak{M}_{H}^{H}\simeq {_{K}\mathfrak{M}}%
}$, we get that ${I}^{r}/{I}^{r+1}$ is a cofree right comodule i.e. ${I}^{r}/%
{I}^{r+1}\simeq V\otimes H$ in ${\mathfrak{M}_{H}^{H},}$ for a suitable $%
V\in {{_{K}\mathfrak{M}}}$. In particular ${I}^{r}/{I}^{r+1}$ is
an injective comodule, so any canonical map
$E/{I}^{r+1}\rightarrow E/{I}^{r}$ has a section in
${\mathfrak{M}^{H}}$. \newline By Theorem \ref{X formally smooth
teo}, we conclude.
\end{proof}

The following remark is due to the referee.

\begin{remark} Let $H$ be a Hopf algebra with antipode $S_H$ and let $A$ be an algebra in ${\mathfrak{M}^{H}}$. Then the existence of an algebra
homomorphism $f:H\rightarrow A$ in ${\mathfrak{M}^{H}}$ is
equivalent to the fact that $A$ is isomorphic as an $H$-comodule
algebra to the smash product $A^{co(H)}\#H.$\\ In fact, using the
terminology of \cite[Definition 7.2.1, page 105]{Montgomery}, the
$H$-extension $A^{co(H)}\subseteq A$ comes out to be $H$-cleft
($f$ is a right $H$-comodule map which is convolution invertible
with inverse $f\circ S_H$). By \cite[Theorem 7.2.2, page
106]{Montgomery}, $A\simeq A^{co(H)}\#_{\sigma}H$, where
$\sigma:H\ot H\to A^{co(H)}$ is defined by $$\sigma(h\ot k)=\sum
f(h_1)f(k_1)fS_H(h_2k_2),$$ for every $h,k\in H$. Since $f$ is an
algebra homomorphism we get that $\sigma(h\ot
k)=\epsilon_H(h)\epsilon_H(k)1_A$ and hence
$A^{co(H)}\#_{\sigma}H=A^{co(H)}\#H$ is the usual smash product.\\
Conversely, for any algebra $R$, the map $H\to R\# H$ is an
algebra homomorphism in ${\mathfrak{M}^{H}}.$
\end{remark}

\begin{example}
Let $H$ be a Hopf algebra and assume that $H$ is formally smooth in ${%
\mathfrak{M}^{H}.}$ Then, by \cite[Corollary 3.30]{AMS}, the
tensor algebra $T:=T_{H}(\Omega ^{1}H)$ is formally smooth as an
algebra in the monoidal category ${\mathfrak{M}^{H}}$. Assume that
$\pi :E\rightarrow T$ is an epimorphism that is also a morphism of
algebras in ${\mathfrak{M}^{H}}$ such that $I:=\mathrm{Ker}\,\pi $
is a nilpotent coideal. By Proposition \ref{lem canonical map},
applied to the case when $f:H\rightarrow T$ is the canonical
injection, $\pi $ has a section which is an algebra homomorphism
in ${\mathfrak{M}^{H}}$. (In particular the projection
$E\rightarrow T\rightarrow H$ has also a section which is an
algebra homomorphism in ${\mathfrak{M}^{H}}$). Observe that, in
general, $T$ is not semisimple because its dimension needs not to
be finite.
\end{example}

\begin{theorem}
\label{te:section} Let $H$ be a Hopf algebra and $\mathcal{M}$ be
either ${\mathfrak{M}^{H}}$ or ${^{H}\mathfrak{M}^{H}}$. Let $E$
be an algebra in $\mathcal{M}$. Let $\pi :E\rightarrow H$ be an
algebra homomorphism in $\M$ which is surjective. Assume that $H$
is formally smooth as an algebra in $\mathcal{M}$ and that the
kernel $I$ of $\pi $ is a nilpotent ideal. Then $\pi $ has a
section which is an algebra homomorphism in $\mathcal{M}$ for

$a)$ $\mathcal{M}={\mathfrak{M}^{H}}$.

$b)$ $\mathcal{M}={^{H}\mathfrak{M}^{H}}$ if any canonical map $%
E/I^{r+1}\rightarrow E/I^{r}$ splits in $\mathcal{M}$.
\end{theorem}

\begin{proof}
Since $\pi $ is a morphism in $\mathcal{M}$, the kernel $I$ of
$\pi$ is a subobject of $E$ in $\mathcal{M}$. Hence, for every
$r>0$, $I^{r}$ is a subobject of $E$ and the canonical maps
$E/I^{r+1}\rightarrow E/I^{r}$ are morphisms in $\mathcal{M}$.
\newline
$a)$\ Apply Proposition \ref{lem canonical map} in the case when $E:=H$ and $%
f:=\mathrm{Id}_{H}.$ \newline $b)$ It follows easily by Theorem
\ref{X formally smooth teo}.
\end{proof}

Proposition \ref{lem canonical map} studies the existence in
${\mathfrak{M}^{H}}$ of algebra sections of morphisms of
algebras $\pi:E\rightarrow A$ where $A$ is a formally smooth algebra in ${%
{\mathfrak{M}^{H}}}$ endowed with a morphism of algebras $f:H\rightarrow A$ in ${%
{\mathfrak{M}^{H}}}$. The following results show that the existence of $ad$%
-invariant integrals provides such a section both in ${\mathfrak{M}%
^{H}}$ and ${^{H}\mathfrak{M}^{H}}$ (without $f$).

\begin{lemma}
\label{coro epi split}Let $H$ be a Hopf algebra with a total
integral $\lambda \in H^{\ast }.$ Let ${\mathcal{M}}$ be either
${\mathfrak{M}^{H}}$ or ${^{H}\mathfrak{M}^{H}}$. Then any
epimorphism in ${\mathcal{M}}$ has a section in ${\mathcal{M}}$.
\end{lemma}

\begin{proof}
Since $\lambda $ is a total integral in $H^{\ast }$, then, by
Theorem \ref{teo integrals}, $H$ is coseparable in
$\mathfrak{M}_{K}$. Therefore any right (resp. two-sided)
$H$-comodule is projective (see Corollary \ref{coro inj}). In
particular any any epimorphism in ${\mathcal{M}}$ has a section in
${\mathcal{M}}$.
\end{proof}

\begin{theorem}\label{teo ad and fs}
Let $H$ be a Hopf algebra with an $ad$-invariant integral $\lambda
\in H^{\ast }$. Let ${\mathcal{M}}$ be either ${\mathfrak{M}^{H}}$
or ${^{H}\mathfrak{M}^{H}}$. Let $A$ and $E$ be algebras in
${\mathcal{M}}$. Let $\pi :E\rightarrow A$ be an algebra
homomorphism in $\M$ which is surjective. Assume that $A$ is
formally smooth as an algebra in ${\mathfrak{M}}_{K}$ and that the
kernel of $\pi $ is a nilpotent ideal. Then $\pi $ has a section
which is an algebra homomorphism in $\mathcal{M}$.
\end{theorem}

\begin{proof}
By Proposition \ref{coro Fsep}, $A$ is formally smooth as an algebra in ${%
\mathcal{M}.}$ Let $n\geq 1$ such that $I^{n}=0$, where $I=%
\mathrm{Ker\,}\pi $. Since, in particular, $\lambda $ is a total
integral, by Lemma \ref{coro epi split}, any epimorphism in the
category ${\mathcal{M}} $ splits in ${\mathcal{M}}$. Thus, for
every $r=1,\cdots ,n-1$ the canonical morphism $\pi
_{r}:E/I^{r}\rightarrow E/I^{r+1}$ has a section in the category
${\mathcal{M}}$. We can now conclude by applying Theorem \ref{X
formally smooth teo} to the homomorphism of algebras $\pi
:E\rightarrow A.$
\end{proof}

\begin{theorem}
\label{coro 5.23} Let $H$ be a Hopf algebra with an $ad$-invariant
integral and such that $H$ is formally smooth as an algebra in
${_{K}\mathfrak{M}}$. Let ${\mathcal{M}}$ be either
${\mathfrak{M}^{H}}$ or ${^{H}\mathfrak{M}^{H}}$. Let $E$ be an
algebra in ${\mathcal{M}}$. Let $\pi :E\rightarrow H$ be a algebra
homomorphism in $\M$ which is surjective and with nilpotent
kernel. Then $\pi$ has a section which is an algebra homomorphism
in ${\mathcal{M}}$.
\end{theorem}

\begin{remark}
By Proposition \ref{coro Fsep}, if $H$ is a Hopf algebra with an $ad$%
-invariant integral and $H$ is formally smooth as an algebra in $({_{K}%
\mathfrak{M}},\otimes ,K),$ then it is formally smooth as an algebra in $( {%
\mathfrak{M}^{H},\otimes ,K}) .$ Then the case ${\mathcal{M}}={\mathfrak{M}%
^{H}}$ of the above corollary can be also deduced by Theorem
\ref{te:section}.
\end{remark}

\section{Formal Smoothness of a Hopf algebra as an
algebra}\label{sec: Formal Smoothness of a Hopf algebra as an
algebra}

In order to apply Theorem \ref{te:section}, it is useful to
characterize when the algebra $H$ is formally smooth either in
${\mathfrak{M}^{H}}$ or in $ {^{H}\mathfrak{M}^{H}}$.

\begin{claim}
Let $H$ be a Hopf algebra with antipode $S$ over a field $K.$ We denote by $%
H^{+}$ the \emph{augmentation ideal}, that is the kernel of the counit $%
\varepsilon :H\rightarrow K.$ \newline
Observe that $\varepsilon $ can be regarded as a morphism in ${_{H}^{H}%
\mathcal{YD}}$, once $H$ is regarded as an object in
${_{H}^{H}\mathcal{YD}}$ via the coaction ${^{H}\varrho }$
(defined in \ref{adj actions}) and the action given by the
multiplication ${m}$. In this way $H^+=\K(\varepsilon) $ inherits
the following structure of left-left Yetter-Drinfeld module:
\begin{equation*}
h\cdot x=hx,\text{\qquad }^{H}\rho (x)=x_{\left( 1\right)
}S(x_{\left( 3\right) })\otimes x_{\left( 2\right) }
\end{equation*}%
for all $h\in H$ and $x\in H^{+}.$\newline We call an\emph{\
}$\emph{fs}$\emph{-section} any map $\tau :H^{+}\rightarrow
H\otimes H^{+}$ such that:

$(i)$ $\tau (hx)=\sum_{i\in I}ha_{i}\otimes b_{i},$

$(ii)$ $\sum_{i\in I}a_{i}b_{i}=x,$\newline for all $h\in H$ and
$x\in H^{+},$ where $\tau (x)=\sum_{i\in I}a_{i}\otimes b_{i}$.

We say that an $fs$-section is \emph{complete} whenever

$(iii)$ $\sum_{i\in I}a_{i\left( 1\right) }b_{i\left( 1\right)
}S(a_{i\left( 3\right) }b_{i\left( 3\right) })\otimes a_{i\left(
2\right) }\otimes b_{i\left( 2\right) }=x_{\left( 1\right)
}S(x_{\left( 3\right) })\otimes \tau (x_{\left( 2\right) }).$

\end{claim}

\begin{lemma}
\label{lem fs-section} $\tau$ is a complete $fs$-section if and
only if $\tau $ is a section in ${_{H}^{H}\mathcal{YD}}$ of the
counit $\varepsilon _{H^{+}}:H\otimes H^{+}\to H^{+}$ of the
adjunction $({F^3},{G^3})$ introduced in \ref{ex. of adjoint}.
\end{lemma}

\begin{proof}
The notion of complete $fs$-section can be read as follows:
condition $(i) $ means that $\tau $ is left $H$-linear, $(iii) $
that $\tau $ is left $H$-colinear and $(ii) $ that $\tau $ is a
section of the counit $\varepsilon _{H^{+}}$ of
the adjunction $({F^3},{G^3})$, i.e. $\varepsilon _{H^{+}}\circ \tau =%
\mathrm{Id}_{H^{+}}$.
\end{proof}

\begin{proposition}
\label{pro: tau and fd} Let $H$ be a finite dimensional Hopf
algebra over a field $K$ and let $H^{+}$ be the augmentation
ideal. Let $\tau :H^{+}\rightarrow H\otimes H^{+}$ be a $K$-linear
map such that $$\tau (hx)=\sum_{i\in I}ha_{i}\otimes b_{i}$$ for
all $h\in H,x\in H^+$, where $\tau (x)=\sum_{i\in I}a_{i}\otimes
b_{i}$. Then $Im\left( \tau \right) \subseteq H^{+}\otimes H^{+}.$
\end{proposition}

\begin{proof}
Since $H$ is finite dimensional, there exists a non-zero right integral $%
t\in H.$

Let $x\in H^{+}.$ Since $Im\left( \tau \right) \subseteq H\otimes
H^{+},$ we can write $\tau (x)=\sum_{i\in I}a_{i}\otimes b_{i},$
$a_i\in H,b_i\in H^{+}.
$ We have%
\begin{eqnarray*}
\sum_{i\in I}ta_{i}\otimes b_{i} &=&\sum_{i\in I}t\varepsilon
\left( a_{i}\right) \otimes b_{i}=t\otimes \sum_{i\in
I}\varepsilon \left(
a_{i}\right) b_{i} \\
\tau (tx) &=&\tau (t\varepsilon \left( x\right) )=0.
\end{eqnarray*}
Therefore, since $\tau (tx)=\sum_{i\in I}ta_{i}\otimes b_{i},$ we get $\sum_{i\in I}\varepsilon \left( a_{i}\right) b_{i}=0.$ Hence $%
Im\left( \tau \right) \subseteq \ker \left( \varepsilon \otimes
H^{+}\right) =H^{+}\otimes H^{+}.$
\end{proof}

\begin{proposition}
\label{formal in bicomodules}Let $H$ be a Hopf algebra with
antipode $S$ over a field $K$ and let $H^{+}$ be the augmentation
ideal. The following assertions are equivalent:

$(a) $ $H$ is formally smooth as an algebra in
${^{H}\mathfrak{M}^{H}.}$

$(b) $ $H^{+}$ is $\mathcal{E}_{G}$-projective where $G$ is the
forgetful functor ${_{H}^{H}\mathcal{YD}}\rightarrow
{^{H}\mathfrak{M}}$.

$(c)$ There exists a complete $fs$-section $\tau :H^{+}\rightarrow
H\otimes H^{+}.$

Moreover, if $H$ finite dimensional, the following assertion is
also equivalent to the others:

$(d)$ The multiplication $H^{+}\otimes H^{+}\rightarrow H^{+}$ has a left $H$%
-linear section $\tau :H^{+}\rightarrow H^{+}\otimes H^{+}$, where $%
H^{+}\otimes H^{+}$ is a left $H$-module via $^{H}\mu
_{H^{+}}\otimes H^{+}$ and such that $\sum_{i\in I}a_{i\left(
1\right) }b_{i\left( 1\right) }S(a_{i\left( 3\right) }b_{i\left(
3\right) })\otimes a_{i\left( 2\right) }\otimes b_{i\left(
2\right) }=x_{\left( 1\right) }S(x_{\left( 3\right) })\otimes \tau
(x_{\left( 2\right) })$ for all $x\in H^{+},$ where $\tau
(x)=\sum_{i\in I}a_{i}\otimes b_{i}$.
\end{proposition}

\begin{proof}
$(b)\Leftrightarrow (c)$ Consider the functor $G^{3}:{_{H}^{H}\mathcal{YD}}%
\rightarrow {^{H}\mathfrak{M}}$ and it's left adjoint $F^{3}:{^{H}\mathfrak{M%
}}\rightarrow {_{H}^{H}\mathcal{YD}}$ (see \ref{ex. of adjoint}).
We know
(see Theorem \ref{teo 4.2.33}) that $H^{+}$ is $\mathcal{E}_{G^{3}}$%
-projective if and only if the counit of the adjunction
$\varepsilon _{H^{+}}:F^{3}G^{3}(H^{+})\rightarrow H^{+}$ has a
section $\tau :H^{+}\rightarrow H\otimes H^{+}$ in
${_{H}^{H}\mathcal{YD}}$: thus, by Lemma \ref{lem fs-section},
$\tau $ is a complete $fs$-section.\newline $(a)\Leftrightarrow
(b)$ In view of Examples \ref{ex of co F}, consider the following
diagrams:
\begin{equation*}
\xymatrix@R=15pt@C=50pt{
  {{^{H}\mathfrak{M}_{H}^{H}}} \ar[d]_{\mathbb{T}=F^b} \ar[r]^{F'=(-)^{coH}}_{\sim} & {^{H}\mathfrak{M}} \ar[d]^{\mathbb{T}'=F^3} \\
  {_{H}^{H}\mathfrak{M}_{H}^{H}} \ar[r]_{F=(-)^{coH}}^{\sim} & {_{H}^{H}\mathcal{YD}}   }
\text{ \qquad } \xymatrix@R=15pt@C=50pt{
  {{^{H}\mathfrak{M}_{H}^{H}}} \ar[r]^{F'=(-)^{coH}}_{\sim} & {^{H}\mathfrak{M}}  \\
  {_{H}^{H}\mathfrak{M}_{H}^{H}} \ar[u]^{\mathbb{H}=G^b} \ar[r]_{F=(-)^{coH}}^{\sim} & {_{H}^{H}\mathcal{YD}} \ar[u]_{\mathbb{H}'=G^3}  }
\end{equation*}
and the forgetful functor $G^{a}:{^{H}\mathfrak{M}_{H}^{H}}\rightarrow {^{H}%
\mathfrak{M}^{H}}$. The second diagram is commutative. Since
$G^{\prime }\circ G^{3}=G^{b}\circ G$, by the uniqueness of the
adjoint, it is
straightforward to prove that the functors $F^{3}\circ F^{\prime }$ and $%
F\circ F^{b}$ are naturally equivalent. By definition, $H$ is
formally
smooth in ${^{H}\mathfrak{M}^{H}}$, if and only if $\Omega ^{1}H$ is ${%
\mathcal{E}_{G^{a}\circ G^{b}}}$-projective (In fact ${_{H}\mathcal{E}_{H}=%
\mathcal{E}_{G^{2}}}$ and $G^{2}={G^{a}\circ G^{b}}$). By
Proposition \ref{2
sep functors}, $G^{a}$ is separable, so that, by Lemma \ref{lemma non so}, ${%
\mathcal{E}_{G^{b}}=\mathcal{E}_{G^{a}\circ G^{b}}}$. Moreover, the functor $%
F$ is separable as an equivalence of categories so that, by
Theorem \ref{teo F and P-project}, $\Omega ^{1}H$ is
${\mathcal{E}_{G^{b}}}$-projective if
and only if $H^{+}\simeq F(\Omega ^{1}H)$ (see \cite[Example 5.8]%
{Schauenburg2} for this isomorphism) is $\mathcal{E}_{G^{3}}$-projective.%
\newline
$(c)\Rightarrow (d)$ By Proposition \ref{pro: tau and fd},
$\text{Im} \left(
\tau \right) \subseteq H^{+}\otimes H^{+}$ so that $\tau ,$ corestricted to $%
H^{+}\otimes H^{+},$ is the required left $H$-linear section of
the multiplication $H^{+}\otimes H^{+}\rightarrow H^{+}$.\newline
$(d)\Rightarrow (c)$ Trivial.
\end{proof}

\begin{proposition}
\label{formal in comodules}Let $H$ be a Hopf algebra over a field
$K$ and let $H^{+}$ be the augmentation ideal. The following
assertions are equivalent:

$(a)$ $H$ is formally smooth as an algebra in
${\mathfrak{M}}_{K}{.}$

$(b)$ $H$ is formally smooth as an algebra in
${\mathfrak{M}^{H}}$.

$(c)$ $H^{+}$ is projective in $_H\mathfrak{M}$.

$(d)$ There exists an $fs$-section $\tau :H^{+}\rightarrow
H\otimes H^{+}$.

$(e)$ $H$ is a hereditary $K$-algebra.

Moreover, if $H$ finite dimensional, the following assertion is
also equivalent to the others:

$\left( f\right) $ The multiplication $H^{+}\otimes
H^{+}\rightarrow H^{+}$
has a left $H$-linear section, where $H^{+}\otimes H^{+}$ is a left $H$%
-module via $^{H}\mu _{H^{+}}\otimes H^{+}.$

\end{proposition}

\begin{proof}
The equivalences between $\left( b\right) ,\left( c\right) ,\left(
d\right) $ and $\left( f\right) $ follow similarly to Proposition
\ref{formal in bicomodules}, but working with the following
diagrams:
\begin{equation*}
\xymatrix@R=15pt@C=50pt{
  {{\mathfrak{M}_{H}^{H}}} \ar[d]_{\mathbb{T}} \ar[r]^{F'=(-)^{coH}}_{\sim} & {_K\mathfrak{M}} \ar[d]^{\mathbb{T}'} \\
  {_{H}\mathfrak{M}_{H}^{H}} \ar[r]_{F=(-)^{coH}}^{\sim} & {_{H}\mathfrak{M}}   }
\text{ \qquad } \xymatrix@R=15pt@C=50pt{
  {{\mathfrak{M}_{H}^{H}}} \ar[r]^{F'=(-)^{coH}}_{\sim} & {_K\mathfrak{M}}  \\
  {_{H}\mathfrak{M}_{H}^{H}} \ar[u]^{\mathbb{H}} \ar[r]_{F=(-)^{coH}}^{\sim} & {_{H}\mathfrak{M}} \ar[u]_{\mathbb{H}'}  }
\end{equation*}
One can check that $(b)\Leftrightarrow (c)\Leftrightarrow H^{+}$ is $%
\mathcal{E}_{\mathbb{H^{\prime }}}$-projective (where
$\mathbb{H^{\prime }}$
is the forgetful functor ${_{H}\mathfrak{M}}\rightarrow {_{K}\mathfrak{M}}$%
). Now, since $K$ is a field, we have that $\mathcal{E}_{\mathbb{H^{\prime }}%
}=\{g\in {_{H}\mathfrak{M}}\mid g\text{ is a surjection}\}$, so
that $H^{+}$ is $\mathcal{E}_{\mathbb{H^{\prime }}}$-projective if
and only if $H^{+}$ is projective in $_H\mathfrak{M}$.\newline
$\left( b\right) \Rightarrow \left( a\right) $ Apply Theorem
\ref{teo
separability of F} In the case when $A=H$ and $\mathcal{M}={\mathfrak{M}^{H}.%
}$\newline $\left( a\right) \Rightarrow \left( e\right) $ See \cite[Proposition 5.1]{CQ}.\newline $\left( e\right) \Rightarrow \left( c\right) $ Every left $H$%
-submodule of a projective left $H$-module is projective. In
particular any left ideal of $H$ is projective in
$_H\mathfrak{M}$.
\end{proof}

\begin{remark}
Let $H$ be a Hopf algebra with antipode $S$ over a field $K$. Then
$H$ is formally smooth as an algebra in
${^H\mathfrak{M}^{H}}\Rightarrow$ $H$ is formally smooth as an
algebra in ${\mathfrak{M}^{H}}$.
\end{remark}

\begin{corollary}
\label{coro: ad-inv + fs} Let $H$ be a Hopf algebra over a field
$K$. Assume that $H$ has an $ad$-invariant integral. Let $H^{+}$
be the augmentation ideal. The following assertions are
equivalent:

$(i)$ ${H}$ is formally smooth as an algebra in
${\mathfrak{M}}_{K}$.

$(ii)$ ${H}$ is formally smooth as an algebra in
${\mathfrak{M}^{H}}$.

$(iii)$ ${H}$ is formally smooth as an algebra in
${^{H}\mathfrak{M}^{H}}$.

$(iv)$ $H^{+}$ is projective in $_H\mathfrak{M}$.

$(v)$ There exists an $fs$-section $\tau :H^{+}\rightarrow
H\otimes H^{+}$.

$(vi)$ $H$ is an hereditary $K$-algebra.

Moreover, if $H$ finite dimensional, the following assertion is
also equivalent to the others:

$\left( vii\right) $ The multiplication $H^{+}\otimes
H^{+}\rightarrow H^{+}$
has a left $H$-linear section, where $H^{+}\otimes H^{+}$ is a left $H$%
-module via $^{H}\mu _{H^{+}}\otimes H^{+}.$
\end{corollary}

\begin{proof}
Let ${\mathcal{M}}$ be either ${\mathfrak{M}^{H}}$ or
${^{H}\mathfrak{M}^{H}}$ and observe that $H$ is an algebra in
$\mathcal{M}$. Then by
Proposition \ref{coro Fsep}, $H$ is formally smooth as an algebra in $%
\mathcal{M}$ if and only if it is formally smooth as an algebra in
${\mathfrak{M}}_{K}$ that is $\left( i\right) ,\left( ii\right) $
and $\left( iii\right) $ are equivalent. By Proposition
\ref{formal in comodules}, we conclude.
\end{proof}

By applying Corollary \ref{coro: ad-inv + fs}, we obtain the
following result to be compared with \cite[Theorem 2]{LB}.

\begin{theorem}
\label{teo: KG is fs}Let $G$ be an arbitrary group an let $KG$ be
the group algebra associated. Then the following assertions are
equivalent:

$(i)$ ${K}G$ is formally smooth as an algebra in
${\mathfrak{M}}_{K}$.

$(ii)$ ${K}G$ is formally smooth as an algebra in
${\mathfrak{M}^{KG}}$.

$(iii)$ ${K}G$ is formally smooth as an algebra in ${^{KG}\mathfrak{M}^{KG}}$%
.

$(iv)$ The augmentation ideal $KG^{+}$ is a projective in
$_{KG}\mathfrak{M}$.

$\left( v\right) $ There exists an $fs$-section $\tau
:KG^{+}\rightarrow KG\otimes KG^{+} $.

$(vi)$ $KG$ is an hereditary $K$-algebra.

$\left( vii\right) $ $G$ is the fundamental group of a connected
graph of finite groups whose orders are invertible in $K.$ (see
\cite[Definition 4.2, page 10]{Dicks}).

Moreover, if $G$ is finite, the following assertion is also
equivalent to the others:

$\left( viii\right) $ The multiplication $KG^{+}\otimes
KG^{+}\rightarrow KG^{+}$ has a left $KG$-linear section, where
$KG^{+}\otimes KG^{+}$ is a left $KG$-module via $^{KG}\mu
_{KG^{+}}\otimes KG^{+}.$
\end{theorem}

\begin{proof}
By the left analogue of \cite[Theorem 2.12, page 118]{Dicks},
$\left( iv\right) $ and $\left( vii\right) $ are equivalent.

By Example \ref{ex ad inv of KG}, the Hopf algebra $KG$ admits an $ad$%
-invariant integral. The conclusion follows by Corollary
\ref{coro: ad-inv + fs}.
\end{proof}

By means of Proposition \ref{formal in comodules}, it is now
possible to rewrite Theorem \ref{te:section} in the following form
which improves Theorem \ref{coro 5.23} in the case
$\M={\mathfrak{M}^{H}}$.

\begin{theorem}\label{teo: precoweak}Let $H$ be a Hopf algebra and let $E$ be an algebra in ${\mathfrak{M}^{H}}$. Let $\pi :E\rightarrow H$ be an algebra
homomorphism in ${\mathfrak{M}^{H}}$ which is surjective. Assume
that $H$ is formally smooth as an algebra in ${\mathfrak{M}}_{K}$
and that the kernel $I$ of $\pi $ is a nilpotent ideal. Then $\pi
$ has a section which is an algebra homomorphism in
${\mathfrak{M}^{H}}$.
\end{theorem}

As a consequence of Theorem \ref{teo: precoweak}, we get the
following result.

\begin{theorem}\label{teo: coweak}
Let $H$ be a Hopf algebra and let $E$ be a bialgebra. Let $\pi
:E\rightarrow H$ be a bialgebra homomorphism which is surjective.
Assume that $H$ is formally smooth as an algebra in
${\mathfrak{M}}_{K}$ and that the kernel $I$ of $\pi $ is a
nilpotent ideal. Then $\pi $ has a section which is an algebra
homomorphism in ${\mathfrak{M}^{H}}$.
\end{theorem}

\begin{remark}
Akira Masuoka pointed out that, in the situation of Theorem
\ref{teo: coweak}, since $H$ is a Hopf algebra so is $E$ (see e.g.
\cite[Lemma 3.52]{AMS-Spliting}).
\end{remark}

\section{Examples}\label{sec: Examples}

\begin{proposition}
\label{pro: groupalg on Integers}Let $K$ be any field. The group algebra $K%
\Z $ over the set of integers admits a complete $fs$-section.
\end{proposition}

\begin{proof}
Let $\left\langle g\right\rangle $ be the multiplicative group
associated to $\Z $ $\left( \left\langle g\right\rangle \simeq \Z
\right) $. Let $H=K\left\langle g\right\rangle .$
Then\begin{equation*} \mathcal{B}\left( H\right) =\left(
g^{n}-g^{n+1}\right) _{n\in \Z }
\end{equation*}%
is a basis for $H^{+}.$
Now define $\tau :H^{+}\rightarrow H\otimes H^{+}$ on generators
by setting
\begin{equation*}
\tau \left( g^{n}-g^{n+1}\right) =g^{n}\otimes \left( 1-g\right) ,
\end{equation*}%
for every $n\in\Z .$ Clearly $g^{n}\cdot \left( 1-g\right)
=g^{n}-g^{n+1}$. Moreover%
\begin{equation*}
\tau \left[ g^{a}\left( g^{n}-g^{n+1}\right) \right] =\tau \left(
g^{a+n}-g^{a+n+1}\right) =g^{a+n}\otimes \left( 1-g\right)
=g^{a}\cdot g^{n}\otimes \left( 1-g\right) =g^{a}\cdot \tau \left[
\left( g^{n}-g^{n+1}\right) \right] .
\end{equation*}%
Since $H$ is cocommutative, this is enough to conclude that $\tau $ is a complete $fs$%
-section of $H.$
\end{proof}

\begin{remarks}\label{rem: GroupAlg}
1) By Proposition \ref{pro: groupalg on Integers} and Theorem
\ref{teo: KG is fs}, $K\Z$ is formally smooth as an algebra in
$\mathfrak{M}_{K}.$ Nevertheless, being not finite dimensional,
$K\Z$ is not separable as an algebra in $\mathfrak{M}_{K}$. More
generally, the group algebra $KG$ is formally smooth but not
separable if and only if $G$ is a free and non-trivial group (see
Remark \ref{rem: GroupCohom}).

2) The complete $fs$-section $\tau $ defined in the proof of
Proposition \ref{pro: groupalg on Integers} is such that
$\Im\left( \tau \right) \nsubseteq K\Z^{+}\otimes K\Z^{+}.$ This
is a counterexample for the last assertion of Proposition
\ref{formal in bicomodules}.
\end{remarks}

\begin{proposition}
\label{pro: cyclic group}Let $C_{n}$ be the cyclic group of order
$n$ and let $KC_{n}$ be the group algebra associated. Then the
following assertions are equivalent:

$(i)$ ${K}C_{n}$ is formally smooth as an algebra in
${\mathfrak{M}}_{K}$.

$(ii)$ ${K}C_{n}$ is separable as an algebra in
${\mathfrak{M}}_{K}$.

$(iii)$ ${n\cdot 1}_{K}\neq 0$.
\end{proposition}

\begin{proof}
$\left( ii\right) \Leftrightarrow \left( iii\right) $ is the well
known Maschke's Theorem.\newline
$\left( ii\right) \Rightarrow \left( i\right) $ follows by Corollary \ref%
{coro: sep=>fs}. \newline $\left( i\right) \Rightarrow \left(
iii\right) $ By Theorem \ref{teo: KG is fs}, the multiplication
$KC_{n}^{+}\otimes KC_{n}^{+}\rightarrow KC_{n}^{+}$ has a
section. In particular the multiplication is surjective, so that $%
KC_{n}^{+}=\left( KC_{n}^{+}\right) ^{2}.$\newline Let $g\in
C_{n}$ be a generator of $C_{n},$ that is $o\left( g\right) =n.$
Then $KC_{n}^{+}=\sum_{i=0}^{n-1}K\left( 1-g^{i}\right) .$ From
$1-g\in KC_{n}^{+}=\left( KC_{n}^{+}\right) ^{2},$ we deduce there
exists $\alpha _{i,j}\in K $ such that
\begin{equation}
1-g=\sum_{0\leq i,j\leq n-1}\alpha _{i,j}\left( 1-g^{i}\right)
\left( 1-g^{j}\right) =\sum_{0\leq i,j\leq n-1}\alpha _{i,j}\left(
1-g^{i}-g^{j}+g^{i+j}\right) .  \label{formula 1-g}
\end{equation}%
Define the $K$-linear map $\varphi :{K}C_{n}\rightarrow K $ by setting $\varphi \left( g^{i}\right) =\left( 1-i\right) 1_{K}$ for every $%
0\leq i\leq n-1.$ Now suppose that ${n\cdot 1}_{K}=0.$ In this case, since $%
n=o\left( g\right) ,$ it is easy to check that $\varphi \left(
g^{i}\right) =\left( 1-i\right) 1_{K}$ for every $i\in \N $ and hence, by (\ref{formula 1-g}), we have%
\begin{equation*}
1=\varphi \left( 1-g\right) =\sum_{0\leq i,j\leq n-1}\alpha
_{i,j}\varphi \left( 1-g^{i}-g^{j}+g^{i+j}\right) =0,
\end{equation*}%
a contradiction.
\end{proof}

\begin{claim}
Implication $\left( i\right) \Rightarrow \left( iii\right) $ of Proposition %
\ref{pro: cyclic group}, can be proved in a different way. In fact
$\left(
i\right) $ implies that the Hochschild cohomology $H^{2}\left( {K}%
C_{n},M\right) $ vanishes for every ${K}C_{n}$-bimodule $M.$ By \cite[%
Theorem 5.5, page 292]{McL} (where the result is proved for
$\mathbb{Z}$ instead of $K$ although the same arguments go through
for any commutative ring),
for every group $G,$ one has a natural isomorphism%
\begin{equation*}
H^{t}\left( {K}G,M\right) \simeq H^{t}\left( G,{_{\chi }M}\right)
\end{equation*}%
where $_{\chi }M$ is $M$ endowed with the left $G$-module
structure given by $g\cdot _{\chi }m=gmg^{-1}$ and $H^{t}\left(
G,{_{\chi }M}\right) $ denotes the group cohomology. Apply this
isomorphism to the case $G=C_{n},t=2$ and let $g$ denote a
generator of $C_{n}$.\newline By \cite[Theorem 7.1, page
122]{McL}, for every left $C_{n}$-module $L,$ one has
\begin{equation*}
H^{2}\left( C_{n},L\right) =\frac{\left\{ l\in L\mid g\cdot l=l\right\} }{%
t\cdot L},
\end{equation*}%
where $t=1+g+g^{2}+\cdots +g^{n-1}$.\newline Now, assume that
$\left( i\right) $ holds. Then
\begin{equation*}
\frac{\left\{ x\in {K}C_{n}\mid g\cdot _{\chi }x=x\right\} }{t\cdot _{\chi }{%
K}C_{n}}=H^{2}\left( C_{n},{_{\chi }{K}C_{n}}\right) \simeq H^{2}\left( {K}%
C_{n},{K}C_{n}\right) =0.
\end{equation*}%
Since $C_{n}$ is commutative, then $g^{i}\cdot _{\chi
}x=g^{i}xg^{-i}=g^{i}g^{-i}x=x$ for every $x\in {K}C_{n}$ so that
\begin{equation*}
{K}C_{n}=\left\{ x\in {K}C_{n}\mid g\cdot _{\chi }x=x\right\}
=t\cdot _{\chi
}{K}C_{n}=\sum_{0\leq i\leq n-1}g^{i}\cdot {K}C_{n}\cdot g^{-i}=n\cdot {K}%
C_{n}.
\end{equation*}%
Therefore $1\in {K}C_{n}=n\cdot {K}C_{n}$ and hence ${n\cdot
1}_{K}\neq 0.$
\end{claim}

\begin{remark}\label{rem: GroupCohom}
By Proposition \ref{pro: cyclic group}, ${K}C_{n}$ is formally
smooth as an algebra in $\mathfrak{M}_{K}$ if and only if it is
separable as an algebra in $\mathfrak{M}_{K}$. \newline
The groups $G$ such that $KG$ is \textbf{formally smooth as an algebra} in $\mathfrak{%
M}_{K}$ \textbf{but not separable} as an algebra in
$\mathfrak{M}_{K}$ are precisely those having
\textbf{cohomological dimension $1.$} This follows in view of the
isomorphism $H^{t}\left( {K}G,M\right) \simeq H^{t}\left( G,{_{\chi }M}%
\right) $ that holds for every $t\in
\mathbb{N}
,$ and for any $KG$-bimodule $M.$ Note also that every left
$G$-module $N$ can be seen as ${_{\chi }}\left(
_{KG}{N}_{KG}\right) $ where $_{KG}{N}_{KG}$
is $N$ itself regarded as a bimodule via $g\cdot n\cdot h:=gn,$ for every $%
g,h\in G,n\in N.$\newline Furthermore (see \cite[%
Example 2, page 185]{Br}) every free group over a non-empty
(possibly infinite) set has cohomological dimension $1$.
Conversely every group of cohomological dimension $1$ is free.
\end{remark}

\section{Ad-coinvariant integrals through separable
functors}\label{sec: Ad-coinvariant integrals through separable
functors}

We want now to treat the dual of all the results of the previous
sections. We just state the main results that can be proved
analogously.\medskip\newline First of all we characterize the
existence of a so-called $ad$-coinvariant integral.\newline A
remarkable fact is that any semisimple and cosemisimple Hopf
algebra $H$ over a field $K$ admits such an integral (see
\cite[Theorem 2.27]{AMS-Spliting}).

\begin{definition}
Let $H$ be a Hopf algebra with antipode $S$ over any field $K$ and let $t\in
H$. \newline
$t$ will be called an $\emph{ad}$\emph{-coinvariant} \emph{integral}
whenever:

$a$) $ht=\varepsilon _{H}(h) t$ for all $h\in H$ (i.e. $t$ is a left
integral in $H$);

$b$) $t_{1}S(t_{3})\otimes t_{2}=1_{H}\otimes t,$ (i.e. $t$ is left
coinvariant with respect to $^{H}\varrho $)$;$

$c$) $\varepsilon _{H}(t) =1_{K}$. \newline
Therefore we have:
\end{definition}

\begin{lemma}
An element $t\in H$ is an $ad$-coinvariant integral if and only if the map $%
\tau :K\rightarrow H:k\mapsto kt$ is a section of the counit $\varepsilon
_{H}:H\to K$ of $H$ in ${_{H}^{H}\mathcal{YD,}}$ where $H$ is regarded as an
object in the category via the left adjoint coaction $^{H}\varrho $ and the
multiplication $m_{H}$.
\end{lemma}

\begin{example}
\label{ex ad coinv of K alla G} 1) Let $G$ be a finite group an
let ${K}^G$ be the algebra of functions from $G$ to $K$. Then
$K^G$ becomes a Hopf algebra which is dual to the group algebra
$KG$. From Example \ref{ex ad inv of KG}, we infer that $K^G$ has
an $ad$-coinvariant integral, namely the map $G\rightarrow
K:g\mapsto \delta _{e,g}$ (the Kronecker symbol), where $e$
denotes the neutral element of $G$.\\2) Every cocommutative
semisimple Hopf algebra has an $ad$-coinvariant integral.
\end{example}

\begin{remark}
It is known that, for any Hopf algebra $H$ with a total integral
$t\in H$, the $K$-linear spaces of left and right integrals in $H$
are both one dimensional and so both generated by $t$. Hence there
can be only one $\emph{ad}$-coinvariant integral, namely the
unique total integral.
\end{remark}

The following lemma shows that in the definition of $ad$-coinvariant
integral we can choose $\varrho ^{H},\overline{\varrho }^{H}$ or $^{H}%
\overline{\varrho }$ instead of $^{H}\varrho .$ Since $t$ is in particular a
total integral, it is both a left integral and a right integral. Thus it is the same to have a retraction of $%
\varepsilon _{H}$ in ${_{H}^{H}\mathcal{YD}},{\mathcal{YD}_{H}^{H},_{H}%
\mathcal{YD}^{H}}$ or ${^{H}\mathcal{YD}_{H}.}$

\begin{lemma}
\label{lem ad-coinvariant}Let $H$ be a Hopf algebra with antipode $S$ over
any field $K$ and let $t\in H$ be a total integral. Then the following are
equivalent:

$(1)$ $t$ is left coinvariant with respect to $^{H}\varrho $.

$(2)$ $t$ is right coinvariant with respect to $\varrho ^{H}$.

$(3)$ $t$ is right coinvariant with respect to $\overline{\varrho }^{H}$.

$(4)$ $t$ is left coinvariant with respect to $^{H}\overline{\varrho }$.
\end{lemma}

\begin{proof}
Analogous to \ref{lem ad-invariant}$.$
\end{proof}

\begin{lemma}
\label{lem ad-coinv}Let $H$ be a Hopf algebra with antipode $S$ over a field
$K.$ Assume there exists an $ad$-coinvariant integral $t\in H.$ Then we have
that:

$i)$ The forgetful functor ${{^{C}\mathfrak{M}_{H}^{C}}\rightarrow {^{C}%
\mathfrak{M}^{C}}}$ is separable for any coalgebra $C$ in ${{\mathfrak{M}_{H}%
}}$.

$ii)$ The forgetful functor ${_{H}^{C}\mathfrak{M}^{C}}\rightarrow {{^{C}%
\mathfrak{M}^{C}}}$ is separable for any coalgebra $C$ in ${_{H}\mathfrak{M}}
$.

$iii)$ The forgetful functor ${_{H}^{C}\mathfrak{M}_{H}^{C}}\rightarrow {^{C}%
\mathfrak{M}^{C}}$ is separable for any coalgebra $C$ in $_{H}\mathfrak{M}%
_{H}$.
\end{lemma}

\begin{proof}
We proceed as in the proof of Lemma \ref{lem ad-inv}.$\newline
i)$ By Examples \ref{ex of co F}, the forgetful functor $G^r:{^{C}%
\mathfrak{M}_{H}^{C}}\rightarrow {^{C}\mathfrak{M}^{C}}$ has a right adjoint
$F^r:{^{C}\mathfrak{M}^{C}}\rightarrow {^{C}\mathfrak{M}_{H}^{C}}%
,F^r(M)=M\otimes H$. Thus by Theorem \ref{teo Rafael}, $G^r$ is
separable if and only if the counit $\varepsilon
^{H}:F^rG^r\rightarrow \mathrm{Id}_{{{^{C}\mathfrak{M}_{H}^{C}}}}$
of the adjunction splits, i.e.
there exists a natural transformation $\sigma ^{H}:\mathrm{Id}_{{{^{C}%
\mathfrak{M}_{H}^{C}}}}\rightarrow F^rG^r$ such that $\varepsilon
_{M}^{H}\circ \sigma _{M}^{H}=\mathrm{Id}_{M}$ for any $M$ in ${{^{C}%
\mathfrak{M}_{H}^{C}}.}$ Using (\ref{left integ 1}), one can
easily check that the following map works: $\sigma
_{M}^{H}:M\rightarrow M\otimes H$, $\sigma
_{M}^{H}(m)=mt_{1}\otimes S(t_{2}).\newline ii)$ Analogous to $i)$
by setting $^{H}\sigma _{M}(m)=t_{1}\otimes S(t_{2}) m.\newline
iii)$ Define $\sigma _{M}:=(^{H}\sigma _{M}\otimes H) \circ \sigma
_{M}^{H}:M\rightarrow H\otimes M\otimes H.$
\end{proof}

We can now consider the main result concerning $ad$-coinvariant
integrals. The equivalence $(1)\Leftrightarrow (3b)$  was proved
in a different way in \cite[Proposition 2.11]{AMS-Spliting}.

\begin{theorem}
\label{teo ad-coinv}Let $H$ be a Hopf algebra over a field $K.$
The following assertions are equivalent:

$(1)$ There is an $ad$-coinvariant integral $t\in H.$

$(2)$ The forgetful functor ${_{H}^{C}\mathfrak{M}_{H}^{C}}\rightarrow {^{C}%
\mathfrak{M}^C}$ is separable for any coalgebra $C$ in $_{H}\mathfrak{M}_{H}$%
.

$(3)$ The forgetful functor ${_{H}^{H}\mathfrak{M}_{H}^{H}\rightarrow {^{H}%
\mathfrak{M}^{H}}}$ is separable.

$(3b)$ $H$ is separable in $({^{H}\mathfrak{M}^{H},\otimes ,K}){.}$

$(4)$ The forgetful functor ${_{H}^{H}\mathcal{YD}}\rightarrow {{^{H}%
\mathfrak{M}}}$ is separable.

$(4b)$ $K$ is $\mathcal{E}_{G}$-projective where $G$ is the forgetful
functor of $(4)$.
\end{theorem}

\begin{proof}
Analogous to that of Theorem \ref{teo ad-inv}.
\end{proof}

\begin{remark}
The following assertions are all equivalent to the existence of an $ad$%
-coinvariant integral $t\in H$:

$(5)$ The forgetful functor ${\mathcal{YD}_{H}^{H}}\rightarrow {{\mathfrak{M}%
^{H}}}$ is separable.

$(6)$ The forgetful functor ${_{H}\mathcal{YD}^{H}}\rightarrow {{\mathfrak{M}%
^{H}}}$ is separable and $S$ is bijective.

$(7)$ The forgetful functor ${^{H}\mathcal{YD}_{H}}\rightarrow {{^{H}%
\mathfrak{M}}}$ is separable and $S$ is bijective.

$(8)$ $K$ is $\mathcal{E}_{G}$-projective where $G$ is the forgetful functor
of $(5)$,$(6)$ or $(7)$.\newline
In fact, note that ${_{H}^{H}\mathfrak{M}_{H}^{H}\simeq \mathcal{YD}_{H}^{H}.%
}$ Since $t$ is in particular a total integral, the antipode $S$ is
bijective and hence, by \cite[Corollary 6.4]{Schauenburg2}, we can also
assume ${{^{H}\mathcal{YD}_{H}}}\simeq {_{H}^{H}\mathfrak{M}_{H}^{H}\simeq {%
_{H}\mathcal{YD}}}^{H}.$ Now, by means of Lemma \ref{lem ad-coinvariant}.
one prove the above equivalences.
\end{remark}

\begin{theorem}
\label{teo fd ad-coinv}Let $H$ be a finite dimensional Hopf
algebra over a field $K$ and let $D(H)$ be the Drinfeld Double.
The following assertions are equivalent:

$(i)$ There is an $ad$-coinvariant integral $t\in H.$

$(ii)$ The forgetful functor ${\mathfrak{M}}^{D(H)^{\ast
}}\rightarrow {{\mathfrak{M}^{H}}}$ (equiv.
$_{D(H)}{\mathfrak{M}}\rightarrow {_{H^*}{\mathfrak{M}}}$) is
separable.

$(iii)$ $D(H)^{\ast }$ is coseparable in
$({^{H}\mathfrak{M}^{H},\square }_{H},H)$ (equiv. $D(H)/H^*$ is
separable).
\end{theorem}

\begin{proof}
It is dual to Theorem \ref{teo fd ad-inv}.
\end{proof}

\begin{proposition}
\label{teo Fcosep}Let $H$ be a Hopf algebra with an
$ad$-coinvariant integral $t$ and let $\mathcal{M}$ be either
${\mathfrak{M}_{H}}$ or ${_{H}\mathfrak{M}_{H}}$. For any
coalgebra $C$ in $\mathcal{M}$, we have:

i) $C$ is coseparable as a coalgebra in $\mathcal{M}$ if and only
if it is coseparable as a coalgebra in ${\mathfrak{M}}_{K}$.

ii) $C$ is formally smooth as a coalgebra in $\mathcal{M}$ if and
only if $it$ is formally smooth in ${\mathfrak{M}}_{K}$.
\end{proposition}

\begin{proof}
Since $H$ has an $ad$-coinvariant integral $t$, by Lemma \ref{lem
ad-coinv},
the forgetful functor $G:{^{C}\mathcal{M}^{C}}\rightarrow {^{C}\mathfrak{M}%
^{C}}$ is separable. By Theorem \ref{teo coseparability of F} we
conclude.
\end{proof}

\section{Splitting coalgebra homomorphisms}\label{sec: Splitting coalgebra
homomorphisms}

\begin{claim}
Let $E$ be a coalgebra in an abelian monoidal category
$\mathcal{M}$. Let us recall,
(see \cite[\S5.2]{Montgomery}), the definition of wedge of two subobject $%
X,Y $ of $E$ in $\M:$%
\begin{equation*}
X\wedge_E Y:=Ker[ (\pi _{X}\otimes \pi _{Y}) \circ \Delta _{E}] ,
\end{equation*}
where $\pi _{X}:E\rightarrow E/X$ and $\pi _{Y}:E\rightarrow E/Y$ are the
canonical quotient maps.
\end{claim}

\begin{claim}Let now $C$ be a subcoalgebra of $E$ in an abelian monoidal category $\mathcal{M}$. Define $(C^{\wedge_E^n})_{n\in \mathbb{N}}$ by
\begin{equation*}C^{\wedge_E^0}:=0,\text{\qquad} C^{\wedge_E^1}:=C,\text{\qquad and \qquad
}C^{\wedge_E^n}:=C^{\wedge_E^{n-1}}\wedge_E C\text{ for any
}n\geq2.
\end{equation*}
Note that $C^{\wedge_E^1}\hookrightarrow \cdots \hookrightarrow
C^{\wedge_E^n}\hookrightarrow C^{\wedge_E^{n+1}}
\hookrightarrow \cdots \hookrightarrow E$ as coalgebras.\\
In the case when $\mathcal{M}$ is one of the monoidal categories ${\mathfrak{M}%
_{K}},{\mathfrak{M}_{H}}$ or ${_{H}\mathfrak{M}_{H}}$, then the
wedge product has the following properties:

\begin{itemize}
\item  $X\wedge_E Y=\Delta ^{-1}(E\otimes Y+X\otimes E)$;

\item  $(X\wedge_E Y) \wedge_E Z=X\wedge_E (Y\wedge_E Z)$;

\item  $X\wedge_E Y$ is a subcoalgebra of $E$ whenever both $X$
and $Y$ are subcoalgebras of $E$.
\end{itemize}
\end{claim}

\begin{remark}\label{rem: Sweedler}
Let $\mathcal{M}$ be one of the monoidal categories ${\mathfrak{M}%
_{K}},{\mathfrak{M}_{H}}$ or ${_{H}\mathfrak{M}_{H}}$. Let $C$ be
a subcoalgebra of  a coalgebra $E$ in $\mathcal{M}$. Then
$C^{\wedge_E^1}\subseteq \cdots \subseteq C^{\wedge_E^n}\subseteq
C^{\wedge_E^{n+1}} \subseteq \cdots \subseteq E.$\newline
Moreover, by \cite[Remark and Proposition, page 226]{Sw}, one has
that $\cup _{n\in \mathbb{N}}C^{\wedge_E^n}=E$ if and only if
$Corad(E)\subseteq C.$ Note that $\cup _{n\in
\mathbb{N}}C^{\wedge_E^n}=\underrightarrow{\lim }C^{\wedge^i_E}.$
\end{remark}

We recall the following important result.

\begin{theorem}
\label{X co-formally smooth teo}(see \cite[Theorem 4.22] {AMS} )
Let $(C,\Delta ,\varepsilon )$ be a coalgebra in an abelian
monoidal category $(\mathcal{M},\otimes , \mathbf{1})$ with direct
limits. Then the following assertions are equivalent:

(a) $C $ is formally smooth as a coalgebra in ${\mathcal{M}}$.

(b) Let $\sigma:D\rightarrow E$ be a coalgebra homomorphism in
$\M$ which is a monomorphism in $\M$. Assume that
$E=\underrightarrow{\lim }C^{\wedge^i_E}$. If for every $r\in \N$
the canonical injection $i_r:C^{\wedge^r_E}\rightarrow
C^{\wedge^{r+1}_E}$ cosplits in $\M$, then $\sigma $ has a
retraction which is a coalgebra homomorphism in $\mathcal{M}$.
\end{theorem}

Then the previous theorem has the following application.

\begin{theorem}
\label{general co fs} Let $H$ be a Hopf algebra. Let $\mathcal{M}$ be one of the monoidal categories ${\mathfrak{M}%
_{K}},{\mathfrak{M}_{H}}$ or ${_{H}\mathfrak{M}_{H}}$. Let $C$ be
a subcoalgebra of a coalgebra $E$ in $\mathcal{M}$. Assume that
$C$ is formally smooth as a coalgebra in ${\mathcal{M}}$ and that
$Corad(E)\subseteq C$. If any inclusion map
$i_r:C^{\wedge^r_E}\rightarrow C^{\wedge^{r+1}_E}$ cosplits in
$\mathcal{M}$, then there exists a coalgebra homomorphism $\pi
:E\rightarrow C$ in ${\mathcal{M}}$ such that $\pi _{\mid
C}=\mathrm{Id}_{C}$.
\end{theorem}

\begin{proof}As observed in Remark \ref{rem: Sweedler}, we have $E=\cup _{n\in
\mathbb{N}}C^{\wedge_E^n}=\underrightarrow{\lim }C^{\wedge^i_E}$.
The conclusion follows by applying Theorem \ref{X co-formally
smooth teo}.
\end{proof}

\begin{proposition}
\label{lem cocanonical map}Let $H$ be a Hopf algebra. Let $C$ be a
subcoalgebra of a coalgebra $E$ in $\mathfrak{M}_{H}$. Assume that
$C$ is formally smooth as a coalgebra in $\mathfrak{M}_{H}$ and
that $Corad(E)\subseteq C$. Given a coalgebra homomorphism
$g:C\rightarrow H$ in $\mathfrak{M}_{H}$, then there exists a
coalgebra homomorphism $\pi :E\rightarrow C$ in $\mathfrak{M}_{H}$
such that $\pi _{\mid C}=\mathrm{Id}_{C}$.
\end{proposition}

\begin{proof} It is similar to \cite[Theorem
2.17]{AMS-Spliting}, where $C=H=Corad(E)$ is cosemisimple and
$g=\Id_H$. In order to apply Theorem \ref{general co fs}, we have
only to prove that any inclusion map
$C^{\wedge_E^n}\hookrightarrow
C^{\wedge_E^{n+1}}$ cosplits in $\mathfrak{M}_{H}$. Since $C^{\wedge_E^{n+1}}=C^{\wedge_E^{n}}%
\wedge_E C=C\wedge_E C^{\wedge_E^{n}}=\Delta _{E}^{-1}(E\otimes C+
C^{\wedge_E^{n}}\otimes E),$ the quotient
$C^{\wedge_E^{n+1}}/C^{\wedge_E^{n}}$ becomes a right
$C$-comodule in $\mathfrak{M}_{H}$ via the map $\rho _{n}^{C},$ given by $%
x+C^{\wedge_E^{n}}\mapsto (x_{1}+C^{\wedge_E^{n}})\otimes x_{2}$.
Since $g:C\rightarrow H$
is a morphism of coalgebras in ${\mathfrak{M}}_{H}$, then $(\mathrm{Id}%
\otimes g)\circ \rho _{n}^{C}$ is a right $H$-comodule structure map for $%
C^{\wedge_E^{n+1}}/C^{\wedge_E^{n}}$ that is right $H$-linear.
Thus $C^{\wedge_E^{n+1}}/C^{\wedge_E^{n}}$ becomes an object in
${\mathfrak{M}}_{H}^{H}$: by the fundamental theorem for Hopf
modules (${\mathfrak{M}_{H}^{H}\simeq {_{K}\mathfrak{M}}}$), we
get that $C^{\wedge_E^{n+1}}/C^{\wedge_E^{n}}\simeq V\otimes H$ in
${\mathfrak{M}_{H}^{H},}$ for a suitable $V\in
{{_{K}\mathfrak{M}}}$, i.e. $C^{\wedge_E^{n+1}}/C^{\wedge_E^{n}}$
is a free right $H$-module. In particular
$C^{\wedge_E^{n+1}}/C^{\wedge_E^{n}}$ is a
projective right $H$-module, so that the inclusion map $i:C^{\wedge_E^{n}}%
\hookrightarrow C^{\wedge_E^{n+1}}$ has a retraction in
${\mathfrak{M}}_{{H}}$.
\end{proof}

\begin{theorem}
\label{te:retraction}Let $H$ be a Hopf algebra and let
$\mathcal{M}$ be either ${\mathfrak{M}_{H}}$ or
${_{H}\mathfrak{M}_{H}}$. Assume that $H$ is a subcoalgebra of a
coalgebra $E$ in $\mathcal{M}$, that $H$ is formally smooth as a
coalgebra in $\mathcal{M}$ and that $Corad(E)\subseteq H$. Then
there exists a coalgebra homomorphism $\pi :E\rightarrow H$ in
$\mathcal{M}$ such that $\pi _{\mid H}=\mathrm{Id}_{H}$ for

$a)$ $\mathcal{M}={{\mathfrak{M}}}_{H}$.

$b)$ $\mathcal{M}={{_{H}\mathfrak{M}_{H}}}$ if any inclusion map $%
H^{\wedge_E^{n}}\hookrightarrow H^{\wedge_E^{n+1}}$ cosplits in
$\mathcal{M}$.
\end{theorem}

\begin{proof}
$H^{\wedge_E^{n}}$ is a subcoalgebra of $E$ in $\mathcal{M}$ and the inclusion map $%
H^{\wedge_E^{n}}\hookrightarrow H^{\wedge_E^{n+1}}$ is obviously a
morphism in $\mathcal{M} $.\newline $a)$\ Apply Proposition
\ref{lem cocanonical map} in the case when $C:=H$ and
$g:=\mathrm{Id}_{H}.$
\newline $b)$ Apply Theorem \ref{general co fs} in the case when
$C=H$.
\end{proof}

\begin{examples}
Let $E$ be a coalgebra in the category of vector spaces. Let
$C=Corad(E)$. In this case, the sequence $(C^{\wedge_E^{n}})_{n\in
\mathbb{N}}$ is simply denoted by $(E_{n})_{n\in \mathbb{N}}$ and
it is the so-called coradical filtration of $E $. \\
Let $H$ be a Hopf algebra and let $\mathcal{M}$ be either ${{\mathfrak{M}}}_{H}$ or ${{_{H}%
\mathfrak{M}_{H}}}$. Assume that $E$ is a coalgebra in
$\mathcal{M}$ and that $H=C=Corad(E) $. We have two cases.\newline
$\mathcal{M}={{_{H}\mathfrak{M}_{H}}}$) If any inclusion $%
E_{n}\hookrightarrow E_{n+1}$ cosplits in
${{_{H}\mathfrak{M}_{H}}}$ and $H$ is formally smooth as a
coalgebra in ${{_{H}\mathfrak{M}_{H},}}$ then, by Theorem
\ref{te:retraction}, there is an homomorphisms of coalgebras $\pi
:E\rightarrow H$ in ${{_{H}\mathfrak{M}_{H}}}$ such that $\pi _{\mid H}=%
\mathrm{Id}_{H}.$\newline $\mathcal{M}={{\mathfrak{M}}}_{H}$) By
\cite[Theorem 2.11] {AMS-Spliting}, since $H$ is cosemisimple in
$\mathfrak{M}_K$, then $H$ is coseparable in
${{\mathfrak{M}}}_{H}.$ In particular $H$ is formally smooth as a
coalgebra in ${{\mathfrak{M}}}_{H}.$ Again, by Theorem \ref
{te:retraction}, there is an homomorphisms of coalgebras $\pi
:E\rightarrow H $ in ${{\mathfrak{M}_{H}}}$ such that $\pi _{\mid
H}=\mathrm{Id}_{H}$ (see also \cite[Theorem 2.17]{AMS-Spliting}).
\end{examples}

Proposition \ref{lem cocanonical map} studies the existence in
$\mathfrak{M}_{H}$ of coalgebra retractions of coalgebras
inclusion $C\hookrightarrow E$ where $C$ is a formally smooth
coalgebras in $\mathfrak{M}_{H}$ endowed with a morphism of
coalgebras $g:C\rightarrow H$ in $\mathfrak{M}_{H}$. The following
results show that the existence of $ad$-coinvariant integrals
provides such a section both in ${\mathfrak{M}_{H}}$ and in
${_{H}\mathfrak{M}_{H}}$ (without $g$).

\begin{lemma}
\label{coro mono cosplit}Let $H$ be a Hopf algebra with a total
integral $t\in H.$ Let $\mathcal{M}$ be either
${\mathfrak{M}_{H}}$ or ${_{H}\mathfrak{M}_{H}}$. Then any
monomorphism in $\mathcal{M}$ has a retraction in ${\mathcal{M}}$.
\end{lemma}

\begin{proof}
Since $t$ is a total integral in $H$, then $H$ is separable by Theorem \ref
{teo integrals}-$2)$. Therefore any right (resp. left, two-sided) $H$-module
is injective (see Corollary \ref{coro proj}). In particular any monomorphism
in ${\mathcal{M}}$ has a retraction in ${\mathcal{M}}$.
\end{proof}

\begin{theorem}
Let $H$ be a Hopf algebra with an $ad$-coinvariant integral $t\in
H$. Let $\mathcal{M}$ be either ${\mathfrak{M}_{H}}$ or
${_{H}\mathfrak{M}_{H}}$. Let $C$ be a subcoalgebra of a coalgebra
$E$ in $\mathcal{M}$. Assume that $C$ is formally smooth as a
coalgebra in ${\mathfrak{M}}_{K}$ and that $Corad(E)\subseteq C$.
Then there exists a coalgebra homomorphism $\pi :E\rightarrow C$
in ${\mathcal{M}}$ such that $\pi _{\mid C}=\mathrm{Id}_{C}$.
\end{theorem}

\begin{proof}
By Proposition \ref{teo Fcosep}, $C$ is formally smooth as a coalgebra in ${%
\mathcal{M}}$. Since $t$ is in particular a total integral in $H,$
by Lemma \ref{coro mono cosplit}, any monomorphism in
${\mathcal{M}}$, in particular the inclusion map
$C^{\wedge_E^{n}}\hookrightarrow C^{\wedge_E^{n+1}}$for any $n\in
\mathbb{N,}$ has a retraction in ${\mathcal{M}}$. Now apply
Theorem \ref {general co fs}.
\end{proof}

\begin{theorem}\label{coro co5.23}
Let $H$ be a Hopf algebra with an $ad$-coinvariant integral and
such that $H$ is formally smooth as a coalgebra in
$\mathfrak{M}_{K}$. Let $\mathcal{M}$ be either
${\mathfrak{M}_{H}}$ or ${_{H}\mathfrak{M}_{H}}$. If $H$ is a
subcoalgebra of a coalgebra $E$ in $\mathcal{M}$ and
$Corad(E)\subseteq H$, then there exists a coalgebra homomorphism
$\pi :E\rightarrow H$ in $\mathcal{M}$ such that $\pi _{\mid
H}=\mathrm{Id}_{H}$.
\end{theorem}

\begin{remark}
By Proposition \ref{teo Fcosep}, if $H$ is a Hopf algebra with an
$ad$-coinvariant integral and $H$ is formally smooth as a
coalgebra in $({\mathfrak{M}_{K}},\otimes ,K),$ then it is
formally smooth as a coalgebra in $({{\mathfrak{M}}}_{H}{,\otimes
,K}) .$ Then the case ${\mathcal{M}}={{\mathfrak{M}}}_{H} $ of the
above corollary can be also deduced by Theorem
\ref{te:retraction}.
\end{remark}

\section{Formal Smoothness of a Hopf algebra as a
coalgebra}\label{sec: Formal Smoothness of a Hopf algebra as a
coalgebra}

In order to apply Theorem \ref{te:retraction}, it is useful to
characterize when the coalgebra $H$ is formally smooth either in ${{\mathfrak{M}}}_{H}$ or in ${{_{H}%
\mathfrak{M}_{H}}}$.

\begin{claim}
Let $H$ be a Hopf algebra with antipode $S$ over a field $K.$ We denote by $%
\overline{H}$ the cokernel of the unit $u:K\rightarrow H.$\newline
Observe that $u$ can be regarded as a morphism in
${_{H}^{H}\mathcal{YD}}$, once $H$ is regarded as an object in
${_{H}^{H}\mathcal{YD}}$ via the action $\vartriangleright $
(defined in \ref{adj actions}) and the coaction given by the
comultiplication ${\Delta }$. In this way $\overline{H}=\C(u)$
inherits the following structure of left-left Yetter-Drinfeld
module:
\begin{equation*}
h\cdot \overline{x}=\overline{h_{1}xS(h_{2})},\text{\qquad }^{H}\rho (%
\overline{x})=x_{1}\otimes \overline{x_{2}}
\end{equation*}
for all $h\in H$ and $x\in H$ (by $\overline{x}$ we denote the
image of $x$ in $\overline{H}$). \newline We call an\emph{\
}$\emph{fs}$\emph{-retraction} any map $\chi :H\otimes
\overline{H}\rightarrow \overline{H}$ such that:

$(i)$ $a_{1}\otimes \overline{a_{2}}=x_1\otimes \chi (x_{2}\otimes
\overline{y}),$

$(ii)$ $\chi (x_{1}\otimes
\overline{x_{2}})=\overline{x},$\newline for all $x,y\in H,$ where
$\chi (x\otimes \overline{y})=\overline{a}$.

We say that an $fs$-retraction is \emph{complete} whenever

$(iii)$ $\chi [ h_{1}xS(h_{4})\otimes \overline{h_{2}yS(h_{3})}] =\overline{%
h_{1}aS(h_{2})},$\newline for all $h,x,y\in H,$ where $\chi
(x\otimes \overline{y})=\overline{a}$.

\end{claim}

\begin{lemma}
\label{lem fs-retraction} $\chi$ is a complete $fs$-retraction if and only $%
\chi $ is a retraction in ${_{H}^{H}\mathcal{YD}}$ of the unit $\eta _{%
\overline{H}}={^{H}\rho _{\overline{H}}}:\overline{H}\to H\otimes \overline{H%
}$ of the adjunction $(F_3,G_3)$ introduced in \ref{ex. of adjoint}
\end{lemma}

\begin{proof}
The notion of complete $fs$-retraction can be read as follows:
condition
$(i)$ means that $\chi $ is left $H$%
-colinear, $(iii)$ that $\chi $ is left $H$-linear and $(ii)$ that
$\chi $ is a retraction of the unit
$\eta _{\overline{H}}$ of the adjunction $%
(F_3,G_3)$, i.e. $\chi \circ \eta
_{\overline{H}}=\mathrm{Id}_{\overline{H}}$.
\end{proof}

\begin{proposition}
\label{pro: chi and fd} Let $H$ be a finite dimensional Hopf
algebra over a
field $K$ and let $\overline{H}$ be the cokernel of the unit $%
u_{H}:K\rightarrow H$. Let $\chi :H\otimes \overline{H}\rightarrow \overline{%
H}$ be a $K$-linear map such that $$a_{1}\otimes \overline{a_{2}}
=x_{1}\otimes \chi \left( x_{2}\otimes \overline{y}\right) $$ for
all $x,y\in
H,$ where $\chi \left( x\otimes \overline{y}\right) =\overline{a}.$ Then $%
\chi :H\otimes \overline{H}\rightarrow \overline{H}$ quotients to a map $%
\overline{\chi }:\overline{H}\otimes \overline{H}\rightarrow
\overline{H}$.
\end{proposition}


\begin{proposition}
\label{formal in bimodules}Let $H$ be a Hopf algebra with antipode $S$ over
a field $K$ and let $\overline{H}$ be the cokernel of the unit $%
u:K\rightarrow H$. The following assertions are equivalent:

$(a)$ $H$ is formally smooth as a coalgebra in ${_{H}\mathfrak{M}_{H}.}$

$(b)$ $\overline{H}$ is $\mathcal{I}_{F}$-injective where $F$ is the
forgetful functor ${_{H}^{H}\mathcal{YD}}\rightarrow {_{H}\mathfrak{M}}$.

$(c)$ There exists a complete $fs$-retraction $\chi :H\otimes \overline{H}%
\rightarrow \overline{H}.$

Moreover, if $H$ finite dimensional, the following assertion is
also equivalent to the others:

$(d) $ The comultiplication $\overline{H}\rightarrow
\overline{H}\otimes \overline{H}$ has a left $H$-colinear
retraction $\overline{\chi} :\overline{H}\otimes
\overline{H}\rightarrow \overline{H}$, where $\overline{H}\otimes
\overline{H}$ is a left $H$-comodule via $^{H}\rho
_{\overline{H}}\otimes \overline{H}$ and such that
$\overline{\chi} [ \overline{h_{1}xS(h_{4})}\otimes \overline{h_{2}yS(h_{3})}] =\overline{%
h_{1}aS(h_{2})},$ for every $h,x,y\in H$, where
$\overline{\chi}(\overline{x}\otimes \overline{y})=\overline{a}$.
\end{proposition}

\begin{proof}Analogous to Proposition \ref{formal in
bicomodules}.
\end{proof}

The referee pointed out that the equivalence $(c)\Leftrightarrow
(e)$ in the following proposition was also proved in \cite[Theorem
1.2]{Mas}.

\begin{proposition}\label{formal in modules}
Let $H$ be a Hopf algebra over a field $K$ and let $%
\overline{H}$ be the cokernel of the unit $u:K\rightarrow H$. The following
assertions are equivalent:

$(a)$ $H$ is formally smooth as a coalgebra in
${\mathfrak{M}_{K}}$.

$(b)$ $H$ is formally smooth as a coalgebra in
${\mathfrak{M}_{H}}$.

$(c)$ $\overline{H}$ is injective in $^{H}\mathfrak{M}$.

$(d)$ There exists an $fs$-retraction $\chi :H\otimes
\overline{H}\rightarrow \overline{H}$.

$(e)$ $H$ is a hereditary $K$-coalgebra.

Moreover, if $H$ finite dimensional, the following assertion is
also equivalent to the others:

$(f) $ The comultiplication $\overline{H}\rightarrow
\overline{H}\otimes \overline{H}$ has a left $H$-colinear
retraction, where $\overline{H}\otimes \overline{H}$ is a left
$H$-comodule via $^{H}\rho _{\overline{H}}\otimes \overline{H}$.
\end{proposition}

\begin{remark}
Let $H$ be a Hopf algebra over a field $K$. Then $H$ is formally
smooth as an coalgebra in ${_H\mathfrak{M}_{H}}\Rightarrow$ $H$ is
formally smooth as an algebra in ${\mathfrak{M}_{H}}$.
\end{remark}

\begin{corollary}\label{coro: ad-coinv + fs}
Let $H$ be a Hopf algebra over a field $K$. Assume that $H$ has an $ad$%
-coinvariant integral. Let $\overline {H}$ be the cokernel of the
unit $u_H:K\to H$. The following assertions are equivalent:

$(i)$ ${H}$ is formally smooth as a coalgebra in
${\mathfrak{M}}_{K}$.

$(ii)$ ${H}$ is formally smooth as a coalgebra in
${\mathfrak{M}_{H}}$.

$(iii)$ ${H}$ is formally smooth as a coalgebra in
${_{H}\mathfrak{M}_{H}}$.

$(iv)$ $\overline {H}$ is injective in $^{H}\mathfrak{M}$.

$(v)$ There exists an $fs$-retraction $\chi :H\otimes
\overline{H}\rightarrow \overline{H}$.

$(vi)$ $H$ is a hereditary $K$-coalgebra.

$(vii) $ The comultiplication $\overline{H}\rightarrow
\overline{H}\otimes \overline{H}$ has a left $H$-colinear
retraction, where $\overline{H}\otimes \overline{H}$ is a left
$H$-comodule via $^{H}\rho _{\overline{H}}\otimes \overline{H}$.
\end{corollary}

\begin{proof}
  It is analogous to Corollary \ref{coro: ad-inv + fs}. Note that here $H$
  is always finite dimensional since we have an ad-coinvariant (in
  particular total) integral in $H$.
\end{proof}

\begin{theorem}\label{teo: K alla G is fs}Let $G$ be a finite group an let $K^G$ be the
Hopf algebra of functions from $G$ to $K$. Then the following
assertions are equivalent:

$(i) $   $K^G$ is formally smooth as a coalgebra in
${\mathfrak{M}}_{K}$.

$(ii) $   $K^G$ is formally smooth as a coalgebra in
${\mathfrak{M}_{K^G}}$.

$(iii) $   $K^G$ is formally smooth as a coalgebra in
${_{K^G}\mathfrak{M}_{K^G}}$.

$(iv) $ $\overline{K^G}$ is injective in $^{K^G}\mathfrak{M}$.

$(v)$ There exists an $fs$-retraction $\chi :K^G\otimes
\overline{K^G}\rightarrow \overline{K^G}$.

$(vi)$ $K^G$ is a hereditary $K$-coalgebra.

$(vii) $ The comultiplication $\overline{K^G}\rightarrow
\overline{K^G}\otimes \overline{K^G}$ has a left $K^G$-colinear
retraction, where $\overline{K^G}\otimes \overline{K^G}$ is a left
$K^G$-comodule via $^{K^G}\rho _{\overline{K^G}}\otimes
\overline{K^G}$.
\end{theorem}

\begin{proof}
By Example \ref{ex ad coinv of K alla G}, the Hopf algebra ${K}^G$
admits an $ad$-coinvariant integral. The conclusion follows by
Corollary \ref{coro: ad-coinv + fs}.
\end{proof}

\begin{remark}
Let $G$ be a finite group. In this case both $KG$ and $K^G$ are
finite dimensional. As observed in Example \ref{ex ad coinv of K
alla G}, $K^G$ becomes a Hopf algebra which is dual to the group
algebra $KG$. In particular, $K^G$ is formally smooth as a
coalgebra in ${\mathfrak{M}}_{K}$ if and only if $KG$ is formally
smooth as an algebra in ${\mathfrak{M}}_{K}$. Hence all the
assertions in Theorem \ref{teo: KG is fs} and in Theorem \ref{teo:
K alla G is fs} are equivalent. In the particular case when $G$ is
$C_n$, the cyclic group of order $n$, then, by Proposition
\ref{pro: cyclic group} $K^G$ is formally smooth as a coalgebra in
${\mathfrak{M}}_{K}$ if and only if $n\cdot 1_K\neq 0.$

\end{remark}

\begin{proposition}\label{pro: poly fs}
Let $K\left[ X\right] $ be the polynomial ring endowed with the
unique Hopf
algebra structure defined by%
\begin{equation*}
\Delta \left( X\right) =1\otimes X+X\otimes 1.
\end{equation*}%
Then $K\left[ X\right] $ is formally smooth as a coalgebra in $\mathfrak{M}%
_{K}$ if and only if $\mathrm{char}\left( K\right) =0.$
\end{proposition}

\begin{proof}
Let $A=K\left[ X\right] $. Assume that $A$ is formally smooth as a
coalgebra
in $\mathfrak{M}_{K}$. Note that $\overline{A}=\sum_{n>0}KX^{n}.$We have%
\begin{equation}
\Delta \left( X^{a}\right) =\left( 1\otimes X+X\otimes 1\right)
^{a}=\sum_{0\leq i\leq a}\binom{a}{i}\left( 1\otimes X\right)
^{a-i}\left( X\otimes 1\right) ^{i}=\sum_{0\leq i\leq
a}\binom{a}{i}X^{i}\otimes X^{a-i}. \label{form: delta poly}
\end{equation}%
By Proposition \ref{formal in modules}, there exists a
$fs$-retraction $\chi :A\otimes \overline{A}\rightarrow
\overline{A}$. For every $a,b\in
\mathbb{N}
,$ $\chi \left( X^{a}\otimes \overline{X^{b}}\right) \in \overline{A}%
=\sum_{n>0}KX^{n}$ so that we can choose $\alpha _{u}^{a,b}\in K$
such that
\begin{equation*}
\chi \left( X^{a}\otimes \overline{X^{b}}\right) =\sum_{u\geq
1}\alpha _{u}^{a,b}X^{u},
\end{equation*}%
where $\alpha _{u}^{a,b}=0,$ for every $u\geq \deg \left( \chi
\left( X^{a}\otimes \overline{X^{b}}\right) \right) .$ By
condition (i) of the definition of $fs$-retraction, we have
\begin{equation*}
a_{1}\otimes \overline{a_{2}}=x_{1}\otimes \chi (x_{2}\otimes
\overline{y})
\end{equation*}%
for every $x,y\in H,$ where $a\in H$ is defined by
$\overline{a}=\chi \left(
x\otimes \overline{y}\right) .$ We apply this, for every $b>0,$ to the case%
\begin{equation*}
x\otimes \overline{y}=X\otimes X^{b},\qquad \overline{a}=\chi
\left( X\otimes X^{b}\right) =\sum_{u\geq 1}\alpha
_{u}^{1,b}X^{u}.
\end{equation*}%
Since%
\begin{equation*}
\left\{
\begin{tabular}{l}
$a_{1}\otimes \overline{a_{2}}\overset{(\ref{form: delta poly})}{=}%
\sum_{u\geq 1}\alpha _{u}^{1,b}\sum_{0\leq i\leq
u}\binom{u}{i}X^{i}\otimes
\overline{X^{u-i}}=\sum_{u\geq 1}\alpha _{u}^{1,b}\sum_{0\leq i\leq u-1}%
\binom{u}{i}X^{i}\otimes X^{u-i}$ ,\\
$x_{1}\otimes \chi (x_{2}\otimes \overline{y})=1\otimes \chi
\left( X\otimes
X^{b}\right) +X\otimes \chi \left( 1\otimes X^{b}\right), $%
\end{tabular}%
\right.
\end{equation*}%
we get%
\begin{equation*}
\sum_{u\geq 1}\alpha _{u}^{1,b}\sum_{0\leq i\leq u-1}\binom{u}{i}%
X^{i}\otimes X^{u-i}=1\otimes \chi \left( X\otimes X^{b}\right)
+X\otimes \chi \left( 1\otimes X^{b}\right) .
\end{equation*}%
Therefore%
\begin{equation*}
\left\{
\begin{tabular}{l}
$\sum_{u\geq 2}\alpha _{u}^{1,b}uX^{1}\otimes X^{u-1}=X\otimes
\chi \left(
1\otimes X^{b}\right) ,$ \\
$\sum_{u\geq 3}\alpha _{u}^{1,b}\sum_{2\leq i\leq u-1}\binom{u}{i}%
X^{i}\otimes X^{u-i}=0,$%
\end{tabular}%
\right.
\end{equation*}%
so that
\begin{equation*}
\left\{
\begin{tabular}{l}
$\sum_{u\geq 2}\alpha _{u}^{1,b}uX^{u-1}=\chi \left( 1\otimes
X^{b}\right) ,$
\\
$\alpha _{u}^{1,b}\binom{u}{i}=0,$ for every $u\geq 3$ and $2\leq i\leq u-1.$%
\end{tabular}%
\right.
\end{equation*}%
Now, from these equalities, where the last one is applied in the case when $%
i=u-1$, we deduce%
\begin{equation*}
\chi \left( 1\otimes X^{b}\right) =\sum_{u\geq 2}\alpha
_{u}^{1,b}uX^{u-1}=\alpha _{2}^{1,b}2X^{2-1}=2\alpha _{2}^{1,b}X.
\end{equation*}%
If $\mathrm{char}\left( K\right) \neq 0,\ $there is a prime $p$ such that $%
\mathrm{char}\left( K\right) =p.$ Since $p\mid
\binom{p}{i},\forall 1\leq i\leq p-1,$ by condition (ii) of the
definition of $fs$-retraction, we have
\begin{equation*}
X^{p}=\chi \left( X_{1}^{p}\otimes \overline{X_{2}^{p}}\right) \overset{(\ref%
{form: delta poly})}{=}\sum_{0\leq i\leq p-1}\binom{p}{i}\chi
\left( X^{i}\otimes X^{p-i}\right) =\chi \left( 1\otimes
X^{p}\right) =2\alpha _{2}^{1,p}X.
\end{equation*}%
that is a contradiction. Therefore $\mathrm{char}\left( K\right)
=0.$
\newline
Conversely, if $\mathrm{char}\left( K\right) =0.$ Consider the vector space $%
C=K\left[ X\right] $ of polynomials in one variable. $C$ can be
regarded as a Hopf algebra with the following structures
\begin{equation*}
\Delta \left( X^{a}\right) =\sum_{i+j=a}X^{i}\otimes X^{j}\qquad \text{and}%
\qquad X^{a}X^{b}=\binom{a+b}{a}X^{a+b}\text{, for every }a,b\geq
0.
\end{equation*}%
By the universal property of the polynomial ring, there exists a
unique algebra homomorphism $\varphi :A\rightarrow C$ such that
$\varphi \left( X\right) =X.$ In fact $\varphi \left( X^{n}\right)
=\varphi \left( X\right) ^{n}=n!X^{n}$, for every $n\geq 0,$ and
$\varphi $ is a Hopf algebra isomorphism (in view of the condition
on the characteristic, one can construct an inverse for $\varphi
$). We conclude by observing that $C$ is exactly the cotensor
coalgebra $T_{K}^{c}\left( K\right) $ which is always formally
smooth as a coalgebra in $\mathfrak{M}_{K}$ (see \cite{JLMS}).
\end{proof}

\begin{remark}
Akira Masuoka pointed out that the "if" part of Proposition
\ref{pro: poly fs} is the same as \cite[Example 1.8]{Mas}, where
it is proved that the polynomial ring $K\left[ X\right] $ is an
hereditary coalgebra when $\mathrm{char}\left( K\right) =0$ (see
also Proposition \ref{formal in modules}).
\end{remark}

By means of Proposition \ref{formal in modules}, it is now
possible to rewrite Theorem \ref{te:retraction} in the following
form which improves Theorem \ref{coro co5.23} in the case
$\M={\mathfrak{M}_{H}}$.

\begin{theorem}\label{teo: preweak}
Let $H$ be a Hopf algebra which is a subcoalgebra of a coalgebra $E$ in ${\mathfrak{M}}_{H}$. Assume that $H$ is formally smooth as a coalgebra in ${%
\mathfrak{M}}_{K}$ and that $Corad(E)\subseteq H$. Then there
exists a coalgebra homomorphism $\pi :E\rightarrow H$ in
${\mathfrak{M}}_{H}$ such that $\pi _{\mid H}=\mathrm{Id}_{H}$.
\end{theorem}

\begin{remark}
  The referee pointed out to our attention \cite[Theorem
1.2]{Mas}. In view of $(i)\Rightarrow(iv)$ of this result, since
any formally smooth coalgebra is also hereditary (see
\cite[Proposition 2.2]{JLMS}), one gets Theorem \ref{teo:
preweak}.
\end{remark}

\begin{definition}\cite[Definition 5.1]{Schauenburg1}
Let $E$ be a bialgebra and let $H$ be a Hopf subalgebra of $E$.
Recall that a \emph{weak projection} (onto $H$) is a retraction
$\pi :E\rightarrow H$ for the inclusion map which is a left
$H$-linear coalgebra map.
\end{definition}

\begin{claim}
Let $E$ be a bialgebra and $H$ a Hopf subalgebra. Given a weak projection $%
\pi :E\rightarrow H$ one can construct a $K$-linear isomorphism
$\psi
:E\rightarrow H\otimes R,$ where $R=E/H^{+}E$. The bialgebra structure that $%
H\otimes R$ inherits via $\psi $ has been described in \cite[Section 5]{Schauenburg1} and in \cite[Section 5]%
{Schauenburg3}.
\end{claim}

As a consequence of the left hand version of Theorem \ref{teo:
preweak}, we get the following result.

\begin{theorem}
\label{coro: fs and linear retraction}\label{teo: weak}Let $H$ be
a Hopf subalgebra of a bialgebra $E$. Assume that $H$ is formally
smooth as a coalgebra in ${\mathfrak{M}}_{K} $ and that
$Corad(E)\subseteq H$. Then $E$ has a weak projection onto $H$.
\end{theorem}

\begin{remark}
Akira Masuoka pointed out that, in the situation of Theorem
\ref{teo: weak}, by Takeuchi's lemma \cite[Lemma
5.2.10]{Montgomery}, $E$ is necessarily a Hopf algebra.
\end{remark}

\begin{proposition}\label{pro: connected}
Let $E$ be a connected Hopf algebra over a field $K$ with $\mathrm{char}%
\left( K\right) =0$. Assume that $E\neq K.$ Then, for every $x\in
P\left( E\right) \backslash \left\{ 0\right\} $, there exists a
weak projection $\pi :E\rightarrow K\left[ x\right] $. In
particular we have a $K$-linear
isomorphism.%
\begin{equation*}
E\simeq K\left[ x\right] \otimes \frac{E}{xE}.
\end{equation*}
\end{proposition}

\begin{proof}
Since $E\neq K,$ we have $P\left( E\right) \neq \left\{ 0\right\} .$ Let $%
x\in P\left( E\right) \backslash \left\{ 0\right\} .$ Note that $K\left[ X%
\right] $ is isomorphic to the tensor algebra $T_{K}\left(
KX\right) $ as a Hopf algebra, the isomorphism being given by the
assignment
\begin{equation*}
X^{n}\mapsto \underset{n}{\underbrace{X\otimes \cdots \otimes X}}.
\end{equation*}%
By the universal property of tensor algebra, there is a unique
Hopf algebra homomorphism $\sigma :K\left[ X\right] \rightarrow
E,$ such that $\sigma
\left( X\right) =x.$ Since $\mathrm{char}\left( K\right) =0,$ we have that $K%
\left[ X\right] $ is a connected coalgebra with $P\left( K\left[
X\right] \right) =KX$. As $\sigma _{\mid KX}$ is injective, by
\cite[Lemma 5.3.3, page 65]{Montgomery}, $\sigma $ is injective
and hence $\text{Im}\left( \sigma \right) \simeq K\left[ X\right]
$ as Hopf algebras. Therefore, by Proposition \ref{pro: poly fs},
$H:=\text{Im}\left( \sigma \right) $ is
formally smooth as a coalgebra in $\mathfrak{M}_{K}.$ Clearly $%
Corad(E)=K\subseteq H.$ We conclude by applying Corollary
\ref{coro: fs and linear retraction} and observing that $K\left[
x\right] ^{+}=\left( x\right) , $ the left ideal of $K\left[
x\right] $ generated by $x$.
\end{proof}

\noindent \textbf{Acknowledgements.} I would like to thank the
referee for many useful comments. My gratitude also goes to Prof.
Akira Masuoka for his helpful remarks.

\end{document}